\documentclass{article}

\usepackage{geometry}
\usepackage{mathtools}
\usepackage{graphicx}
\usepackage{amsmath}
\usepackage{amssymb}
\usepackage{dsfont}
\usepackage{amsthm}
\usepackage{cases}
\usepackage{hyperref}
\usepackage{abstract}
\usepackage{float}
\usepackage{titling}
\usepackage{extarrows}
\usepackage[dvipsnames]{xcolor}
\usepackage[ruled,vlined]{algorithm2e}
\usepackage{enumerate}
\usepackage{bm}
\usepackage{multirow}
\usepackage[shortlabels]{enumitem}
\usepackage[
	safeinputenc,
	backend=biber,
	style=alphabetic,
	maxnames = 10,
	minnames = 1,
	isbn=false,
	url=false,
	eprint=false
  ]{biblatex}
\hypersetup{
    colorlinks=True,
    linkcolor=blue,
    filecolor=magenta,      
    urlcolor=cyan,
    citecolor=ForestGreen
}
\newtheorem{Theorem}{Theorem}[section]
\newtheorem{Proposition}{Proposition}[section]
\newtheorem{Lemma}{Lemma}[section]
\newtheorem{Corollary}{Corollary}[section]

\theoremstyle{definition}

\newtheorem{Definition}{Definition}[section]
\newtheorem{Remark}{Remark}[section]
\newtheorem{Observation}{Observation}

\def\ra{\rangle}
\def\la{\langle}
\def\intOm{\int\limits_\Omega}

\begin{document}
\title{Regularity lost: the fundamental limitations and constraint qualifications in the problems of elastoplasticity}
\author{Ivan Gudoshnikov \thanks{Institute of Mathematics of the Czech Academy of Sciences, \v{Z}itn\'{a} 25, 110 00, Praha 1, Czech Republic, \url{gudoshnikov@math.cas.cz}}}
\thanksmarkseries{arabic}
\date{2024, revised in 2025}

\renewcommand{\thefootnote}{\fnsymbol{footnote}} 
\renewcommand{\thefootnote}{\arabic{footnote}} 

\renewcommand{\thefootnote}{\fnsymbol{footnote}} 
\renewcommand{\thefootnote}{\arabic{footnote}} 

\maketitle
\begin{abstract}
We investigate the existence and non-existence of a function-valued strain solution in various models of elastoplasticity from the perspective of the constraint-based ``dual'' formulations. We describe abstract frameworks for linear elasticity, elasticity-perfect plasticity and elasticity-hardening plasticity in terms of adjoint linear operators and convert them to equivalent formulations in terms of differential inclusions (the sweeping process in particular). Within such frameworks we consider several manually solvable examples of discrete and continuous models. Despite their simplicity, the examples show how for discrete models with perfect plasticity it is possible to find the evolution of stress and strain (elongation), yet continuum models within the same framework may not possess a function-valued strain. Although some examples with such phenomenon are already known, we demonstrate that it may appear due to displacement loading. The central idea of the paper is to explain the loss of strain regularity in the dual formulation by the lack of additivity of the normal cones and the failure of Slater's constraint qualification.

In contrast to perfect plasticity, models with hardening are known to be well-solvable for strains. We show that more advanced constraint qualifications can help to distinguish between those cases and, in the case of hardening, ensure the additivity of the normal cones, which means the existence of a function-valued strain rate.
\end{abstract}
\tableofcontents
\newpage
\section{Introduction}
\label{sect:intro}
The treatment of elastoplasticity as a variational inequality problem in the PDE setting goes back to the book by Duvaut and Lions \cite{Duvaut1972}. Already then it was understood that for the seemingly simplest elastoplastic constitutive law, namely for perfect plasticity, the quasistatic evolution of stress is well-solvable, but the evolution of strain is problematic. J.-J. Moreau approached this issue with the technique of differential inclusions with the normal cone operator, which worked well in a finite-dimensional setting \cite[Ch. 6]{Moreau1973}, \cite{Moreau1973-3}. Onward, Moreau could solve the strain problem for the one-dimensional continuous rod \cite{Moreau1976} by using the space of bounded measures as the state space, instead of a space of Lebesgue-integrable functions. In order to apply this idea to elastic-perfectly plastic three-dimensional continua, the theory of the space of bounded deformations ($BD$) was developed in the 1980's, with major  contributions by P.-M. Suquet \cite{Suquet1978a}, \cite{Suquet1980}, \cite{Suquet1981}, \cite{Suquet1988}, R. Temam and G. Strang \cite{Strang1980}, \cite{TemamStrangBD}, \cite{Temam1985}. More recent treatments of the problem with perfect plasticity include \cite{Ebobisse2004}, \cite{DalMaso2006}, \cite{Demyanov2009}, \cite{Francfort2016}, \cite{Mora2016}, which are also based on $BD$-spaces.


In the current text we primarily explore the necessity for the measure-valued solutions in the continuous setting, i.e. we go back to the transition from  \cite{Moreau1973} to \cite{Moreau1976}. We begin with a rigorous presentation of abstract functional-analytic frameworks for quasistatic evolution with elasticity and with Prandtl-Reuss elasticity-perfect plasticity. These frameworks are suitable for discrete as well as continuous models of media, and we provide toy examples for both. Using the ideas of  \cite{Moreau1973} and \cite{Kunze2000} we convert the abstract problem into a pair of differential inclusions of the form
\begin{numcases}{}
-\frac{d}{dt}\bm{y}\in N_{\left(\Sigma-\bm{\widetilde{\sigma}}(t)\right)\cap \mathcal{V}}(\bm{y}), \label{eq:sp-elastoplastic-intro}\\
\frac{d}{dt} \bm{\varepsilon} \in M\left(\left(N_{\Sigma-\bm{\widetilde{\sigma}}(t)}(\bm{y})+ \frac{d}{dt} \bm{y}\right)\cap \mathcal{V}^\perp\right), 
\label{eq:inclusion2-elastoplastic-intro}
\end{numcases}
where $\bm{y}=\bm{\sigma}-\bm{\widetilde{\sigma}}(t)$ is the difference between the stress in the elastoplastic problem and the known stress $\bm{\widetilde{\sigma}}$ in the corresponding elastic problem, $\bm{\varepsilon}$ is the strain in the elastoplastic problem, $N$ denotes a normal cone, $M$ is a known affine map, $\Sigma$ is the set of stress distributions which fit in the elastic range and $\mathcal{V}$ is the subspace of stress distributions which satisfy the homogeneous equilibrium equation. On this initial step we consider the classical sufficient condition of  solvability of \eqref{eq:sp-elastoplastic-intro}--\eqref{eq:inclusion2-elastoplastic-intro}, which is
\begin{equation}
\left({\rm int}\, \Sigma-\bm{\widetilde{\sigma}}(t)\right)\cap \mathcal{V} \neq \varnothing.
\label{eq:Slater-I-intro}
\end{equation}

Our main goal is to compare the outcomes for discrete (finite-dimensional) and continuous models being plugged into {\it the same} abstract framework \eqref{eq:sp-elastoplastic-intro}--\eqref{eq:Slater-I-intro}. More specifically,
\begin{itemize}
\item for finite-dimensional state space we provide easily-solvable toy examples of discrete mechanical systems and observe that the abstract framework admits solutions for both strain and stress. In such a case the sufficient condition \eqref{eq:Slater-I-intro} is adequate.
\item for the continuous models in $L^2$ we provide another simple example, for which the inclusion \eqref{eq:inclusion2-elastoplastic-intro} admits no solution with values in $L^2$. We will explain that, in contrast to the discrete case, the sufficient condition \eqref{eq:Slater-I-intro} does not hold for any of the models with plasticity in $L^2$.

\end{itemize}
 While the examples to merely observe the lack of function-valued strain solution are long present in the literature (see \cite[Sect. 3.g, pp.~69--70]{Moreau1976}, \cite[Sect. 3.4.a), pp.~312--313]{Suquet1988}), we aim to understand {\it how} Moreau's method \eqref{eq:sp-elastoplastic-intro}--\eqref{eq:inclusion2-elastoplastic-intro} fails in $L^2$ and what is the mathematical reason for that.

For the dissipation-based strain formulation of the elastoplasticity problem (``primal formulation'' in terms of \cite{Han2012}) the mathematical reason for the unsolvability is generally known: the dissipation function, defined on $L^2$, has linear growth at infinity and does not admit a minimum (see e.g. \cite[pp. vi, 79--80]{Temam1985}). We follow Moreau's constraint-based formulation  \eqref{eq:sp-elastoplastic-intro}--\eqref{eq:inclusion2-elastoplastic-intro} (``dual formulation''), and it turns out that in such a case the unsolvability is due to the {\it lack of  additivity of the normal cones}, i.e. in \eqref{eq:inclusion2-elastoplastic-intro} we may have
\begin{equation}
N_{\Sigma-\bm{\widetilde{\sigma}}(t)}(\bm{y}) +N_{\mathcal{V}}(\bm{y}) \subsetneqq  N_{\left(\Sigma-\bm{\widetilde{\sigma}}(t)\right)\cap \mathcal{V}}(\bm{y}).
\label{eq:strict-subadditivity-intro}
\end{equation}

\noindent Moreover, as we will show by an example in Section \ref{ssect:ex21-no-strain}, a continuous model with elasticity-perfect plasticity may get stuck in the situation of \eqref{eq:strict-subadditivity-intro} for a time-interval of positive length, during which the right-hand side of \eqref{eq:inclusion2-elastoplastic-intro} remains an empty set. 

It is known, however, that similar models with the {\it hardening} type of plasticity (replacing perfect plasticity),  can be solved for both strain and stress in $L^2$, see e.g. \cite{Han2012}. We can extend the abstract framework \eqref{eq:sp-elastoplastic-intro}--\eqref{eq:inclusion2-elastoplastic-intro} to include such models. In the extended version we follow the same idea to use the additivity of the normal cones to obtain $\bm{\varepsilon}$. However, the condition of the type \eqref{eq:Slater-I-intro} fails in $L^2$ for the models with hardening as well. Thus, in Section \ref{ssect:constraint-qualifications} we turn to explore various sufficient conditions for the additivity of the normal cones, generally known as {\it constraint qualifications}. 

Our second goal is to show that
\begin{itemize}
\item constraint qualifications can be used to distinguish solvable models (discrete, continuous with hardening) among those which fit in the abstract framework.
\end{itemize}

Such constraint qualifications are usually written for a pair of {\it Fenchel-Rockefaller dual problems} to ensure the existence and equality of their minimizers ({\it strong duality}). The proof of equivalence between the additivity of the normal cones and strong duality of specially constructed dual problems is included in the text. While constraint qualifications in optimization go back to Slater in the 1950's \cite{Slater1950} and Rockafellar in the 1970's \cite{Rockafellar1974}, the topic matured and was reflected on in the context of infinite-dimensional spaces only in the late 1980's--early 1990's  \cite{Attouch1986}, \cite{Gowda1990}, \cite{Jeyakumar1992}. This is much later than the pioneering works of Moreau   \cite{Moreau1973}, \cite{Moreau1976} on elastoplasticity, whose method we can now analyze with the help of the constraint qualifications machinery. 

The paper is organized as follows. Section \ref{sect:convex_prelims} is a reference of the basic definitions and facts about convex analysis, which we will need throughout the paper. This also includes the definition and well-posedness of the general type of problem, known as the ``sweeping process'', to which  \eqref{eq:sp-elastoplastic-intro} belongs. Section \ref{sect:elasticity} contains a formal framework for the problems in linear elasticity in terms of adjoint operators between abstract spaces, together with examples of finite-dimensional (simple discrete networks of springs) and infinite-dimensional (continuous rod with lateral displacements) models. We solve the abstract problem and the examples with the method of orthogonal subspaces (cf. \cite{Gudoshnikov2025}). The components of the framework and the method, such as the fundamental subspaces, as well as the final formulas for elastic solutions are reused in Section \ref{sect:perfect-plasticity}, where we consider the problem of elasticity-perfect plasticity in the abstract form and discuss its well-posedness for discrete examples. In contrast, in Section \ref{sect:regularity-lost} we uncover how for the continuous rod example, although it is written within the same abstract framework, there exists no strain solution in $L^2$ due to the lack of additivity of normal cones \eqref{eq:strict-subadditivity-intro}. In Section \ref{sect:constraint-qualification-and-hardening} we discuss constraint qualifications, which can guarantee the additivity of normal cone in infinite-dimensional spaces. Also there we extend the abstract framework to cover plasticity with hardening. We show that the continuous elastoplastic rod satisfies some of the constraint qualifications and, therefore, possesses additivity of normal cones in $L^2$ when it has hardening with a uniform linear growth estimate. Finally, Section \ref{sect:conclusions} contains summarizing conclusions and the ideas for future work, and Appendix \ref{sect:appendix} contains some additional mathematical facts we rely on in the paper. 

Let us clarify some conventions about notation:
\begin{enumerate}[{\it i)}]
\item $I$ is a compact time-interval
\[
I=[0, T]
\]
with $T>0$ given.
\item
Because we study evolution problems in vector spaces, and, particularly, in spaces of functions, for the values in such state spaces we will use a bold font. For example, for a time-dependent solution $\sigma\in W^{1,\infty}(I, L^2(\Omega))$ we write its value at a time instant $t\in I$ as $\bm{\sigma}(t)\in L^2(\Omega)$, while numerical values for a.a. $x\in \Omega$ we write as $\sigma(t,x)$. We also use a bold font for special quantities ${\bf E},{\bf E'}, {\bf D}, {\bf C}, \bm{C}, \bm{H}$ which are operators between the state spaces and real functions to define such operators. 
\item With the tilde symbol we indicate the state variables in the problems with elasticity, e.g. $\bm{\widetilde{\sigma}}$.  The state variables without tilde correspond to the problems with elastoplasticity.
\item Some quantities play a similar role in the problems with perfect plasticity, as well as in the problems with hardening. With the hat symbol (e.g. $\widehat{\Sigma}$ vs $\Sigma$) we indicate the context of a problem with hardening. 
\end{enumerate}
\section{Preliminaries from convex analysis} 
\label{sect:convex_prelims}
\begin{Definition}
Let $\mathcal{H}$ be a Hilbert space, and $\mathcal{C} \subset \mathcal{H}$ be a closed, convex, nonempty set.
\begin{enumerate}[{\it i)}]
\item
Let a point $\bm{x}\in \mathcal{C}$. Then the {\it outward normal cone} to $\mathcal{C}$ at $x$ is the set
\begin{equation}
N_{\mathcal{C}}(\bm{x}) = \left\{\bm{y}\in \mathcal{H}: \text{for any }\bm{c}\in \mathcal{C}\text{ we have } \la\bm{y}, \bm{c}-\bm{x}\ra_{\mathcal{H}}\leqslant 0\right\},
\label{eq:nc-abstract-def}
\end{equation}
see Fig. \ref{fig:fig-normal-cone} for an illustration.
\begin{figure}[H]\center
\includegraphics{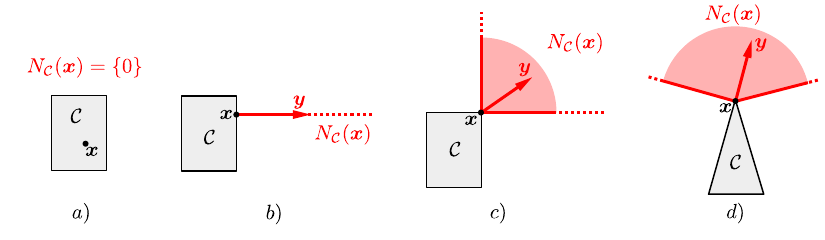}
\caption{\footnotesize Examples of typical normal cones in a finite-dimensional setting. Here we illustrate it with $\mathcal{C}$ taken as a polygon in $\mathcal{H}=\mathbb{R}^2$. We depict the normal cones $N_{\mathcal{C}}(\bm{x})$ as translated from $0$ to $\bm{x}$. Vectors $\bm{y}$ are generic elements of $N_{\mathcal{C}}(\bm{x})$, called {\it supporting vectors} to the set $\mathcal{C}$ at $\bm{x}$. 
} 
\label{fig:fig-normal-cone}
\end{figure}
\item Let a point $\bm{y}\in \mathcal{H}$. Then the {\it metric projection} of $\bm{y}$ on $\mathcal{C}$ is a point $\bm{x}\in \mathcal{C}$, denoted
\[
{\rm proj}(\bm{y}, \mathcal{C}) 
\]
such that
\[
\|\bm{x}-\bm{y}\| = \inf_{\bm{c}\in \mathcal{C}} \|\bm{c}-\bm{y}\|. 
\]
Such point always exists and it is unique \cite[Th. 5.2, p. 132]{Brezis2011}.
\item The {\it indicator function} of $\mathcal{C}$ is the function
\begin{equation}
\delta_\mathcal{C}: \mathcal{H} \to \mathbb{R}\cup\{+\infty\},\qquad 
\delta_\mathcal{C}(\bm{x})= \begin{cases} 0&\text{when }\bm{x}\in \mathcal{C},\\ +\infty & \text{when }\bm{x}\notin \mathcal{C}. \end{cases}
\label{eq:def-indicator-function}
\end{equation}
\item The {\it support function} of $\mathcal{C}$ is the function
\begin{equation}
\delta^*_{\mathcal{C}}: \mathcal{H} \to \mathbb{R}\cup\{+\infty\},\qquad
\delta^*_{\mathcal{C}}(\bm{y})= \sup_{\bm{x}\in \mathcal{C}}\,\la \bm{y}, \bm{x} \ra_{\mathcal{H}}.
\label{eq:def-support-function}
\end{equation}
\end{enumerate}
\end{Definition}
\noindent We are going to use the following well-known properties.
\begin{Proposition}
\label{prop:nc-equivalent}
Let $\mathcal{H}$ be a Hilbert space and $\mathcal{C} \subset \mathcal{H}$ be a closed, convex, nonempty set. For vectors $\bm{x}\in \mathcal{C}$ and $\bm{y}\in \mathcal{H}$ the following statements are equivalent:
\begin{enumerate}[{\it i)}]
\item $\bm{y}\in N_\mathcal{C}(\bm{x})$,
\item \label{enum:abstract-nc-in-proj-notation} ${\rm proj}(\bm{y}+\bm{x}, \mathcal{C})=\bm{x}$,
\item \label{enum:abstract-nc-in-support-func-notation} $\la\bm{y}, \bm{x}\ra_\mathcal{H} = \delta^*_\mathcal{C}(\bm{y})$.
\end{enumerate}
\end{Proposition}

\begin{Proposition}
\label{prop:nc-abstract-properties}
Formulas \eqref{eq:nc-abstract-def} and \eqref{eq:def-support-function} yield the following properties:
\begin{enumerate}[{\it i)}] 
\item \label{enum:prop:nc-abstract-properties-nc-translation} For any $\bm{x_1}, \bm{x_2}\in \mathcal{C}$
\begin{equation}
N_{\mathcal C}(\bm{x_1}+\bm{x_2}) = N_{\mathcal C-\bm{x_2}}(\bm{x_1})
\label{eq:argument-sum-nc}
\end{equation}
\item  \label{enum:prop:nc-abstract-properties-nc-subadditivity} \cite[Section 2.f, p. 206]{Moreau1973}, \cite[Lemma 1 (b), pp.~4--5]{Kunze2000} For any closed convex nonempty sets $\mathcal{C}_1,\mathcal{C}_2\subset\mathcal{H}$ and any $\bm{x}\in \mathcal{C}_1\cap\mathcal{C}_2$ the subadditivity of normal cones holds:
\begin{equation}
N_{\mathcal{C}_1}(\bm{x}) + N_{\mathcal{C}_2}(\bm{x}) \subseteqq  N_{\mathcal{C}_1 \cap \mathcal{C}_2}(\bm{x}).
\label{eq:nc-subadditivity}
\end{equation}
Additionally, if the so-called Slater's constraint qualification holds:
\begin{equation}
{\mathcal{C}_1}\cap {{\rm int}\, \mathcal{C}_2}\neq \varnothing,
\label{eq:nontrivial-intersection-condition}
\end{equation}
then we have the equality (i.e. the additivity of normal cones):
\begin{equation}
N_{\mathcal{C}_1}(\bm{x}) + N_{\mathcal{C}_2}(\bm{x}) =  N_{\mathcal{C}_1 \cap \mathcal{C}_2}(\bm{x}).
\label{eq:nc-additivity}
\end{equation}
\item  \label{enum:prop:nc-abstract-properties-support-function} For closed, convex nonempty sets $\mathcal{C}_1, \mathcal{C}_2\subset \mathcal{H}$
\begin{equation}
\mathcal{C}_1\subset \mathcal{C}_2 \qquad\Longleftrightarrow \qquad \text{for any }y\in \mathcal{H}\quad \delta^*_{\mathcal{C}_1}(y)\leqslant\delta^*_{\mathcal{C}_2}(y). 
\label{eq:sp-subset-property}
\end{equation}
\end{enumerate}
\end{Proposition}

\noindent The normal cone can be used to define the following nonsmooth evolution problem.

\begin{Definition}{\cite[Section 5.f]{Moreau1973}, \cite[p. 9]{Kunze2000}} \label{def:abstract-sp} Let $\mathcal{H}$ be a Hilbert space and $\mathcal{C}:I\rightrightarrows \mathcal{H}$ be a set-valued map such that its values $\mathcal{C}(t)$ are closed convex nonempty subsets of $\mathcal{H}$ for every $t$. Given a $\bm{y_0}\in C(0)$, we call the initial value problem 
\begin{numcases}{}
 -\frac{d}{dt} \bm{y} \in N_{\mathcal{C}(t)}(\bm{y}), \label{eq:abstract-sp}\\
\bm{y}(0)=\bm{y_0}, \label{eq:abstract-sp-ic}
\end{numcases}
a {\it (Moreau's) sweeping process}. We say that a function $y:I\to \mathcal{H}$ is a {\it solution of the sweeping process}  \eqref{eq:abstract-sp}--\eqref{eq:abstract-sp-ic} if it satisfies the following conditions:
\begin{enumerate}[{\it i)}]
\item $\bm{y}(0)=\bm{y_0}$,
\item $\bm{y}(t)\in \mathcal{C}(t)$ for all $t\geqslant 0$, \label{enum:abstract-sp-solution-def-constraint}
\item $\frac{d}{dt}\bm{y}(t)$ exists and the inclusion \eqref{eq:abstract-sp} holds for a.a. $t\in I$.
\label{enum:abstract-sp-solution-def-nc}
\end{enumerate}
\end{Definition}

The sweeping process gets its name from the intuition about the solution $\bm{y}$ being ``swept'' by the time-dependent unilateral constraint of Definition \ref{def:abstract-sp} \ref{enum:abstract-sp-solution-def-constraint}. For more details we refer to \cite[Sections 1 and  3.1]{Kunze2000} as an excellent introduction to the topic and to Moreau's original works \cite[Section 5.f]{Moreau1973}, \cite{Moreau1977}. Here we only include Definition \ref{def:abstract-sp} and the following existence and uniqueness theorem, which we will use to find the evolution of stress in problems with elastoplasticity.

\begin{Theorem}{\bf \cite[Th. 2]{Kunze2000}}
\label{th:abstract-sp-existence} \newline 
Let $\mathcal{H}$ be a Hilbert space and $\mathcal{C}:I\rightrightarrows \mathcal{H}$ be a set-valued map. Assume that
\begin{enumerate}[{\it i)}]
\item  values $\mathcal{C}(t)$ are closed convex nonempty sets for every $t\in I$,
\item map $\mathcal{C}$ is Lipschitz-continuous with a constant $L$ with respect to Hausdorff distance, i.e. there is $L\geqslant 0$ such that
\begin{equation}
d_H(\mathcal{C}(t_1),\mathcal{C}(t_2))\leqslant L|t_1-t_2| \qquad \text{for any } t_1,t_2\in I,
\label{eq:h-d-Lipschitz-def}
\end{equation}
where for any nonempty subsets $\mathcal{C}_1, \mathcal{C}_2\subset \mathcal{H}$ the Hausdorff distance between the subsets is
\begin{multline}
d_H(\mathcal{C}_1,\mathcal{C}_2)=\max\left(\sup_{x\in \mathcal{C}_2}{\rm dist}(x, \mathcal{C}_1),\,\sup_{x\in \mathcal{C}_1}{\rm dist}(x, \mathcal{C}_2)\right)=\\
=\max\left(\sup_{x\in \mathcal{C}_1}\inf_{y\in \mathcal{C}_2} \|x-y\|,\,\sup_{x\in \mathcal{C}_2}\inf_{y\in \mathcal{C}_1} \|x-y\|\right)=\\
=\inf\left\{r>0: B_r(0)+\mathcal{C}_1\subset \mathcal{C}_2 \text{ and }B_r(0)+\mathcal{C}_2\subset \mathcal{C}_1\right\},
\label{eq:h-d-def}
\end{multline}
\end{enumerate}
Then for any $\bm{y_0}\in \mathcal{C}(0)$ there is a unique solution $y$ of  \eqref{eq:abstract-sp}--\eqref{eq:abstract-sp-ic}  and $y$ is Lipschitz-continuous with the same constant $L$, i.e. $y\in W^{1,\infty}(I, \mathcal{H})$.
\end{Theorem}

\section{The abstract geometric framework and examples for linear elasticity}
\label{sect:elasticity}
\subsection{An abstract linear operator and its adjoint}
\label{ssect:adjoint-operatorsED}
We will use the concept of an adjoint to a linear operator to unite discrete and continuum models in one generalization, but we need to account for a non densely-defined operator. The following definition introduces the operators and the setting, which we will use throughout the paper.

\begin{Definition}
\label{def:abstract-adjoint}
Let $\mathcal{H}$ be a Hilbert space and $\mathcal{X}$ be a reflexive Banach space, in which $\mathcal{W}_0$ is a linear, but not necessarily closed, subspace. Let a linear operator ${\bf E}$ be defined on $\mathcal{W}_0$, i.e.
\begin{equation}
{\bf E}: D({\bf E})=\mathcal{W}_0\subset \mathcal{X} \to \mathcal{H}.
\label{eq:ae-abstract-E-def}
\end{equation}
We assume that ${\bf E}$ is a closed operator (i.e. its graph is closed in $\mathcal{X}\times \mathcal{H}$) with its image also closed.
Let 
\[
\mathcal{X}_0=\overline{\mathcal{W}_0},
\]
which is a reflexive Banach space itself \cite[Prop. 3.20, p. 70]{Brezis2011}. Subspace $\mathcal{X}_0$ may coincide with $\mathcal{W}_0$ or $\mathcal{X}$. To clarify, here and anywhere in the paper the closure $\overline{\mathcal{W}_0}$ is in the sense of the norm of $\mathcal{X}$.
Consider now the closed {\it densely defined} linear operator with closed image
\begin{align*}
&{\bf E'}: \mathcal{W}_0\subset {\mathcal{X}_0} \to \mathcal{H},\\
&{\bf E'}: \bm{u} \mapsto {\bf E}\,\bm{u},
\end{align*}
and denote by ${\bf D}$ its adjoint ${\bf E'}^*$ in the sense of the unbounded operator theory \cite[Section 2.6, pp. 43--48]{Brezis2011} :
\[
{\bf D}:D({\bf D})\subset \mathcal{H}\to \mathcal{X}_0^*.
\]
In other words, the operator ${\bf D}$ is defined over the set
\[
D({\bf D})= \left\{\bm{\tau}\in \mathcal{H}: \text{there exists } c\geqslant 0: \text{ for all } \bm{u}\in \mathcal{W}_0:|\left\la \bm{\tau},  {\bf E}\,  \bm{u} \right\ra_{\mathcal{H}}|\leqslant c\|\bm{u}\|_\mathcal{X}\right\},
\]
and its value ${\bf D}\bm{\tau}$ is a continuous linear functional, defined by the following formula for $\bm{u}\in \mathcal{W}_0$:
\[
{\bf D}\bm{\tau}:  \bm{u}\mapsto \la \bm{\tau}, {\bf E}\, \bm{u}\ra_{\mathcal{H}},
\]
and by the continuity limit for $\bm{u}\in \mathcal{X}_0  \setminus \mathcal{W}_0$.
\end{Definition}
\begin{Lemma}
\label{lemma:adjoint_operators_lemma}
Under the assumptions of Definition \ref{def:abstract-adjoint} the following properties hold:
\begin{enumerate}[{\it i)}]
\item \label{enum:adjoint_operators_lemma_prop1} The operator ${\bf D}$ is also densely defined, and ${\bf E'} = {\bf D}^*$, see \cite[Th. 3.24, p.~72]{Brezis2011}. The image of ${\bf D}$ is also closed, see \cite[Th. 2.19, p.~46]{Brezis2011}.
\item The image of ${\bf E}$ and the kernel of ${\bf D}$ are orthogonal complements in the space $\mathcal{H}$, see \cite[pp.~45--46]{Brezis2011}.
\end{enumerate}
\end{Lemma}
\begin{Remark} 
\label{rem:abstract-elasticity-nonunique-ext}
By the Hahn-Banach theorem \cite[Corollary 1.2, p. 3]{Brezis2011} it is possible to extend each functional ${\bf D}\bm{\tau}$ to take its argument from the entire space $\mathcal{X}$. When ${\bf E}$ happens to be densely defined (i.e. ${\bf E}={\bf E'}$) such extension is unique, see \cite[p. 44]{Brezis2011}. However, in some of the examples  we will have $\mathcal{X}_0$ being a proper closed subspace of $\mathcal{X}$, which  implies that multiple different extensions of ${\bf D}\bm{\tau}$ are possible. 
\end{Remark}

\subsection{The abstract geometric framework for linear elasticity}
\label{ssect:abstract-elasticity}
We start by considering the quasistatic evolution of elastic models, which can be written via the following abstract framework.
\begin{Definition}
\label{def:ae} Let spaces $\mathcal{X}, \mathcal{X}_0, \mathcal{W}_0, \mathcal{H}$ and operators ${\bf E}, \bf{D}$ be as in the previous Section \ref{ssect:adjoint-operatorsED}. Additionally, let us be given a linear operator
\begin{equation}
{\bf C}:\mathcal{H}\to \mathcal{H}
\label{eq:ae-C-def}
\end{equation}
which is
\begin{enumerate}[{\it i)}]
\item bounded, i.e. there is $C_0>0$ s.t. for all $\bm{\tau} \in \mathcal{H}$
\[
\|{\bf C} \,\bm{\tau}\| \leqslant C_0 \|\bm{\tau}\|,
\]
\item coercive, (see \cite[Section A.2.4, pp. 473--474]{Ern2004}, \cite[p. 269]{Kress2014}) i. e. there is $c_0>0$ s.t. for all $\bm{\tau} \in \mathcal{H}$
\[
\la{\bf C} \,\bm{\tau}, \bm{\tau}\ra \geqslant c_0 \|\bm{\tau}\|^2,
\]
\item symmetric, i.e. for all $\bm{\tau}_1, \bm{\tau}_2 \in \mathcal{H}$
\[
\la {\bf C}\, \bm{\tau}_1, \bm{\tau}_2 \ra = \la \bm{\tau}_1, {\bf C} \,\bm{\tau}_2 \ra.
\]
\end{enumerate}
Finally, let us be given the following functions of time, which we call, respectively, {\it Dirichlet offset} and {\it external force} (together we call them {\it loads}):
\begin{equation}
g\in W^{1, \infty}(I, \mathcal{H}), \qquad f\in W^{1, \infty}(I,{\rm Im}\, {\bf D}).
\label{eq:loads-abstract-def}
\end{equation}
We say that unknown variables
\[
\widetilde{\varepsilon}, \,\widetilde{\sigma} \in W^{1,\infty}(I, \mathcal{H})
\]
solve the {\it abstract problem of quasi-static evolution in elasticity} if they satisfy for all $t\in I$
\begin{align}
\bm{\widetilde{\varepsilon}}& \in {\rm Im}\, {\bf E}+\bm{g}(t) , \label{eq:ae-1} \tag{EL1}\\
\bm{\widetilde{\sigma}}& = {\bf C} \, \bm{\widetilde{\varepsilon}}, \label{eq:ae-2} \tag{EL2}\\
\bm{\widetilde{\sigma}}\in D(\bf D)\quad \text{ and }\quad {\bf D}\, \bm{\widetilde{\sigma}}&=\bm{f}(t). \label{eq:ae-3} \tag{EL3}
\end{align}
From now on we will refer to $\widetilde{\sigma}$ as the {\it stress solution for elasticity}, and finding it is the main goal of the current Section \ref{sect:elasticity}.
\end{Definition}

\begin{figure}[H]\center
\includegraphics{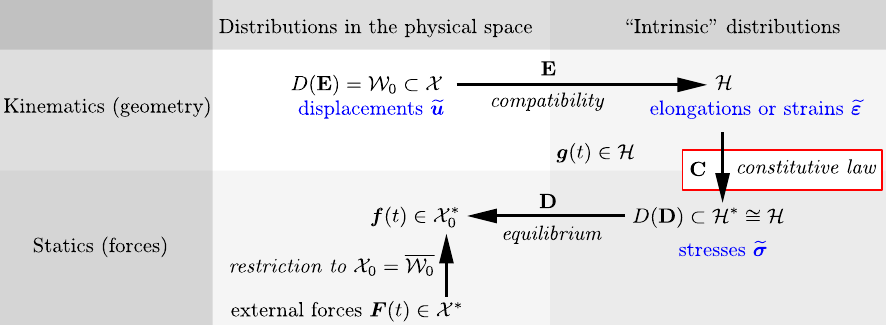}
\caption{\footnotesize Schematic representation of the problem of Definition \ref{def:ae} and Remark \ref{rem:abstract-operator-to-restrict}. The unknown variables are indicated by blue color. In the problem of Definition \ref{def:ae} we are only looking for the unknowns $\bm{\widetilde{\varepsilon}}$ and $\bm{\widetilde{\sigma}}$.
} 
\label{fig:elasticity-scheme}
\end{figure}

As we will show in the examples below, \eqref{eq:ae-1} is an abstract form of the {\it compatibility equation} and {\it boundary conditions}, \eqref{eq:ae-2} is the {\it constitutive law} of linear elasticity ({\it Hooke's law}), and \eqref{eq:ae-3} is an abstract form of the {\it equation of equilibrium}, see also the scheme of the Fig. \ref{fig:elasticity-scheme}. The notation of ${\bf E}$ and ${\bf D}$ is motivated by mechanics as well, as one can write the linear elasticity problem for a three-dimensional continuous body in the form \eqref{eq:ae-1}--\eqref{eq:ae-3}, with ${\bf E}$ being the {\it strain operator} and ${\bf D}$ being the {\it divergence operator}, see \cite{Gudoshnikov2025}. In the current text we will not consider such a problem and resort to simpler examples, but our conclusions will be meaningful for it nevertheless.

\begin{Remark}
\label{rem:abstract-operator-to-restrict}
On top of Definition \ref{def:ae} we would like to add some details, which are common in mechanical models. Although these details can be left outside of the abstract Definition   \ref{def:ae}, they will appear in the construction of the concrete examples. 
\begin{enumerate}[{\it i)}]
\item
As we will see, the operator ${\bf E}$ is usually defined as a restriction of another operator ${\rm E}$, which has its domain larger than $D({\bf E})=\mathcal{W}_0$. I.e. we begin with a given operator ${\rm E}$ acting as 
\begin{equation}
{\rm E}: D({\rm E})=\mathcal{W}\subset \mathcal{X} \to \mathcal{H},
\label{eq:ae-general-strain}
\end{equation}
where $\mathcal{W}$ is a linear space, such that $\mathcal{W}_0\subset\mathcal{W}\subset \mathcal{X}$. The subspace $\mathcal{W}$ may or may not coincide with the entire $\mathcal{X}$. Given $\mathcal{W}_0$ and such ${\rm E}$ we then define operator ${\bf E}$ acting as \eqref{eq:ae-abstract-E-def} by
\begin{equation}
{\bf E}: \bm{u} \mapsto {\rm E}\, \bm{u}\qquad \text{for any }\bm{u}\in \mathcal{W}_0.
\label{eq:ex11-E-def}
\end{equation}
\item
In turn, equation  \eqref{eq:ae-1} usually appears in the form
\begin{align}
\bm{\widetilde{\varepsilon}} &= {\rm E}\, \bm{\widetilde{u}}, \label{eq:ae-compatibility-abstract}
\\ \bm{\widetilde{u}} & \in \mathcal{W}_0 + \bm{u}_{\bf D}(t)\subset \mathcal{W}, \label{eq:ae-bc-abstract}
\end{align}

where $u_{\rm D}\in W^{1\infty}(I,\mathcal{W})$ is given, and $\bm{\widetilde{u}}\in \mathcal{W}$ is the unknown variable, which represents in mechanical models the vector distribution of {\it displacements} from a reference configuration. The purpose of equation \eqref{eq:ae-compatibility-abstract}  is to ensure the {\it geometric compatibility}, and \eqref{eq:ae-bc-abstract} is an abstract form of the {\it boundary condition}. Both equations \eqref{eq:ae-compatibility-abstract}--\eqref{eq:ae-bc-abstract} are combined in \eqref{eq:ae-1}, where 
\begin{equation}
\bm{g}(t) = {\rm E}\, \bm{u}_{\bf D}(t).
\label{eq:ae-abstract-g-def}
\end{equation}

 In this paper we focus on solving problems in terms of $\bm{\widetilde{\varepsilon}}$ and $\bm{\widetilde{\sigma}}$, and do not attempt to determine the concrete values of $\bm{\widetilde{u}}$. In general, provided with $\bm{\widetilde{\varepsilon}}$, the corresponding $\bm{\widetilde{u}}$ can be found from \eqref{eq:ae-compatibility-abstract}--\eqref{eq:ae-bc-abstract}. 
\item 
\label{rem:abstract-operator-to-restrict-about-forces}
We think of $\mathcal{X}$ as the space of {\it geometric configurations of the system in the physical space}, such as a a finite collection of displacement vectors in a system of particles or a vector field of displacements in a continuous body.  
Therefore, it is natural to think of the external force to be given for each $t\in I$ as a functional 
\[
\bm{F}(t) \in \mathcal{X}^*,
\]
i.e. as a counterpart of $\bm{\widetilde{u}}\in \mathcal{W}\subset \mathcal{X}$. Thus for any two displacement distributions $\bm{\widetilde{u}}_1, \bm{\widetilde{u}}_2\in \mathcal{W}\subset \mathcal{X}$ the expression 
\[
\la\bm{F}(t),\bm{\widetilde{u}}_2- \bm{\widetilde{u}}_1\ra_{\mathcal{X}^*, \mathcal{X}}
\]
is precisely how the {\it virtual work} of the force $\bm{F}(t)$ is defined over a path from $\bm{\widetilde{u}}_1$ to $\bm{\widetilde{u}}_2$. We should make this compatible with the operator ${\bf D}$ taking values in $\mathcal{X}_0^*$. For a given external force load $\bm{F}\in \mathcal{X}^*$ (or even $\bm{F}\in \mathcal{X}$ when $\mathcal{X}$ is a Hilbert space and the Riesz theorem is applicable), we define the functional $\bm{f}\in \mathcal{X}_0^*$ as the restriction of $\bm{F}$, i.e.
\begin{equation}
\bm{f}: \bm{u}\mapsto \la \bm{F},\bm{u} \ra_{\mathcal{X}^*, \mathcal{X}} \qquad \text{for }\bm{u} \in \mathcal{X}_0,
\label{eq:ae-force-functional-abstract}
\end{equation}
and this is the given load, which we use in the right-hand side of \eqref{eq:ae-3}. However, some caution is required, as the mapping $\bm{F}\mapsto \bm{f}$ will not always be injective, see also Remark \ref{rem:abstract-elasticity-nonunique-ext} on the same issue.
\item \label{rem:abstract-operator-resolvable-forces-def} 
We stress that the values of the force load $\bm{f}(t)$ in \eqref{eq:loads-abstract-def} must not be from the entire space $\mathcal{X}_0^*$, but must be given from its subspace ${\rm Im}\,{\bf D}$ (closed due to Lemma \ref{lemma:adjoint_operators_lemma} \ref{enum:adjoint_operators_lemma_prop1}). In terms of mechanics this means that given external force $\bm{F}(t)$ must be such that it could be compensated (i.e. brought to an equilibrium state) by the elastic restorative force ${\bf D}\, \bm{\sigma}$ corresponding to at least one stress configuration $\bm{\sigma}$. This is trivially true when ${\bf D}$ is surjective onto $\mathcal{X}_0^*$, but such surjectivity is not always present in mechanical models, as we will see.  We borrow the terminology of \cite[p. 424]{RothWhiteley1981} 
(see also \cite[Sect. 3]{Gudoshnikov2023preprint}) and call $\bm{F}\in \mathcal{X}^*$ and $\bm{f}\in \mathcal{X}_0^*$ (corresponding  to it by \eqref{eq:ae-force-functional-abstract})  {\it resolvable} force loads whenever
\begin{equation}
\bm{f}\in {\rm Im}\,{\bf D}.
\label{eq:resolvable-load-def}
\end{equation}
\end{enumerate}
These discussions allow us to represent the abstract problem of quasi-static evolution in elasticity as a scheme of the Fig. \ref{fig:elasticity-scheme}.
\end{Remark}

\subsection{Stress solution for elasticity in the abstract geometric framework}

To find the stress solution we need to redefine the spaces and operators to take into account ${\bf C}$. Recall, that by Lax-Milgram lemma (see e.g. \cite[Th. 13.29, p. 269]{Kress2014}, \cite[Sect. 1.1, pp. 1--5]{Gatica2014}), the operator ${\bf C}$ has an inverse ${\bf C}^{-1}$, which is also bounded, coercive and symmetric.
\begin{Definition}
\label{def:ae-weighted}
In the setting of Definition \ref{def:ae} denote by $\mathcal{H}_{{\bf C}^{-1}}$ the space $\mathcal{H}$ equipped with the following inner product and norm:
\begin{gather}
\la\bm{\sigma}, \bm{\tau}\ra_{{\bf C}^{-1}} = \la \bm{\sigma}, {\bf C}^{-1} \bm{ \tau} \ra_\mathcal{H} \qquad \text{for any } \bm{\sigma}, \bm{\tau}\in \mathcal{H}, \label{eq:weighted-ip} \\
\|\bm{\sigma}\|_{{\bf C}^{-1}}=\left(\la\bm{\sigma}, \bm{\sigma}\ra_{{\bf C}^{-1}}\right)^{\frac{1}{2}} \qquad \text{for any } \bm{\sigma}\in \mathcal{H}. \label{eq:weighted-norm}
\end{gather}
Further, define the operators
\begin{equation}
\begin{array}{cc}
{\bf E_{C}}:\mathcal{W}_0\subset \mathcal{X} \to \mathcal{H}_{{\bf C}^{-1}}, & 
{\bf D_{C}}:D({\bf D})\subset \mathcal{H}_{{\bf C}^{-1}} \to \mathcal{X}_0^* ,\\[2mm]
{\bf E_{C}} = {\bf C\, E}, & 
{\bf D_{C}}\, \bm{\sigma} ={\bf D}\, \bm{\sigma} \text { for any } \bm{\sigma} \in D({\bf D}).
\end{array}
\label{eq:EC-DC-def}
\end{equation}
The operators ${\bf E_{C}}$ and ${\bf D_{C}}$ themselves fit into the framework of Section \ref{ssect:adjoint-operatorsED}, so that they are also closed, have closed ranges and ${\bf D_{C}}$ is densely defined.  Furthermore, the subspaces
\begin{equation}
\mathcal{U} = {\rm Im}\, {\bf E_C}, \qquad \mathcal{V} = {\rm Ker}\, {\bf D_{C}},
\label{eq:ae-UV-def}
\end{equation}
are closed in $\mathcal{H}_{{\bf C}^{-1}}$, orthogonal in the sense of \eqref{eq:weighted-ip} and
\[
\mathcal{H}_{{\bf C}^{-1}}=\mathcal{U}\oplus\mathcal{V}.
\]
We call $\mathcal{U}$ and $\mathcal{V}$ the {\it fundamental subspaces}, and thoughout the paper we consider them as equipped with the norm of $\mathcal{H}_{{\bf C}^{-1}}$.
Finally, define the operators of orthogonal (in the sense of \eqref{eq:weighted-ip}) projections
\begin{equation}
P_{\mathcal U}:  \mathcal{H}_{{\bf C}^{-1}} \to \mathcal{U}, \qquad P_{\mathcal V}:  \mathcal{H}_{{\bf C}^{-1}} \to \mathcal{V},
\label{eq:ae-PU-PV-def}
\end{equation}
i.e. for any $\bm{\tau}\in  \mathcal{H}_{{\bf C}^{-1}} $ we have a unique Helmholtz- or Beltrami-type decomposition
\begin{equation}
{\bm \tau} = P_{\mathcal{U}}\,{\bm \tau} + P_{\mathcal{V}}\,\bm{\tau}.
\label{eq:projection_decomposition}
\end{equation}
\end{Definition}
The next proposition will be instrumental in finding the stress solution.
\begin{Proposition}
\label{prop:elasticity-affine-planes}
Let $G$ and $Q$ be the following bounded operators. Operator $G$ to be defined by
\begin{equation}
G: \mathcal{H} \to \mathcal{V}\subset  \mathcal{H}_{{\bf C}^{-1}} , \qquad G=P_{\mathcal{V}}\, {\bf C}, 
\label{eq:ae-G-def}
\end{equation}
and
\[
Q: {\rm Im}\,{\bf D} \to \mathcal{U} \subset  \mathcal{H}_{{\bf C}^{-1}} 
\]
be the inverse of the restriction ${\bf D}_{{\bf C},r}$ of ${\bf D_C}$ to $\mathcal{U}$, acting as
\[
{\bf D}_{{\bf C},r}: \mathcal{U}\cap D({\bf D})\subset \mathcal{U}\to {\rm Im}\, {\bf D}.
\]
We have the following facts.
\begin{enumerate}[\it i)]
\item A vector $\bm{\widetilde{\varepsilon}}\in \mathcal{H}$ satisfies \eqref{eq:ae-1} if and only if
\begin{equation}
{\bf C}\, \bm{\widetilde{\varepsilon}} \in \mathcal{U} + G\, {\bm g}.
\label{eq:elasticity-affine-plane1}
\end{equation}
\item A vector $\bm{\widetilde{\sigma}}\in \mathcal{H}$ satisfies \eqref{eq:ae-3} if and only if 
\begin{equation}
\bm{\widetilde{\sigma}} \in \mathcal{V} + Q\, {\bm f}.
\label{eq:elasticity-affine-plane2}
\end{equation}
\end{enumerate}
\end{Proposition}
\noindent {\bf Proof.} {\it i)} Since the operator ${\bf C}$ is invertible, equation \eqref{eq:ae-1} is equivalent to 
\[
{\bf C}\, \bm{\widetilde{\varepsilon}} \in \mathcal{U} + {\bf C}\, {\bm g},
\]
which, by the decomposition \eqref{eq:projection_decomposition} of ${\bf C}\, {\bm g}$, is equivalent to \eqref{eq:elasticity-affine-plane1}.

\noindent {\it ii)} As ${\bf D}_{\bf C}$ is closed and densely defined, we follow the lines of \cite[Appendix A]{Beutler1976} to show that the operator ${\bf D}_{{\bf C},r}$ is closed, densely defined and bijective. Hence there exists the inverse $Q$ of ${\bf D}_{{\bf C},r}$, and the continuity of $Q$ follows from the closed graph theorem (see e.g. \cite[Th. 2.9, p. 37]{Brezis2011}). Thus  \eqref{eq:ae-3} is equivalent to \eqref{eq:elasticity-affine-plane2}. $\blacksquare$

At last, we have a convenient expression of the stress solution for elasticity. 
\begin{Theorem} 
\label{th:elasticity-solution}
Function $\bm{\widetilde{\sigma}}(t)$ is a stress solution for elasticity if and only if
\begin{equation}
\bm{\widetilde{\sigma}}(t)\in ( \mathcal{U}+G\, {\bm g}(t)) \cap( \mathcal{V}+Q\,{\bm f}(t)).
\label{eq:abstract_intersection}
\end{equation}
Inclusion \eqref{eq:abstract_intersection} has exactly one solution 
\begin{equation}
\bm{\widetilde{\sigma}}(t) = G\, {\bm g}(t) + Q\,{\bm f}(t).
\label{eq:abstract_linear_solution}
\end{equation} 
\end{Theorem}

\noindent {\bf Proof.} 
Inclusion \eqref{eq:abstract_intersection} is equivalent to \eqref{eq:ae-1}--\eqref{eq:ae-3} by Proposition \ref{prop:elasticity-affine-planes}. To derive \eqref{eq:abstract_linear_solution} we notice, that \eqref{eq:abstract_intersection} is also equivalent to
\[
\bm{\widetilde{\sigma}}- G\, \bm{g} -Q\, \bm{f} \in (\mathcal{U}+G\, \bm{g}) \cap(\mathcal{V}+Q\, \bm{f})-G\,  \bm{g} -Q\, \bm{f},
\]
In turn, the right-hand side is
\begin{multline*}
(\mathcal{U}+G\,  \bm{g}) \cap(\mathcal{V}+Q\, \bm{f})-G\, \bm{g} -Q\, \bm{f}  =\\= (\mathcal{U}+G\,  \bm{g}-G\, \bm{g} -Q\, \bm{f}) \cap(\mathcal{V}+Q\, \bm{f}-G\,  \bm{g} -Q\, \bm{f}) =\\= (\mathcal{U}-Q\, \bm{f}) \cap(\mathcal{V}-G\,  \bm{g})=\mathcal{U}\cap\mathcal{V}=\{0\},
\end{multline*}
where the last two equalities are due to $Q\, \bm{f} \in \mathcal{U},\,G\, \bm{g}\in \mathcal{V}$ and the triviality of the intersection of the orthogonal subspaces, respectively. Therefore, \eqref{eq:abstract_linear_solution} is indeed equivalent to \eqref{eq:abstract_intersection}. $\blacksquare$

\begin{Remark}
 While we wrote the problem of Definition~\ref{def:ae} as time-dependent and {\it quasi-static}, its solution \eqref{eq:abstract_linear_solution} is independently determined at each time-instant $t\in I$ by the values of the loads $\bm{f}(t)$ and $\bm{g}(t)$.
In other words, one can see Definition~\ref{def:ae} as a family of uncoupled {\it static} problems, depending on $t$ as a parameter. For this reason, in the context of elasticity, we use the terms {\it static} and {\it quasi-static} interchangeably in the current Section \ref{sect:elasticity}. This will no longer be valid for the problems with Prandtl-Reuss elastoplasticity, which we will consider in Section \ref{sect:perfect-plasticity} and beyond, because those are proper history-dependent evolution problems, and they cannot be reduced to {\it uncoupled} problems for different $t$. However, stress solution for elasticity $\widetilde{\sigma}$ in its time-dependent form \eqref{eq:abstract_linear_solution} will play a major role for solving the elastoplasticity problems. In particular, we will use its time-regularity, described in the next remark.
\end{Remark}

\begin{Remark} 
\label{rem:ae-time-regularity}
When the loads have time-regularity \eqref{eq:loads-abstract-def}, the stress solution for elasticity also has time-regularity $W^{1, \infty}(I, \mathcal{H}_{{\bf C}^{-1}}) \cong W^{1, \infty}(I, \mathcal{H})$. Indeed, by continuity of operators $G$ and $Q$ we have by Proposition \ref{prop:prelim-classical-derivatives-ae} that for a.a. $t\in I$
\[
\frac{d}{dt}\, \bm{\widetilde{\sigma}}(t) = \lim_{\Delta\to 0} \frac{\bm{\widetilde{\sigma}}(t+\Delta)-\bm{\widetilde{\sigma}}(t)}{\Delta}=G \frac{d}{dt}\, \bm{g}(t)+ Q\,\frac{d}{dt}\, \bm{f}(t) 
\] 
and, therefore,
\[
\left\|\frac{d}{dt}\, \bm{\widetilde{\sigma}}(t)\right\|_{{\bf C}^{-1}} \leqslant \|G\|_{\rm op} \,\left\|\frac{d}{dt}\bm{g}(t)\right\|_{\mathcal{H}}+\|Q\|_{\rm op}\,\left\|\bm{f}(t)\right\|_{\mathcal{X}_0^*},
\]
where the operator norms are
\[
\|G\|_{\rm op} = \sup_{\bm{\tau}\in \mathcal{H}\setminus\{0\}} \frac{\|G \,\bm{\tau}\|_{{\bf C}^{-1}}}{\|\bm{\tau}\|_{\mathcal{H}}}, \qquad \|Q\|_{\rm op} = \sup_{\bm{\varphi}\in {\rm Im}\, {\bf D}\setminus \{0\}}\frac{ \|Q\, \bm{\varphi}\|_{{\bf C}^{-1}}}{ \|\bm{\varphi}\|_{\mathcal{X}_0^*}}.
\]
We conclude that
\begin{equation}
\left\|\widetilde{\sigma}\right\|_{W^{1,\infty}\left(I, \mathcal{H}_{{\bf C}^{-1}}\right)}\leqslant \|G\|_{\rm op}\, \|g\|_{W^{1, \infty}(I, \mathcal{H})} +  \|Q\|_{\rm op} \, \|f\|_{W^{1, \infty}\left(I, \mathcal{X}_0^*\right)}.
\label{eq:ae-time-regularity}
\end{equation}
\end{Remark}

\begin{Remark}
\label{remark:to-find-elastic-strain}
 In the problem of elasticity the unknown $\bm{\widetilde{\varepsilon}}$ can be recovered from \eqref{eq:ae-2} by using bounded operator ${\bf C}^{-1}$.
\end{Remark}

\subsection{Example 1 --- a discrete model with one-dimensional \texorpdfstring{$\mathcal{V}$}{V}}
\label{ssect:ex11-elasticity}
\begin{figure}[H]\center
\includegraphics{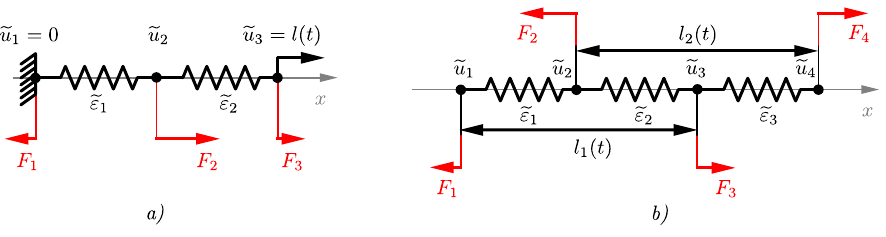}
\caption{\footnotesize Discrete models of Examples 1 ({\it a}) and 2 ({\it b}). Red arrows denote the external forces $\bm{F}$, applied at the nodes.
} 
\label{fig:discrete-models-elastic}
\end{figure} 
As the first example, consider the discrete model of Fig. \ref{fig:discrete-models-elastic} {\it a)}, which consists of two elastic springs connected in series between three nodes. The left endpoint is held at $x=0$ and the displacement of the right endpoint from its position in the relaxed configuration is prescribed as a known function $l\in W^{1,\infty}(I)$.

\subsubsection{Example 1 --- kinematics}
In terms of the abstract framework of Section \ref{ssect:abstract-elasticity} we have the following quantities. The ambient spaces are given as
\[
\mathcal{X} = \mathbb{R}^3, \qquad \mathcal{H} = \mathbb{R}^2,
\]
where $\mathcal{X}$ represents all longitudinal {\it displacements} of the three nodes from a stress-free reference configuration, and $\mathcal{H}$ represents all {\it elongations} of the two springs from the stress-free reference lengths.  We write the boundary condition as the following constraint on the vector of  displacements $\bm{\widetilde u}\in \mathcal{X}$:
\begin{equation}
R\, \bm{\widetilde u} + r(t) =0,
\label{eq:ex11-bc}
\end{equation}
where
\[
R=\begin{pmatrix}
1 & 0 & 0\\
0 & 0 & 1
\end{pmatrix}, \qquad r(t)= \begin{pmatrix} 0\\ -l(t) \end{pmatrix}.
\]
Further, denote the following matrix
\begin{equation}
{\rm E} = \begin{pmatrix}
-1 & 1 & 0 \\
0 & -1 & 1
\end{pmatrix}.
\label{ae:ex11-E-matrix}
\end{equation}
The elongations $\bm{\widetilde{\varepsilon}}  \in \mathcal{H}$ of the springs are computed for each $t\in I$ via the expression \eqref{eq:ae-compatibility-abstract}. We define the operator ${\bf E}$ from \eqref{eq:ae-abstract-E-def} by the formula \eqref{eq:ex11-E-def} with 
\[
\mathcal{W} = \mathcal{X}=\mathbb{R}^3,\qquad \mathcal{W}_0=D({\bf E}) = {\rm Ker}\, R.
\]

\noindent Since $\mathcal{W}_0$ is closed we have $\mathcal{X}_0=\mathcal{W}_0$. It is a one-dimensional space, parametrized by $u_2\in \mathbb{R}$, so that
\begin{equation}
\mathcal{W}_0=\mathcal{X}_0=D({\bf E})={\rm Ker}\, R= {\rm Im}\, R_0 = \left\{R_0\, 
 u_2:u_2\in \mathbb{R}
\right\} \qquad \text{with }R_0=\begin{pmatrix} 0 \\ 1 \\ 0 \end{pmatrix}.
\label{eq:ex11-DE-parametrization}
\end{equation}
To derive \eqref{eq:ae-1} we write \eqref{eq:ex11-bc} in the  equivalent form \eqref{eq:ae-bc-abstract} by using the Moore-Penrose pseudoinverse $R^+$ of $R$ (see Corollary \ref{cor:mp-solving-linear-systems})
\[
R^+= R^{\top}\left(RR^{\top}\right)^{-1}= \begin{pmatrix}1 & 0 \\ 0 & 0 \\ 0 & 1\end{pmatrix}
\]
and setting
\[
\bm{u}_{\bf D} (t) = -R^+ r(t) = \begin{pmatrix} 0\\0\\ l(t)\end{pmatrix}.
\]
Thus, we follow \eqref{eq:ae-abstract-g-def} and let 
\[
\bm{g}(t) ={\rm E}\,\bm{u}_{\bf D} (t)  = -  {\rm E}\, R^+ r(t) = \begin{pmatrix}  0 \\  l(t)\end{pmatrix}.  
\]

\noindent We have set up all of the quantities, which appear in \eqref{eq:ae-1} and \eqref{eq:ae-general-strain}--\eqref{eq:ae-abstract-g-def}.

\subsubsection{Example 1 --- Hooke's law}
In the two-spring the model of Fig. \ref{fig:discrete-models-elastic} {\it a)} we have the operator ${\bf C}$ from \eqref{eq:ae-C-def} and its inverse defined by the respective diagonal matrices
\[
{\bf C} =\begin{pmatrix} 
k_1 & 0\\[1mm] 0 &  k_2
\end{pmatrix}, \qquad {\bf C}^{-1} =\begin{pmatrix} 
k_1^{-1} & 0\\[1mm] 0 &  k_2^{-1}
\end{pmatrix},
\]
in which $k_1, k_2>0$ are the {\it stiffness} parameters of the springs. Equation \eqref{eq:ae-2} in this model is Hooke's law, connecting the {\it stresses} $\bm{\widetilde{\sigma}}\in \mathbb{R}^2$ of the springs with the elongations $\bm{\widetilde{\varepsilon}}$.

\subsubsection{Example 1 --- statics}
Assume that a given {\it external force} $F\in W^{1,\infty}(I, \mathbb{R}^3)$ is applied in the model of Fig. \ref{fig:discrete-models-elastic} {\it a)}. Specifically, for each $t\in I$, each of the three components of $\bm{F}$ denotes a longitudinal force, applied at the respective node. 
From the {\it principle of virtual work} (see e.g.  \cite[Section 1.4, pp.~16--17]{Goldstein2002}) we write the equation of equilibrium on the unknown $\bm{\widetilde{\sigma}}\in \mathbb{R}^2$ as 
\begin{equation}
\left\la - {\rm E}^\top \,\bm{\widetilde{ \sigma}} + \bm{F}, \bm{\delta u} \right\ra_{\mathbb{R}^3} = 0 \qquad \text{for any }\bm{\delta u}\in {\rm Ker}\, R.
\label{eq:ex11-equilibrium-basic-formulation}
\end{equation}
We use the convention, that positive values of $\bm{\widetilde{\varepsilon}}$ denote a configuration with stretched springs (with respect to the stress-free reference configuration). Therefore, Hooke's law \eqref{eq:ae-2} implies that positive values of $\bm{\widetilde{\sigma}}$ correspond to the forces acting towards {\it contraction}. Substitute ${\rm E}$ from \eqref{ae:ex11-E-matrix} to observe that the term  $ - {\rm E}^\top \,\bm{\widetilde{ \sigma}}$ in \eqref{eq:ex11-equilibrium-basic-formulation} indeed expresses the forces at the nodes, produced by a stress vector $\bm{\widetilde{\sigma}}$.

Now we transform equation \eqref{eq:ex11-equilibrium-basic-formulation} to the abstract form \eqref{eq:ae-3}. For each value $\bm{F}(t)\in \mathbb{R}^3$ we define a functional $\bm{f}(t)\in \mathcal{X}_0^*$ acting as
\begin{equation}
{\bm f}: \bm{u} \mapsto \left\la \bm{F},  \bm{u}  \right\ra_{\mathbb{R}^3} \qquad \text{for any } \bm{u}\in \mathcal{X}_0,
\label{eq:ex11-functionals-def}
\end{equation}
which coincides with \eqref{eq:ae-force-functional-abstract} in the current setting. 
In fact, due to the representation \eqref{eq:ex11-DE-parametrization}, we can write \eqref{eq:ex11-functionals-def} explicitly as
\begin{equation}
{\bm f}: \bm{u} \mapsto F_2\,  u_2 \qquad \text{for any } \bm{u}\in \mathcal{X}_0.
\label{eq:ex11-functionals-explicit}
\end{equation}

\noindent We can now rewrite the equation of equilibrium \eqref{eq:ex11-equilibrium-basic-formulation} as
\begin{equation}
\left\la \bm{\widetilde{\sigma}}, {\bf E}\, \bm{\delta u}\right\ra_{\mathcal{H}} = \la\bm{f}, \bm{\delta u}\ra_{\mathcal{X}_0^*, \mathcal{X}_0} \qquad \text{for any } \bm{\delta u}\in \mathcal{X}_0,
\label{eq:ex11-equilibrium-functional-formulation}
\end{equation}
which becomes \eqref{eq:ae-3} when we take the operator ${\bf D}$ as the adjoint, according to the Definition \ref{def:abstract-adjoint}. The adjoint operator has its domain 
\[
D({\bf D})= \mathcal{H}^*\cong \mathcal{H}=\mathbb{R}^2
\]
and it is defined by
\begin{gather}
{\bf D}: \mathcal{H}^* \to \mathcal{X}_0^*,
\label{eq:ex11-D-def1}
\\
\left\la{\bf D}\, \bm{\sigma}, \bm{\delta u}\right \ra_{\mathcal{X}_0^*,\mathcal{X}_0} = \left\la \bm{\sigma}, {\bf E}\, \bm{\delta u}\right\ra_{\mathcal{H}} \qquad \text{for any }  \bm{\sigma}\in \mathcal{H},\, \bm{\delta u} \in \mathcal{X}_0.
\label{eq:ex11-D-def2}
\end{gather}

\begin{Remark}
\label{rem:ex11-all-forces-are-resolvable}
We can use the representation \eqref{eq:ex11-DE-parametrization} to observe for any $\bm{\sigma} = \begin{pmatrix} \sigma_1\\ \sigma_2 \end{pmatrix} \in \mathbb{R}^2,\, \bm{u}\in  \mathcal{X}_0\subset \mathbb{R}^3$  that
\begin{equation}
\left\la
{\bf D} \bm{\sigma}, \bm{u}
\right\ra_{\mathcal{X}_0^*,\mathcal{X}_0} = 
\left\la \begin{pmatrix} \sigma_1\\ \sigma_2 \end{pmatrix}, {\rm E} R_0\, u_2 \right\ra_{\mathbb{R}^2}
=  \begin{pmatrix} 1 & -1 \end{pmatrix} \begin{pmatrix} \sigma_1 \\ \sigma_2 \end{pmatrix} u_2 =(\sigma_1-\sigma_2)\,u_2,
\label{eq:ex11-explicit-D}
\end{equation}
from which we deduce that ${\bf D}$ is surjective onto $\mathcal{X}_0^*$, since both ${\rm Im}\, {\bf D}$ and $\mathcal{X}_0^*$ are one-dimensional (respectively, due to \eqref{eq:ex11-explicit-D}, and due to $\mathcal{X}_0$ being one-dimensional, see \eqref{eq:ex11-DE-parametrization}). In terms of mechanics of the model of Fig. \ref{fig:discrete-models-elastic} {\it a)}, the surjectivity of ${\bf D}$ means  that any external force $\bm{F}\in \mathbb{R}^3$ is resolvable, see Remark \ref{rem:abstract-operator-to-restrict} \ref{rem:abstract-operator-resolvable-forces-def}.
\end{Remark}

\begin{Remark}
\label{rem:ex11-elasticity-nonunique-ext}
We observed earlier that the space $\mathcal{X}_0^*$ is one-dimensional, as the dual to the one-dimensional space $\mathcal{X}_0$. At the same time, we just have shown that any $\bm{F}$ from the three-dimensional space $\mathbb{R}^3$ can be compensated. This can be explained by the fact that the mapping 
\[
\bm{F} \mapsto \bm{f}
\]
\[ 
\mathbb{R}^3 \to \mathcal{X}_0^*,
\]
defined by \eqref{eq:ex11-functionals-def} is not injective. Indeed, for any 
\[
\bm{F_0} \in \mathcal{X}_0^\perp = \left({\rm Ker}\, R\right)^\perp = {\rm Im}\, R^\top
\]
both $\bm F$ and $\bm{F}+\bm{F_0}$ define the same functional by \eqref{eq:ex11-functionals-def}. Put differently, a functional $\bm{f}\in \mathcal{X}_0^*$ admits many such extensions $\bm{F}+\bm{F_0}\in \mathcal{X}^*$, cf. Remark \ref{rem:abstract-elasticity-nonunique-ext} and Remark \ref{rem:abstract-operator-to-restrict} \ref{rem:abstract-operator-to-restrict-about-forces}. Notice, that from a physical point of view such $\bm{F_0}$ represents a {\it reaction} of the constraint \eqref{eq:ex11-bc}, which gets ``absorbed'' when we pass to the formalism of $\bm{f}\in \mathcal{X}_0^*$
\end{Remark}

\subsubsection{Example 1 --- the fundamental subspaces and the stress solution}
In the above we have provided all of the data required in Definition \ref{def:ae} to set up a problem of the type \eqref{eq:ae-1}-\eqref{eq:ae-3}. Now we will give exact expressions of the quantities from Definition \ref{def:ae-weighted} and Theorem \ref{th:elasticity-solution} for the model of Fig. \ref{fig:discrete-models-elastic} {\it a)}.

\begin{figure}[H]\center
\includegraphics{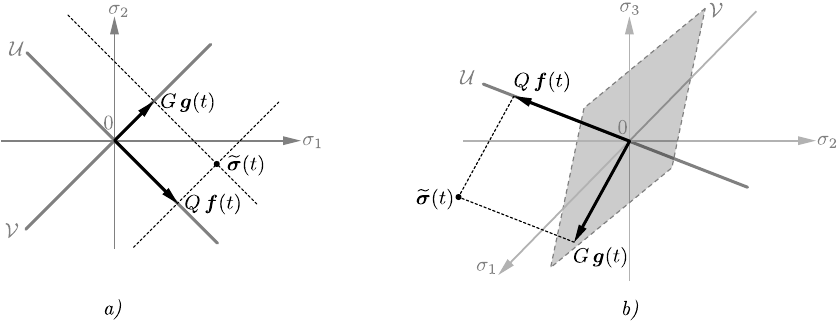}
\caption{\footnotesize The fundamental spaces and the stress solution for elasticity in the discrete models of Example 1 ({\it a}) and Example 2 ({\it b}). The figure shows the situation with the stiffness parameters $k_i$ equal to $1$.
} 
\label{fig:discrete-models-spaces}
\end{figure}

Space $\mathcal{H}_{{\bf C}^{-1}}$ is $\mathbb{R}^2$ equipped with the inner product
\begin{equation}
\left\la \begin{pmatrix} \sigma_1\\ \sigma_2\end{pmatrix}, \begin{pmatrix}\tau_1\\ \tau_2\end{pmatrix} \right\ra_{{\bf C}^{-1}} = \sigma_1 k_1^{-1} \tau_1+ \sigma_2 k_2^{-1} \tau_2 \qquad \text{for all }\begin{pmatrix} \sigma_1\\ \sigma_2\end{pmatrix}, \begin{pmatrix}\tau_1\\ \tau_2\end{pmatrix}\in \mathbb{R}^2.
\label{eq:ex11-inner-product}
\end{equation}
Subspace  $\mathcal{U}$ is
\[
\mathcal{U}  ={\rm Im}\,{\bf C} {\rm E} R_0  = {\rm Im}\, \begin{pmatrix}k_1 & 0 \\ 0 & k_2\end{pmatrix}  \begin{pmatrix} -1 & 1 & 0\\ 0 & -1 & 1 \end{pmatrix}  \begin{pmatrix} 0 \\ 1 \\ 0 \end{pmatrix}  = {\rm Im}\,  \begin{pmatrix} k_1 \\ -k_2  \end{pmatrix},
\]
and, from \eqref{eq:ex11-explicit-D} used in the definition of ${\bf D}_{\bf C}$ in \eqref{eq:EC-DC-def}, we have that
\begin{equation}
\mathcal{V} = {\rm Im}\, \begin{pmatrix} 1 \\ 1\end{pmatrix},
\label{eq:ex11-space-V-exp}
\end{equation}
see Fig. \ref{fig:discrete-models-spaces} {\it a)}.

By the formula for the projection on a subspace in terms of a weighted inner product (see e.g. \cite[(5.31), p. 247]{OlverLinearAlgebra2018}) we have the matrices of projections
\begin{multline*}
P_{\mathcal{U}}= \begin{pmatrix} k_1 \\ -k_2\end{pmatrix}\left(\begin{pmatrix} k_1 & -k_2\end{pmatrix}\begin{pmatrix} k_1^{-1} & 0\\ 0 & k_2^{-1}\end{pmatrix}\begin{pmatrix} k_1 \\ -k_2\end{pmatrix}\right)^{-1}\begin{pmatrix} k_1 & -k_2\end{pmatrix}\begin{pmatrix} k_1^{-1} & 0\\ 0 &  k_2^{-1} \end{pmatrix}=\\ = \frac{1}{k_1+k_2}\begin{pmatrix} k_1 &  -k_1 \\ -k_2 & k_2\end{pmatrix},
\end{multline*}
\[
P_{\mathcal{V}} = \begin{pmatrix} 1 \\ 1\end{pmatrix}\left(\begin{pmatrix} 1 & 1\end{pmatrix}\begin{pmatrix} k_1^{-1}& 0 \\ 0 & k_2^{-1}\end{pmatrix}\begin{pmatrix} 1 \\ 1\end{pmatrix}\right)^{-1}\begin{pmatrix} 1 & 1\end{pmatrix}\begin{pmatrix} k_1^{-1}& 0 \\ 0 & k_2^{-1}\end{pmatrix} = \frac{1}{k_1^{-1}+k_2^{-1}}\begin{pmatrix} k_1^{-1}&  k_2^{-1}\\[1mm] k_1^{-1} & k_2^{-1}\end{pmatrix}.
\]
Further, operator $G$ is represented by the matrix
\[
G= P_{\mathcal{V}}{\bf C} = \frac{1}{k_1^{-1}+k_2^{-1}}\begin{pmatrix} k_1^{-1}&  k_2^{-1}\\ k_1^{-1} & k_2^{-1}\end{pmatrix} \begin{pmatrix} k_1& 0 \\ 0 & k_2\end{pmatrix} =\frac{1}{k_1^{-1}+k_2^{-1}}\begin{pmatrix} 1& 1 \\ 1 &1\end{pmatrix}.
\]

Now we only have to find a representation for the operator $Q$ in terms of an  external force $\bm{F}$. Notice, that for any $\bm{\sigma}\in \mathcal{U}$  there exists a unique $\sigma\in \mathbb{R}$, such that 
\[
\bm{\sigma}= \begin{pmatrix}k_1\\ -k_2\end{pmatrix} \sigma. 
\]
Plug that into  \eqref{eq:ex11-explicit-D}, and then plug the result together with \eqref{eq:ex11-functionals-explicit} into \eqref{eq:ae-3}  to see that
\[
(k_1+k_2)\sigma = F_2,
\]
i. e.
\[
\sigma = \frac{F_2}{k_1+k_2}.
\]
Therefore, for a given external force $\bm{F}\in \mathbb{R}^{3}$ and the corresponding (by \eqref{eq:ex11-functionals-def}) functional $\bm{f}\in \mathcal{X}_0^*$ we have the value of the operator $Q$ 
\[
Q\, \bm{f} = \frac{F_2}{k_1+k_2}\begin{pmatrix}k_1 \\ -k_2\end{pmatrix}\in \mathcal{U}.
\]

By Theorem \ref{th:elasticity-solution} we conclude that the stress solution for elasticity in the model of Fig.  \ref{fig:discrete-models-elastic} {\it a)} has the closed form
\begin{equation}
\bm{\widetilde{\sigma}}(t) = G\, \bm{g}(t) + Q\, \bm{f}(t) = \frac{l(t)}{k_1^{-1}+ k_2^{-1}} \begin{pmatrix}1 \\ 1\end{pmatrix} + \frac{F_2(t)}{k_1+k_2}\begin{pmatrix} k_1\\-k_2\end{pmatrix},
\label{eq:ex11-linear-solution}
\end{equation}
see Fig. \ref{fig:discrete-models-spaces} {\it a)}.
\subsection{Example 2 --- a discrete model with two-dimensional \texorpdfstring{$\mathcal{V}$}{V}}
\label{ssect:ex12-elasticity}
Our second example is the discrete model of Fig. \ref{fig:discrete-models-elastic} {\it b)} with three elastic springs connected in series and prescribed constraints on elongations
\begin{equation}
\begin{array}{c}
\widetilde{\varepsilon}_1+ \widetilde{\varepsilon}_2 = l_1(t),\\
\widetilde{\varepsilon}_2+ \widetilde{\varepsilon}_3 = l_2(t),
\end{array}
\label{eq:ex12-constraint}
\end{equation}
where $l_1, l_2\in W^{1, \infty}(I)$ are given.

\subsubsection{Example 2 --- kinematics}
We proceed similarly to the previous example and start with the ambient spaces
\[
\mathcal{X}=\mathbb{R}^4, \qquad \mathcal{H} = \mathbb{R}^3
\]
representing displacements $\bm{\widetilde{u}}\in \mathcal{X}$ of the four nodes and elongations $\bm{\widetilde{\varepsilon}}\in \mathcal{H}$ of the three springs respectively. The relation \eqref{eq:ae-compatibility-abstract} between them is now defined via the matrix
\[
{\rm E} =  \begin{pmatrix}
-1 & 1 & 0 & 0\\
0 & -1 & 1 & 0\\
0 & 0 & -1 & 1
\end{pmatrix}.
\]
But, instead of the boundary condition \eqref{eq:ex11-bc}, we write the constraint \eqref{eq:ex12-constraint} in terms of $\bm{\widetilde{\varepsilon}}\in \mathcal{H}$ as
\begin{equation}
R\,\bm{\widetilde{\varepsilon}} + r(t) = 0
\label{eq:ex12-constraint-matrix-form}
\end{equation}
with
\[
R= \begin{pmatrix}
1 & 1 & 0\\
0 & 1 & 1
\end{pmatrix}, \qquad r(t) = -\begin{pmatrix}l_1(t) \\ l_2(t) \end{pmatrix}.
\]
In terms of $\bm{\widetilde{u}}\in \mathcal{X}$ the constraint can be written as
\begin{equation}
R {\rm E}\,\bm{\widetilde{u}} +r(t) =0,
\label{eq:ex12-constrain-u-form}
\end{equation}
where 
\[
R {\rm E} = \begin{pmatrix}-1 & 0 & 1& 0\\ 0 & -1& 0 & 1 \end{pmatrix}.
\]

\noindent We set the spaces 
\[
\mathcal{W}=\mathcal{X}=\mathbb{R}^4, \qquad \mathcal{W}_0= D({\bf E}) ={\rm Ker}\, (R {\rm E}),
\]
and by \eqref{eq:ex11-E-def} we set the operator ${\bf E}$ from \eqref{eq:ae-abstract-E-def}.

We again have $\mathcal{W}_0=\mathcal{X}_0$ with a parametrization
\begin{equation}
\mathcal{W}_0=\mathcal{X}_0=D({\bf E}) = {\rm Im}\, R_0 = \left\{R_0 \begin{pmatrix}u_1\\ u_2 \end{pmatrix}: \begin{pmatrix} u_1\\ u_2 \end{pmatrix}\in \mathbb{R}^2\right\}\qquad \text{with } R_0 = \begin{pmatrix}1& 0\\ 0 & 1\\ 1& 0 \\ 0 & 1\end{pmatrix}.
\label{eq:ex12-DE-parametrization}
\end{equation}

Similarly to the previous example, to derive \eqref{eq:ae-1} we write \eqref{eq:ex12-constrain-u-form} in the equivalent form \eqref{eq:ae-bc-abstract} with the use (as in Corollary \ref{cor:mp-solving-linear-systems}) of Moore-Penrose pseudoinverse $(R {\rm E})^+$:
\[
(R {\rm E})^+= (R{\rm E})^{\top}\left(R{\rm E}(R{\rm E})^{\top}\right)^{-1}= \frac{1}{2}\begin{pmatrix}-1 & 0 \\ 0 & -1 \\ 1 & 0 \\ 0 & 1\end{pmatrix}.
\]
We set
\[
\bm{u}_{\bf D} (t) = -(R{\rm E})^+ r(t) = \frac{1}{2}\begin{pmatrix}-l_1(t)\\-l_2(t)\\l_1(t)\\l_2(t)\end{pmatrix}.
\]
Again, we follow \eqref{eq:ae-abstract-g-def} and let 
\begin{equation}
\bm{g}(t) ={\rm E}\,\bm{u}_{\bf D} (t)  = - {\rm E}\, (R{\rm E})^+ r(t) =  \frac{1}{2}\begin{pmatrix}1 & -1 \\ 1 & 1 \\ -1 & 1\end{pmatrix} \begin{pmatrix}  l_1(t) \\  l_2(t)\end{pmatrix}.  
\label{eq:ex12-g-def}
\end{equation}
\noindent For the model of Fig. \ref{fig:discrete-models-elastic} {\it b)}  we have set up all of the quantities in \eqref{eq:ae-1} and \eqref{eq:ae-general-strain}--\eqref{eq:ae-abstract-g-def}.

\begin{Remark} In the particular (and the currently present) case of the matrix ${\rm E}$ being surjective onto $\mathcal{H}$, instead of \eqref{eq:ex12-g-def} one can use the simpler expression
\begin{equation*}
\bm{g}(t)= - R^+ r(t),
\end{equation*}
which leads to 
\begin{equation}
\bm{g}(t) =\frac{1}{3}\begin{pmatrix}2 & -1 \\ 1 & 1 \\ -1 & 2\end{pmatrix} \begin{pmatrix}  l_1(t) \\  l_2(t)\end{pmatrix}. 
\label{eq:ex12-g-def-alt-concrete}
\end{equation}
when the current values are plugged in. One can verify the equivalence between \eqref{eq:ex12-g-def} and \eqref{eq:ex12-g-def-alt-concrete} in terms of \eqref{eq:ae-1} by checking that their difference belongs to the subspace ${\rm Im}\,{\bf E} = {\rm Im}\, {\rm E}R_0$.
\end{Remark}
\subsubsection{Example 2 --- Hooke's law}
In the three-spring the model of Fig. \ref{fig:discrete-models-elastic} {\it b)} we have the operator ${\bf C}$ from \eqref{eq:ae-C-def} and its inverse defined by the respective diagonal matrices
\[
{\bf C} =\begin{pmatrix} 
k_1 & 0 & 0\\[1mm] 0 &  k_2 & 0\\[1mm] 0 & 0 & k_2
\end{pmatrix}, \qquad {\bf C}^{-1} =\begin{pmatrix} 
k_1^{-1} & 0 & 0\\[1mm] 0 &  k_2^{-1} & 0\\[1mm] 0 & 0 & k_3^{-1}
\end{pmatrix},
\]
in which $k_1, k_2, k_3>0$. The stress vectors $\bm{\widetilde{\sigma}}$, defined via \eqref{eq:ae-2}, are now from $\mathbb{R}^3$.
\subsubsection{Example 2 --- statics}
Similarly to the previous example, we assume that a given external force $F\in W^{1,\infty}(I, \mathbb{R}^4)$ is applied in the model of Fig. \ref{fig:discrete-models-elastic} {\it b)}, so that for each $t\in I$ each of the four components of $\bm{F}$ denotes a longitudinal force, applied at the respective node. The equilibrium equation with an unknown stress $\bm{\widetilde{\sigma}}\in \mathbb{R}^3$ is 
\begin{equation}
\left\la - {\rm E}^\top \,\bm{\widetilde{ \sigma}} + \bm{F}, \bm{\delta u} \right\ra_{\mathbb{R}^4} = 0 \qquad \text{for any }\bm{ \delta u}\in {\rm Ker}\, {\rm E}R,
\label{eq:ex12-equilibrium-basic-formulation}
\end{equation}
Again, we need to transform equation \eqref{eq:ex12-equilibrium-basic-formulation} to the abstract form \eqref{eq:ae-3}. For each $\bm{F}\in \mathbb{R}^4$ we define the corresponding functional $\bm{f}\in \mathcal{X}_0^*$ by the formula
\begin{equation}
{\bm f}: \bm{u} \mapsto \left\la \bm{F},  \bm{u}  \right\ra_{\mathbb{R}^4} \qquad \text{for any } \bm{u}\in \mathcal{X}_0.
\label{eq:ex12-functionals-def}
\end{equation}
 The parametrization \eqref{eq:ex12-DE-parametrization} yields the following explicit formula for $\bm{f}$:
\begin{equation}
{\bm f}: \bm{u} \mapsto \begin{pmatrix}F_1+ F_3& F_2+ F_4\end{pmatrix}\begin{pmatrix}u_1\\ u_2\end{pmatrix} \qquad \text{for any } \bm{u}\in \mathcal{X}_0.
\label{eq:ex12-functionals-explicit}
\end{equation}

Again, we rewrite the equilibrium equation as \eqref{eq:ex11-equilibrium-functional-formulation}, which becomes \eqref{eq:ae-3} when we take the operator ${\bf D}$ as the adjoint, according to the Definition \ref{def:abstract-adjoint}.
In other words, the adjoint operator is defined by \eqref{eq:ex11-D-def1}--\eqref{eq:ex11-D-def2} on its domain 
\[
D({\bf D}) =\mathcal{H}^*\cong\mathcal{H}=\mathbb{R}^3.
\]

\begin{Remark}
We can use the representation \eqref{eq:ex12-DE-parametrization} to observe for any $\bm{\sigma}=\begin{pmatrix} \sigma_1\\ \sigma_2\\ \sigma_3 \end{pmatrix} \in \mathbb{R}^3,\, \bm{u}\in  \mathcal{X}_0\subset \mathbb{R}^4$  that
\begin{multline}
\left\la
{\bf D} \bm{\sigma}, \bm{u}
\right\ra_{\mathcal{X}_0^*,\mathcal{X}_0} = 
\left\la \begin{pmatrix} \sigma_1\\ \sigma_2\\ \sigma_3 \end{pmatrix}, {\rm E} R_0 \begin{pmatrix}u_1\\ u_2\end{pmatrix} \right\ra_{\mathbb{R}^3}
=  \begin{pmatrix} \sigma_1 &  \sigma_2 & \sigma_3 \end{pmatrix}\begin{pmatrix}-1 & 1\\ 1 & -1 \\ -1 & 1\end{pmatrix}\begin{pmatrix}u_1\\u_2\end{pmatrix}=\\=
 \begin{pmatrix} \sigma_1 &  \sigma_2 & \sigma_3 \end{pmatrix}\begin{pmatrix}-1 \\ 1 \\ -1 \end{pmatrix}\left(u_1 - u_2\right)= (-\sigma_1+\sigma_2-\sigma_3)\left(u_1 - u_2\right).
\label{eq:ex12-explicit-D}
\end{multline}
We can see that ${\bf D}$ is not surjective, since ${\rm Im}\, {\bf D}$ is one-dimensional  (due to \eqref{eq:ex12-explicit-D}), but $\mathcal{X}_0^*$ is two-dimensional (due to two-dimensional $\mathcal{X}_0$ in \eqref{eq:ex12-DE-parametrization}), cf. Remark \ref{rem:ex11-all-forces-are-resolvable}. Therefore, we indeed must specifically require the values of $\bm{f}(t)$ to satisfy \eqref{eq:resolvable-load-def}, otherwise the equation of equilibrium \eqref{eq:ae-3} cannot be solved. Notice, that to have \eqref{eq:ex12-functionals-explicit} of the type \eqref{eq:ex12-explicit-D} we must take $\bm{F}$ such that
\begin{equation}
F_1 + F_3 = -(F_2+ F_4)
\label{eq:ex12-F-restriction-form1}
\end{equation}
i.e.
\begin{equation}
F_1+F_2+F_3+F_4 =0.
\label{eq:ex12-F-restriction-form2}
\end{equation}
This agrees with the physical reasoning, that since the model of Fig.  \ref{fig:discrete-models-elastic} {\it b)} is not fixed in place, an external force can be compensated by internal stresses (i.e. the force is resolvable) only when the total applied momentum is zero. 
\end{Remark}

\begin{Remark}
\label{rem:ex12-elasticity-nonunique-ext}
Notice that, while the space $\mathcal{X}_0^*$ is two-dimensional, external forces $\bm{F}$ are taken from the three-dimensional subspace of $\mathbb{R}^4$, given by \eqref{eq:ex12-F-restriction-form2}. In other words, we again have the situation with non-injective mapping $\bm{F}\mapsto \bm{f}$, i.e. a functional $\bm{f}\in \mathcal{X}_0^*$ admits many extensions $\bm{F}\in \mathcal{X}^*$, which differ by a reaction term
\[
\bm{F_0}\in \mathcal{X}_0^\perp = ({\rm Ker}\, (R{\rm E}))^\perp ={\rm Im}\, (R{\rm E})^\top,
\]
compare this to Remarks \ref{rem:abstract-elasticity-nonunique-ext}, \ref{rem:abstract-operator-to-restrict} \ref{rem:abstract-operator-to-restrict-about-forces} and \ref{rem:ex11-elasticity-nonunique-ext}. Observe, that all such $\bm{F_0}$ satisfy \eqref{eq:ex12-F-restriction-form2}, thus \eqref{eq:ex12-F-restriction-form2} is an extension-independent condition on $\bm{f}$. One would expect this, as \eqref{eq:ex12-F-restriction-form2} it is really \eqref{eq:resolvable-load-def}, written for $\bm{F}$ in this particular example. Moreover, we can notice that 
\begin{equation}
{\rm Ker}\, {\bf E}={\rm Ker}\,{\rm E}\cap D({\bf E})={\rm Ker}\,{\rm E}\cap{\rm Ker}\,R{\rm E} = {\rm Im}\begin{pmatrix}1\\1\\1\\1\end{pmatrix},
\label{eq:ex2-kerE}
\end{equation}
while \eqref{eq:ex12-F-restriction-form2} means that
\[
\bm{F}\in \left({\rm Im}\begin{pmatrix}1\\1\\1\\1\end{pmatrix}\right)^\perp = ({\rm Ker}\, {\bf E})^\perp.
\]
Such representation of condition \eqref{eq:resolvable-load-def} is not limited to this particular example, as we shall discuss it in the general setting. 
\end{Remark}

\subsubsection{Resolvable loads is terms of kinematics in the abstract setting}
\label{sssect:on_rigid_motions}
Recall, that to construct the subspaces $\mathcal{U}$ and $\mathcal{V}$ and solve the problem of Definition \ref{def:ae} we exploited the orthogonality of ${\rm Im}\, {\bf E}={\rm Im}\, {\bf E'}$ and ${\rm Ker}\, {\bf D}$ in $\mathcal{H}$, which, respectively, describe kinematically and statically admissible ``intrinsic'' states, see Fig.~\ref{fig:elasticity-scheme}. Just above we have arrived to another connection between kinematics and statics, which is a consequence of the orthogonality between ${\rm Ker}\, {\bf E'}={\rm Ker}\, {\bf E}$ and ${\rm Im}\, {\bf D}$ with its memebrs' extensions to $\mathcal{X}$. Specifically, we mean the following fact.
\begin{Lemma}
\label{lemma:orthogonality-in-physical-space}
In the setting of Definition \ref{def:abstract-adjoint}
\[
({\rm Ker}\, {\bf E'})^\perp = {\rm Im}\, {\bf D}\subset \mathcal{X}_0^*
\]
and
\[
({\rm Ker}\, {\bf E})^\perp = \{\bm{F}\in \mathcal{X}^*:\text{the corresponding by \eqref{eq:ae-force-functional-abstract} }\bm{f}\in{\rm Im}\, {\bf D}\},
\]
where the orthogonality is in the sense of a Banach space \cite[p.~9]{Brezis2011}.
\end{Lemma}
\noindent{\bf Proof.} The equality in $\mathcal{X}^*_0$ is another direct application of \cite[Th.~2.19, p.~46]{Brezis2011}. Thus, if $\bm{F}\in X^*$ takes zero values on ${\rm Ker}\, {\bf E} = {\rm Ker}\, {\bf E'}$ then its restriction $\bm{f}$ must be in $ {\rm Im}\, {\bf D}$, and vice versa. $\blacksquare$

\noindent Lemma \ref{lemma:orthogonality-in-physical-space} means that an external force given as a functional $\bm{F}$ over the configurations in the physical space is resolvable (see \eqref{eq:resolvable-load-def} and Remark \ref{rem:abstract-operator-to-restrict} \ref{rem:abstract-operator-to-restrict-about-forces}, \ref{rem:abstract-operator-resolvable-forces-def} in general) if and only if 
\begin{equation}
\bm{F}\in ({\rm Ker}\, {\bf E})^\perp.
\label{eq:resolvable-load-kinematic-def}
\end{equation}
So far we have encountered two common situations:
\begin{itemize}
\item {\it Kinematic determinacy}, which means that ${\rm Ker}\, {\bf E}=\{0\}$, therefore any $\bm{F}\in \mathcal{X}^*$ is resolvable, as in Example 1. The model \eqref{eq:ae-1}--\eqref{eq:ae-3} ``registers'' any admissible motion $\bm{\delta u}\in D({\bf E})$ when it is kinematically determinate. It also means that we can uniquely recover displacement $\bm{\widetilde{u}}$ satisfying \eqref{eq:ae-compatibility-abstract}--\eqref{eq:ae-bc-abstract} after solving the problem for $\bm{\widetilde{\sigma}}$ and $\bm{\widetilde{\varepsilon}}$ via Theorem \ref{th:elasticity-solution} and Remark \ref{remark:to-find-elastic-strain} respectively.   The correspondence between \eqref{eq:resolvable-load-def} and \eqref{eq:resolvable-load-kinematic-def} for kinematically determinate discrete models is stated and discussed in \cite[Th.~3.3]{Gudoshnikov2023preprint}.
\item {\it Infinitesimal rigidity} means that ${\rm Ker}\, {\bf E}$ is the smallest we can hope for when the sample (e.g. a system of particles or a continuous body) is not anchored in the physical space $\mathbb{R}^d, d\in \{1,2,3\}$. Specifically, we call a model infinitesimally rigid when ${\rm Ker}\, {\bf E}$ consists of the configurations, generated from the reference configuration by the actions of infinitesimal rigid motions of the physical space (i.e. the members of Lie algebra of Euclidean group $E(d)$, see \cite[Sect.~3.8.4.2, pp.~216--222]{Ivancevic2007}, \cite[Sect.~4.1, pp.~48--50]{Jensen2016}). 

In the current Example 2 we have $\mathbb{R}^1$ as the physical space (see Fig.~\ref{fig:discrete-models-elastic} b), so its rigid motions are just translations along a single direction. Applied to the four nodes, these translations generate exactly the space \eqref{eq:ex2-kerE}, i.e. the model of Example 2 is infinitesimally rigid. As we noted previously, the corresponding criterion on resolvable loads \eqref{eq:ex12-F-restriction-form2} merely means that the total applied momentum is zero. 

More generally, when $d\in\{2,3\}$ and the model is infinitesimally rigid, the resolvability criterion \eqref{eq:resolvable-load-kinematic-def} means that both total applied momentum and total applied torque is zero. For the discussions about discrete unanchored models and their infinitesimal rigidity see \cite[Ch.~9]{Alfakih2018_rigidity}, \cite{RothWhiteley1981}, \cite[Sect.~3.2,~3.5]{Gudoshnikov2023preprint}.
\end{itemize}
Naturally, kinematic determinacy and infinitesimal rigidity are not all possible situations, as one can consider an example where only some, but not all of the rigid motions in $\mathbb{R}^d, d\in\{2,3\}$ are restrained by \eqref{eq:ae-bc-abstract}, or  an example where one part of the sample can move relative to another without changing the value of ${\bf E}$ (this is called a {\it mechanism/floppy mode}, see e.g. \cite{Lubensky2015}).

\subsubsection{Example 2 --- the fundamental subspaces and the stress solution}
Coming back to Example 2, above we have provided all of the data required in Definition \ref{def:ae} to set up a problem of the type \eqref{eq:ae-1}--\eqref{eq:ae-3}. Now we will give exact expressions of the quantities from Definition \ref{def:ae-weighted} and Theorem \ref{th:elasticity-solution} for the model of Fig. \ref{fig:discrete-models-elastic} {\it b)}.

Space $\mathcal{H}_{{\bf C}^{-1}}$ is $\mathbb{R}^3$ equipped with the inner product
\begin{equation}
\left\la \begin{pmatrix} \sigma_1\\ \sigma_2\\ \sigma_3\end{pmatrix}, \begin{pmatrix}\tau_1\\ \tau_2\\ \tau_3\end{pmatrix} \right\ra_{{\bf C}^{-1}} = \sigma_1 k_1^{-1} \tau_1+ \sigma_2 k_2^{-1} \tau_2+ \sigma_3 k_3^{-1} \tau_3 \qquad \text{for all } \begin{pmatrix} \sigma_1\\ \sigma_2\\ \sigma_3\end{pmatrix}, \begin{pmatrix}\tau_1\\ \tau_2\\ \tau_3\end{pmatrix}\in \mathbb{R}^3.
\label{eq:ex12-inner-product}
\end{equation}
Subspace  $\mathcal{U}$ is
\begin{multline*}
\mathcal{U}  ={\rm Im}\,{\bf C} {\rm E} R_0  = {\rm Im}\, \begin{pmatrix}k_1 & 0 & 0\\ 0 & k_2 & 0 \\ 0 & 0 & k_3\end{pmatrix}   \begin{pmatrix}
-1 & 1 & 0 & 0\\
0 & -1 & 1 & 0\\
0 & 0 & -1 & 1
\end{pmatrix} \begin{pmatrix}1& 0\\ 0 & 1\\ 1& 0 \\ 0 & 1\end{pmatrix}  =\\= {\rm Im}\,  \begin{pmatrix} -k_1 & k_1 \\ k_2 & -k_2 \\ -k_3 & k_3  \end{pmatrix}= {\rm Im}\, \begin{pmatrix} -k_1 \\ k_2 \\ -k_3  \end{pmatrix},
\end{multline*}
and, from \eqref{eq:ex12-explicit-D} used in the definition of ${\bf D}_{\bf C}$ in \eqref{eq:EC-DC-def}, we have that
\begin{equation}
\mathcal{V} = {\rm Ker}\, \begin{pmatrix} -1 & 1 & -1\end{pmatrix} = {\rm Im}\,  \begin{pmatrix} 1 & 0 \\1 & 1 \\ 0  & 1\end{pmatrix},
\label{eq:ex12-space-V-exp}
\end{equation}
see Fig. \ref{fig:discrete-models-spaces} {\it b)}.

Again, by the formula \cite[(5.31), p. 247]{OlverLinearAlgebra2018} we have the matrices of projections
\[
P_{\mathcal{U}} 
= \frac{1}{k_1+k_2+k_3}\begin{pmatrix} k_1 &  -k_1 & k_1 \\ -k_2 & k_2 & -k_2\\ k_3 & -k_3 & k_3\end{pmatrix},
\]
\[
P_{\mathcal{V}} = \frac{1}{k_1^{-1}k_2^{-1}+k_2^{-1}k_3^{-1}+k_1^{-1}k_3^{-1}} \begin{pmatrix} k_1^{-1}k_2^{-1}+k_1^{-1}k_3^{-1} & k_2^{-1}k_3^{-1}& -k_2^{-1}k_3^{-1} \\[1mm] k_1^{-1}k_3^{-1} & k_1^{-1}k_2^{-1}+k_2^{-1}k_3^{-1} & k_1^{-1}k_3^{-1} \\[1mm]-k_1^{-1}k_2^{-1}& k_1^{-1}k_2^{-1}&k_1^{-1}k_3^{-1}+ k_2^{-1}k_3^{-1} \end{pmatrix}.
\]
Further, operator $G$ is represented by the matrix
\[
G= P_{\mathcal{V}}{\bf C} = \frac{1}{k_1^{-1}k_2^{-1}+k_2^{-1}k_3^{-1}+k_1^{-1}k_3^{-1}}\begin{pmatrix} k_2^{-1}+k_3^{-1} & k_3^{-1}& -k_2^{-1} \\[1mm]k_3^{-1} & k_1^{-1}+k_3^{-1} & k_1^{-1} \\[1mm]-k_2^{-1}& k_1^{-1}&k_1^{-1}+ k_2^{-1} \end{pmatrix}.
\]

We now find a representation for the operator $Q$ in terms of an  external force $\bm{F}$. Notice, that for any $\bm{\sigma}\in \mathcal{U}$  there exists a unique $\sigma\in \mathbb{R}$, such that 
\[
\bm{\sigma}= \begin{pmatrix}-k_1\\ k_2\\-k_3\end{pmatrix} \sigma. 
\]
Plug that into  \eqref{eq:ex12-explicit-D}, and then plug the result together with \eqref{eq:ex12-functionals-explicit} into \eqref{eq:ae-3}  to see that
\[
(k_1+k_2+k_3)\,\sigma \,(u_1-u_2) = \begin{pmatrix}F_1+ F_3& F_2+ F_4\end{pmatrix}\begin{pmatrix}u_1\\ u_2\end{pmatrix}.
\]
Due to \eqref{eq:ex12-F-restriction-form1} we get that
\[
(k_1+k_2+k_3)\sigma  = F_1+F_3,
\]
i. e.
\[
\sigma = \frac{F_1+F_3}{k_1+k_2+k_3}.
\]
Therefore, for a given external force $\bm{F}\in \mathbb{R}^{4}$ and the corresponding (by \eqref{eq:ex12-functionals-def}) functional $\bm{f}\in \mathcal{X}_0^*$ we have the value of the operator $Q$ 
\[
Q\, \bm{f} = \frac{F_1+F_3}{k_1+k_2+k_3}\begin{pmatrix}-k_1 \\ k_2\\ -k_3\end{pmatrix}.
\]

By Theorem \ref{th:elasticity-solution} we conclude that the stress solution for elasticity in the model of Fig.  \ref{fig:discrete-models-elastic} {\it b)} exists whenever \eqref{eq:ex12-F-restriction-form2} holds and it has the closed form
\begin{multline}
\bm{\widetilde{\sigma}}(t) = G\, \bm{g}(t) + Q\, \bm{f}(t) = \\ 
=\frac{1}{k_1^{-1}k_2^{-1}+k_2^{-1}k_3^{-1}+k_1^{-1}k_3^{-1}}\begin{pmatrix} k_2^{-1}+k_3^{-1} & -k_2^{-1} \\[1mm]k_3^{-1} & k_1^{-1} \\[1mm]-k_2^{-1}& k_1^{-1}+ k_2^{-1} \end{pmatrix}\begin{pmatrix}  l_1(t) \\  l_2(t)\end{pmatrix} + \frac{F_1(t)+F_3(t)}{k_1+k_2+k_3}\begin{pmatrix}-k_1 \\ k_2\\ -k_3\end{pmatrix},
\label{eq:ex12-linear-solution}
\end{multline}
see Fig. \ref{fig:discrete-models-spaces} {\it b)}.

\subsection{Example 3 --- a continuum model}
\label{ssect:ex21-elasticity}
\begin{figure}[H]\center
\includegraphics{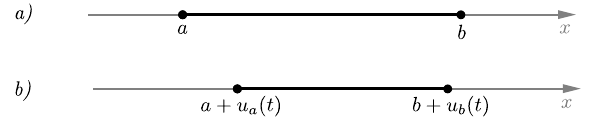}
\caption{\footnotesize A continuum model of Example 3 in the relaxed reference configuration ({\it a}) and an arbitrary current configuration ({\it b}). 
} 
\label{fig:continuous-rod-model-elastic}
\end{figure} 
As the third example we describe the quasi-static evolution of a one-dimensional elastic rod subjected to Dirichlet boundary conditions at both endpoints and a body force load distributed along its length. We use the Lagrangian description with respect to a relaxed reference configuration of the rod occupying the interval $\overline{\Omega}$ with 
\[
\Omega = (a,b)
\] for given reference values $a<b$.
\subsubsection{Example 3 --- kinematics}
The displacement variable is now an unknown function 
\[
\widetilde{u}\in W^{1,\infty}(I, W^{1,2}(\Omega))
\]
of time $t\in  I$  and Lagrangian coordinate $x\in \Omega$.
Let the Dirichlet boundary conditions at $x=a$ and $x=b$ be given by scalars $u_a(t), u_b(t)\in \mathbb{R}$ respectively. We assume that 
\begin{equation}
u_a, u_b\in W^{1, \infty}(I)
\label{eq:ex21-bc-1}
\end{equation}
and that $a+u_a(t)<b+u_b(t)$ for all $t>0$.
Using $u_a(t)$ and $u_b(t)$ we choose a function
\begin{equation}
u_{\rm D} \in W^{1, \infty}(I, W^{1,2}(\Omega)), 
\label{eq:ex21-bc-2}
\end{equation}
such that 
\begin{equation}
u_{\rm D}(t,a)= u_a (t), \qquad u_{\rm D}(t,b)= u_b (t), \qquad \text{for all }t\in I,
\label{eq:ex21-bc-3}
\end{equation}
e. g. we can choose the linear interpolation
\[
u_{\rm D}(t,x) = \frac{b-x}{b-a}\,u_a(t)+\frac{x-a}{b-a}\,u_b(t).
\]
Using such $u_{\rm D}$ we write the Dirichlet boundary condition as
\begin{equation}
\bm{\widetilde{u}}(t) \in W_0^{1,2}(\Omega) +{\bm{u}_{\bf D}}(t).
\label{eq:ex21-bc}
\end{equation}

Define the following operator
\begin{align}
&{\rm E}: W^{1,2}(\Omega) \to L^2(\Omega),
\label{eq:strain-operator-def1}\\
&{\rm E}: \bm{u}\mapsto \frac{d}{d x}\bm{u} \label{eq:strain-operator-def2}
\end{align}
and call {\it strain} an unknown function 
\[
\widetilde{\varepsilon}\in W^{1, \infty}(I, L^{2}(\Omega)),
\]
 defined for each $t\in I$ by formula \eqref{eq:ae-compatibility-abstract}.

In terms of the framework of Definition \ref{def:ae} and Remark \ref{rem:abstract-operator-to-restrict} we set
\[
\mathcal{H} =  L^2(\Omega) \qquad \mathcal{X}=\mathcal{X}_0=L^2(\Omega),
\]
\begin{equation}
\mathcal{W} = D({\rm E})=W^{1,2}(\Omega), \qquad \mathcal{W}_0 = D({\bf E}) = W^{1,2}_0(\Omega),
\label{eq:ex21-E-domain-def}
\end{equation}
and we define operator {\bf E} by the formula \eqref{eq:ex11-E-def}. 
We follow \eqref{eq:ae-abstract-g-def} and set
\begin{equation}
\bm{g}(t) = {\rm E}\, {\bm{u}_{\bf D}}(t) = \frac{\partial}{\partial x}{u_ {\rm D}}(t,x),
\label{eq:ex21-g-def}
\end{equation}
so that we have set up \eqref{eq:ae-compatibility-abstract}, \eqref{eq:ae-bc-abstract} (which is \eqref{eq:ex21-bc}), and \eqref{eq:ae-1}. It is a standard exercise in functional analysis to check that the operator ${\bf E}$, as we defined it here, is unbounded, densely defined, closed and has closed image. Clearly 
\begin{align*}
{\rm Ker}\, {\rm E} =& \{u\in W^{1,2}(\Omega):\text{there exists } c\in\mathbb{R} \text{ such that for all }x\in \Omega \text{ we have } u(x)\equiv c\},\\[2mm]
{\rm Ker}\, {\bf E} = &{\rm Ker}\, {\rm E} \cap W_0^{1,2}(\Omega)= \{0\},
\end{align*}
i.e. the continuum rod model with Dirichlet boundary condition is kinematically determinate in terms of Section \ref{sssect:on_rigid_motions}. In higher dimensions a continuous body occupying an open set $\Omega$ is kinematically determinate when the Dirichlet boundary condition is enforced on a part of the boundary of $\Omega$ of positive measure (on the boundary), see e.g. \cite[Lemma~2.1, p.~ 91]{Necas1981}, \cite[Prop.~3.3 {\it ii)}, p.~12]{Gudoshnikov2025}. In turn, a model with an open set $\Omega$ representing an unanchored body with pure traction boundary condition is infinitesimally rigid, see \cite[Sect.~7.3, pp.~95--100]{Necas1981}, \cite{Ciarlet2005}. However, an example with several continuous bodies joined together via hinges could have floppy modes.

\subsubsection{Example 3 --- Hooke's law}
\label{sssect:ex21-Hooke-law}
The elastic property of the rod is characterized by the {\it stiffness}, which is a given function $\bm{C}\in L^{\infty}(\Omega)$, for which there exist $0<c_0\leqslant C_0$ such that  
\[
c_0\leqslant \bm{C}(x)\leqslant C_0 \qquad \text{for a.a. }x\in \Omega.
\]
We define operator ${\bf C}$, acting as \eqref{eq:ae-C-def}, by the formula
\[
{\bf C}: \bm{\widetilde{\varepsilon}} \mapsto \bm{C}\bm{\widetilde{\varepsilon}},
\]
and then \eqref{eq:ae-2} is the constitutive law of elasticity, which connects strain $\widetilde{\varepsilon}$ and {\it stress} $\widetilde{\sigma}\in W^{1,\infty}(I,L^2(\Omega))$ for each $t\in I$.
\subsubsection{Example 3 --- statics}
Consider given {\it body force load} 
\begin{equation}
F\in W^{1, \infty}(I, L^{2}(\Omega)),
\label{eq:ex21-body-force-load}
\end{equation}
applied at each $x\in \Omega$ along the direction of the rod.
 The principle of virtual work (for the continuum mechanics formulations see \cite[Section 2.2.2, pp.~63--65]{Altenbach2018}, \cite[Appendix 1, pp.~589--591]{Anandarajah2011}, \cite[Remark 5.1.1, p.~61]{Necas1981}) yields the following equation of equilibrium of the rod
\begin{equation}
-\intOm \bm{\widetilde{\sigma}}(t)\, {\rm E}(\delta u)\, dx + \intOm \bm{F}(t)\, \delta u\, dx = 0 \qquad \text{for all }\delta u\in W^{1,2}_0(\Omega).
\label{eq:ex21-equilibrium-basic-formulation}
\end{equation}

On the other hand, when ${\bf E}$ is defined via \eqref{eq:ae-abstract-E-def}, \eqref{eq:ex11-E-def}, \eqref{eq:ex21-E-domain-def}, its adjoint ${\bf D} ={\bf E}^*$ is defined by the formula
\begin{equation}
{\bf D}: \bm{\sigma} \mapsto -\frac{d}{dx} \bm{\sigma}
\label{eq:ex21-D-def}
\end{equation}
for any $\bm{\sigma}$ from the domain 
\begin{equation}
D({\bf D}) = W^{1,2}(\Omega).
\label{eq:ex21-D-domain-def}
\end{equation}
To prove the adjacency, and to check that ${\bf D}$ is unbounded, densely defined, closed and has closed image is again a standard exercise in functional analysis (we also refer to \cite[Th. 4.1]{Gudoshnikov2025}, which is a similar statement for a more complicated case of many spatial dimensions). Notice, in particular, that the operator ${\bf D}$ is surjective i.e. every $\bm{F}\in L^2(\Omega)$ is resolvable, as one would expect for a kinematically determinate model.
\begin{Remark}
Naturally, we use the Riesz representation theorem (see e.g. \cite[Th. 5.5, p. 135]{Brezis2011}) to understand the values in \eqref{eq:ex21-D-def} as functionals over $\mathcal{W}_0 = W^{1,2}_0(\Omega)$ and, by extension, over $\mathcal{X}=\mathcal{X}_0=L^{2}(\Omega)$. In particular, due to the density of $W_0^{1,2}(\Omega)$ in $L^2(\Omega)$ we have this time  $\mathcal{X}=\mathcal{X}_0$, cf. Remarks \ref{rem:abstract-elasticity-nonunique-ext}, \ref{rem:abstract-operator-to-restrict} \ref{rem:abstract-operator-to-restrict-about-forces}, \ref{rem:ex11-elasticity-nonunique-ext}, \ref{rem:ex12-elasticity-nonunique-ext}. This also means that the representation of the reaction forces to the boundary condition as members of $\mathcal{X}_0^\perp$ is no longer useful, as we have
\[
\mathcal{X}_0^\perp =\{0\}.
\]
\end{Remark}

For each $t\in I$ define $\bm{f}(t)\in \mathcal{H^*}$ by $\bm{F}(t)$ in accordance with the Riesz representation theorem, i.e.
\[
\bm{f}: \bm{\delta u} \mapsto \intOm \bm{F}\, \bm{\delta u} \,dx \qquad \text{for any }\bm{\delta u}\in \mathcal{X}_0=\mathcal{X}= L^2(\Omega).
\]
By plugging this and the formula \eqref{eq:ex11-E-def} into the equation of equilibrium \eqref{eq:ex21-equilibrium-basic-formulation} we get exactly \eqref{eq:ae-3}.
\subsubsection{Example 3 --- the fundamental subspaces and the stress solution}
In the continuum model we described above we have the space $\mathcal{H}_{\bf C^{-1}}$ being $L^2(\Omega)$, equipped with the inner product
\begin{equation}
\la\bm{\tau}_1, \bm{\tau}_2\ra_{{\bf C}^{-1}} = \intOm \bm{\tau}_1\, {\bf C}^{-1}\,\bm{\tau}_2\, dx =\intOm \frac{\tau_1(x)\, \tau_2(x)}{\bm{C}(x)}\, dx.
\label{eq:ex21-ip}
\end{equation}
We denote the space by $L^2_{{\bf C}^{-1}}(\Omega)$.

From the definition of ${\mathcal{V}}$ in \eqref{eq:ae-UV-def}, the definition of ${\bf D_C}$ in \eqref{eq:EC-DC-def}, and the formulas \eqref{eq:ex21-D-def}--\eqref{eq:ex21-D-domain-def} for ${\bf D}$ we deduce with the help of the fundamental theorem of calculus for Lebesgue integrals (see e.g. \cite[Section 3.2, pp.~71--84]{Leoni2017}) that the subspace $\mathcal{V}$ is the one-dimensional space of constant functions:
\begin{equation}
\mathcal{V} = \left\{\bm{\sigma} \in L^2_{{\bf C}^{-1}}(\Omega): \text{ there exists }c\in \mathbb{R} \text{ such that for a.a. }x\in \Omega \text{ we have } \sigma(x)\equiv c \right\}.
\label{eq:ex21-space-V}
\end{equation}

Therefore, it orthogonal complement in the sense of the inner product \eqref{eq:ex21-ip} is the subspace of functions with zero weighted average:
\[
\mathcal{U} = \left\{\bm{\sigma}\in L^2_{{\bf C}^{-1}}(\Omega): \intOm \frac{1}{\bm{C}(x)}\,\sigma(x)\, dx=0\right\}.
\]

To find the formulas for the operators of projection  \eqref{eq:ae-PU-PV-def} 
we take arbitrary $\bm{\sigma}\in L^2_{{\bf C}^{-1}}(\Omega)$ and solve for $c\in \mathbb{R}$ such that
\[
\intOm \frac{\sigma(x)\lambda}{\bm{C}(x)} dx= \intOm \frac{c \lambda}{\bm{C}(x)} dx\qquad \text{for any }\lambda\in \mathbb{R},
\]
from which it follows that
\[
c=\frac{\intOm \frac{\sigma(x)}{\bm{C}(x)} dx}{\intOm \frac{1}{\bm{C}(x)} dx}.
\]
Thus the operators of projection are
\[
(P_{\mathcal{V}}\, \bm{\sigma})(x) \equiv w_0 \intOm \frac{1}{\bm{C}(y)}\,\sigma(y)\, dy, \qquad (P_{\mathcal{U}}\, \bm{\sigma})(x) = \sigma(x)-w_0 \intOm \frac{1}{\bm{C}(y)}\,\sigma(y)\, dy,
\] 
where $w_0\in \mathbb{R}$ is the constant
\begin{equation}
w_0 = \frac{1}{\intOm\frac{1}{\bm{C}(y)}\, dy}.
\label{eq:w0-const-def}
\end{equation}
It follows that the operator $G$ from \eqref{eq:ae-G-def} is 
\begin{equation}
(G\, \bm{\varepsilon})(x) \equiv w_0 \intOm \varepsilon(y)\, dy \qquad \text{for any }\bm{\varepsilon}\in \mathcal{H}  =L^2(\Omega).
\label{eq:ex21-operator-G-def}
\end{equation}

To find the representation on the operator $Q$ in terms of a given $\bm{F}\in L^2(\Omega)$ and the corresponding $\bm{f}\in \mathcal{X}_0^*$ we must find $\bm{\sigma}\in \mathcal{U}$ satisfying \eqref{eq:ae-3} (since the values of ${\bf D}$ and ${\bf D_C}$ coincide by construction, see \eqref{eq:EC-DC-def}). This means that
\[
\intOm \frac{1}{\bm{C}(x)}\,\sigma(x)\, dx=0
\]
and, by \eqref{eq:ae-3}, \eqref{eq:ex21-D-def} and the fundamental theorem of calculus for Lebesgue integrals, we have that
\[
\sigma(x) - \sigma(a) = -\int\limits_a^x F(x)\, dx\qquad\text{for any } x\in \Omega=(a,b),
\]
i.e.
\begin{equation}
\sigma(x) =\sigma(a) -\int\limits_a^x F(y)\, dy.
\label{eq:ex21-stress-intermediate-expression}
\end{equation}
Therefore,
\[
\intOm \frac{1}{\bm{C}(x)}\left(\sigma(a) -\int\limits_a^x F(y)\, dy\right) dx=0,
\]
i.e.
\[
\sigma(a)\intOm \frac{1}{\bm{C}(x)}\,dx= \intOm \frac{1}{\bm{C}(x)}\int\limits_a^x F(y)\, dy\, dx
\]
and, if we rename $x$ to $z$, we can compute
\[
\sigma(a)  = \omega_0 \intOm \frac{1}{\bm{C}(z)}\int\limits_a^z F(y)\, dy\, dz.
\]
Thus, using \eqref{eq:ex21-stress-intermediate-expression}, we can write that
\begin{equation}
(Q\bm{f})(x)= \sigma(x)  =\omega_0 \intOm \frac{1}{\bm{C}(z)}\int\limits_a^z F(y)\, dy\, dz  -\int\limits_a^x F(y)\, dy.
\label{eq:ex21-operator-Q-def}
\end{equation}

We plug \eqref{eq:ex21-g-def}, \eqref{eq:ex21-operator-G-def} and \eqref{eq:ex21-operator-Q-def} into \eqref{eq:abstract_linear_solution} to write the stress solution for elasticity in terms of given time-dependent displacement boundary data \eqref{eq:ex21-bc-1}--\eqref{eq:ex21-bc-3} and a time-dependent longitudinal body force load \eqref{eq:ex21-body-force-load} in the closed form
\begin{equation}
\begin{aligned}
\widetilde{\sigma}(t,x)& = \left(G\, \bm{g}(t)\right)(x) + \left(Q\, \bm{f}(t)\right)(x)=\\[1mm] &= \omega_0 \intOm \frac{d}{dy} u_{\rm D}(t,y)\, dy + \omega_0 \intOm \frac{1}{\bm{C}(z)}\int\limits_a^z F(t,y)\, dy\, dz  -\int\limits_a^x F(t,y)\, dy = \\[1mm] &= \omega_0 \left(u_b (t)- u_a (t)+ \intOm \frac{1}{\bm{C}(z)}\int\limits_a^z F(t,y)\, dy\, dz\right) - \int\limits_a^x F(t,y)\, dy.
\end{aligned}
\label{eq:ex21-linear-solution}
\end{equation}
For a reader's convenience we collect the quantities from all three examples of linear elasticity in Table \ref{tab:linear_elasticity_examples}.

\begin{table}[H]
\begin{center}
\renewcommand{\arraystretch}{1.5} 
\begin{tabular}{|c|c|c|c|}
\hline
\textbf{Quantity} & \textbf{Example 1 (Fig. \ref{fig:discrete-models-elastic} {\it a})} & \textbf{Example 2 (Fig. \ref{fig:discrete-models-elastic} {\it b})} & \textbf{Example 3 (Fig. \ref{fig:continuous-rod-model-elastic})}\\
\hline
$\mathcal{X}$ & $\mathbb{R}^3$ & $\mathbb{R}^4$ &  $L^2(\Omega), \, \Omega=(a,b)$
\\
\hline
$\mathcal{H}$ & $\mathbb{R}^2$ & $\mathbb{R}^3$ & $L^2(\Omega)$\\
\hline
$\mathcal{W}$ & $\mathbb{R}^3$ & $\mathbb{R}^4$ &$W^{1,2}(\Omega)$ \\
\hline
${\rm E}:\mathcal{W}\to \mathcal{H}$ & \renewcommand{\arraystretch}{1} $\begin{pmatrix}
-1 & 1 & 0 \\
0 & -1 & 1
\end{pmatrix}$ & \begin{minipage}{31mm}\begin{center}\renewcommand{\arraystretch}{1} $\begin{pmatrix}
-1 & 1 & 0 & 0\\
0 & -1 & 1 & 0\\
0 & 0 & -1 & 1
\end{pmatrix}$ \end{center}\end{minipage} &  $ \bm{u}\mapsto \frac{d}{d x}\bm{u}$
\\ \hline
$\mathcal{W}_0$ & $\begin{array}{c} {\rm Ker}\, R, \\ R={\renewcommand{\arraystretch}{1}\begin{pmatrix}
1 & 0 & 0\\
0 & 0 & 1
\end{pmatrix}}
\end{array}
$  & 
$\begin{array}{c} {\rm Ker}\, ({R} {\rm E}),\\ {R} ={\renewcommand{\arraystretch}{1} \begin{pmatrix}
1 & 1 & 0\\
0 & 1 & 1
\end{pmatrix}}
\end{array}
$
& $W_0^{1,2}(\Omega)$ \\ 
\hline
$\mathcal{X}_0 =\overline{\mathcal{W}_0}$ & ${\rm Ker}\, R$ & ${\rm Ker}\, ({R} {\rm E})$ & $L^2(\Omega)$\\ \hline
``Reactions space'' $\mathcal{X}_0^\perp$ & ${\rm Im}\,R^\top$ & ${\rm Im}\,(R{\rm E})^\top$ & $\{0\}$ \\      
\hline
${\bf E}: \mathcal{W}_0 \subset \mathcal{X} \to  \mathcal{H}$ & ${\rm E} \text{ restricted to } \mathcal{W}_0$ & ${\rm E} \text{ restricted to } \mathcal{W}_0$ & ${\rm E} \text{ restricted to } \mathcal{W}_0$ \\ 
\hline
${\bf D}: D({\bf D})\subset \mathcal{H}\to \mathcal{X}_0^* $ & $\begin{array}{c} \bm{\sigma} \mapsto \text{the functional}\\\left(\bm{u}\in \mathcal{X}_0\mapsto \left({\rm E}^\top \bm{\sigma}\right)\cdot \bm{u}\right)\end{array}$ & $\begin{array}{c} \bm{\sigma} \mapsto \text{the functional}\\\left(\bm{u}\in \mathcal{X}_0\mapsto \left({\rm E}^\top \bm{\sigma}\right)\cdot \bm{u}\right)\end{array}$ & $\bm{u}\mapsto -\frac{d}{d x}\bm{u}$
\\ 
\hline
$D({\bf D})$ & $\mathbb{R}^2$ & $\mathbb{R}^3$ & $W^{1,2}(\Omega)$\\
\hline
${\bf C}:\mathcal{H}\to \mathcal{H}$ & \renewcommand{\arraystretch}{1} $\begin{pmatrix} 
k_1 & 0\\[1mm] 0 &  k_2
\end{pmatrix}$ & \begin{minipage}{23mm}\begin{center}\renewcommand{\arraystretch}{1} $\begin{pmatrix} 
k_1 & 0 & 0\\[1mm] 0 &  k_2 & 0\\[1mm] 0 & 0 & k_2
\end{pmatrix}$ \end{center}\end{minipage}  & $\bm{C}\in L^\infty(\Omega)$
\\ 
\hline
 $\begin{array}{c}\mathcal{U}={\rm Im}\, {\bf E_C}=\\={\rm Im}\, {\bf C E}\subset \mathcal{H}_{\bf C^{-1}}\end{array}$& \renewcommand{\arraystretch}{1} ${\rm Im}\,  \begin{pmatrix} k_1 \\ -k_2  \end{pmatrix}$ & \renewcommand{\arraystretch}{1} $ {\rm Im}\, \begin{pmatrix} -k_1 \\ k_2 \\ -k_3  \end{pmatrix}$ &\begin{minipage}{34mm}\begin{center} $\begin{array}{c}\left\{\bm{\sigma}\in L^2_{{\bf C}^{-1}}(\Omega):\vphantom{\intOm}\right.\\ \left. \intOm \frac{1}{\bm{C}(x)}\,\sigma(x)\, dx=0\right\}\end{array}$ \end{center}\end{minipage}
\\
\hline
$\begin{array}{c}\mathcal{V}={\rm Ker}\, {\bf D_C}=\\
={\rm Ker}\, {\bf D}\subset \mathcal{H}_{\bf C^{-1}}\end{array} $&\renewcommand{\arraystretch}{1} ${\rm Im}\, \begin{pmatrix} 1 \\ 1\end{pmatrix}$ &\begin{minipage}{18mm}\begin{center} \renewcommand{\arraystretch}{1} $ {\rm Im}\,  \begin{pmatrix} 1 & 0 \\1 & 1 \\ 0  & 1\end{pmatrix}$ \end{center}\end{minipage} & $\begin{array}{c}\left\{\bm{\sigma} \in L^2_{{\bf C}^{-1}}(\Omega): \exists c\in \mathbb{R}\right.\\\left. \text{for a.a. }x\in \Omega \quad \sigma(x)\equiv c \right\}\end{array}$
\\
\hline
$\bm{g}(t)$&\renewcommand{\arraystretch}{1} $\begin{pmatrix}  0 \\  l(t)\end{pmatrix}$ & \begin{minipage}{33mm}\begin{center} \renewcommand{\arraystretch}{1} $\frac{1}{2}\begin{pmatrix}1 & -1 \\ 1 & 1 \\ -1 & 1\end{pmatrix} \begin{pmatrix}  l_1(t) \\  l_2(t)\end{pmatrix}$ \end{center}\end{minipage} &$ \frac{\partial}{\partial x}{u_ {\rm D}}(t,x)$
\\
\hline
$\begin{array}{c}\text{Stress solution }\\
\bm{\widetilde{\sigma}}(t)=  G\, {\bm g}(t) + Q\,{\bm f}(t) \end{array}$ &\text{formula \eqref{eq:ex11-linear-solution}}& \text{formula \eqref{eq:ex12-linear-solution}} & \text{formula \eqref{eq:ex21-linear-solution}}.
\\
\hline
\end{tabular}
\end{center}
\caption{\footnotesize The quantities in the examples of linear elasticity of the current paper.}  
\label{tab:linear_elasticity_examples}
\end{table}
\newpage
\section{The sweeping process framework for elasticity-perfect plasticity}
\label{sect:perfect-plasticity}
\subsection{The abstract geometric framework and the abstract sweeping process}
Now we consider the quasistatic evolution with elasticity and perfect plasticity combined, see Definition \ref{def:aepp} and Fig. \ref{fig:elasticity-pp-scheme} below. 
\begin{Definition} 
\label{def:aepp}
Let spaces $\mathcal{H}, \mathcal{X}, \mathcal{W}_0$, operators ${\bf E}, {\bf D}, {\bf C}$ and functions $\bm{g}, \bm{f}$ be as in Section \ref{ssect:adjoint-operatorsED} and Definition \ref{def:ae}, but we require $\mathcal{H}$ to be separable. In addition, let us be given a bounded closed convex nonempty set $\Sigma\subset \mathcal{H}$. We say, that the unknown variables
\begin{equation}
\varepsilon, \varepsilon_{\rm el}, \varepsilon_{\rm p}, \sigma \in W^{1, \infty}(I, \mathcal{H}) 
\label{eq:aepp-unknowns}
\end{equation}
solve the {\it abstract problem of quasi-static evolution in elasticity-perfect plasticity} if they satisfy
\begin{align}
\bm{\varepsilon}& \in {\rm Im}\, {\bf E}+\bm{g} , \label{eq:aepp-1} \tag{EPP1}\\
\bm{\varepsilon}& = \bm{\varepsilon_{\bf el}}+\bm{\varepsilon_{\bf p}}, \label{eq:aepp-2} \tag{EPP2}\\
\bm{\sigma}& = {\bf C} \, \bm{\varepsilon_{\bf el}}, \label{eq:aepp-3} \tag{EPP3}\\
\frac{d}{dt} \bm{\varepsilon_{\bf p}}& \in N_\Sigma (\bm{\sigma}), \label{eq:aepp-4} \tag{EPP4}\\ 
\bm{\sigma}\in D(\bf D)\quad \text{ and }\quad {\bf D}\, \bm{\sigma}&=\bm{f}, \label{eq:aepp-5} \tag{EPP5}
\end{align}
and the initial condition
\[
(\bm{\varepsilon}(0), \bm{\varepsilon}_{\bf  el}(0), \bm{\varepsilon}_{\bf  p}(0), \bm{\sigma}(0)) =(\bm{\varepsilon}_{\bf 0}, \bm{\varepsilon}_{\bf  el0}, \bm{\varepsilon}_{\bf  p0}, \bm{\sigma}_{\bf 0})
\]
with some given right-hand side from $\mathcal{H}^4$ satisfying \eqref{eq:aepp-1}--\eqref{eq:aepp-3}, \eqref{eq:aepp-5} at $t=0$ and 
\begin{equation}
{\bm \sigma}_{\bf 0} \in \Sigma.
\label{eq:aepp-sigma-ic-compatible}
\end{equation}
Equations \eqref{eq:aepp-1}--\eqref{eq:aepp-3}, \eqref{eq:aepp-5} are required to hold for all $t\in I$. Equation \eqref{eq:aepp-4} is required to hold for a.a. $t\in I$, as its left-hand side is the following vector-valued derivative (see Section \ref{ssect:prelim-bochner-sobolev})
\[
\frac{d}{dt} \varepsilon_{\rm p}\in L^\infty(I,\mathcal{H}).
\]
\noindent In turn, the right-hand side of \eqref{eq:aepp-4} is the outward normal cone to $\Sigma$ at $\bm{\sigma}$, as defined by formula \eqref{eq:nc-abstract-def}.
\end{Definition}

\begin{figure}[H]\center
\includegraphics{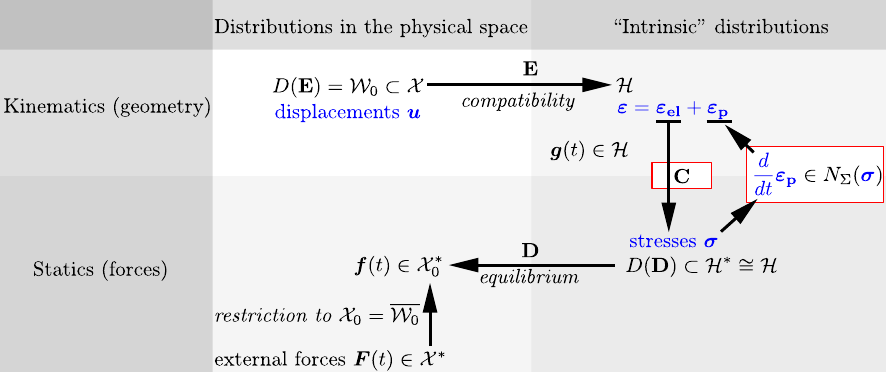}
\caption{\footnotesize Schematic representation of the problem of Definition \ref{def:aepp}. The unknown variables are indicated by blue color. In the problem of Definition \ref{def:aepp} we are only looking for the unknowns $\varepsilon, \varepsilon_{\rm el}, \varepsilon_{\rm p}$ and $\sigma$. Red rectangles indicate the constitutive relations.
} 
\label{fig:elasticity-pp-scheme}
\end{figure} 

We are going to write the abstract problem of Definition \ref{def:aepp} in terms of differential inclusions. In the current setting we have all the data, required by Definition \ref{def:ae-weighted}, so that we have the space $\mathcal{H}_{{\bf C}^{-1}}$, the operators ${\bf E_C}, {\bf D_C}$ and the fundamental subspaces $\mathcal{U}, \mathcal{V}$ defined accordingly.  Moreover, from Proposition \ref{prop:elasticity-affine-planes} we have the operators $G, Q$, and from Theorem \ref{th:elasticity-solution} we have the stress solution for elasticity $\widetilde{\sigma}$, which we will treat further as a {\it known} function of time, defined by \eqref{eq:abstract_linear_solution}.

In space $\mathcal{H}_{{\bf C}^{-1}}$ consider the outward normal cone \eqref{eq:nc-abstract-def} in the sense of the inner product \eqref{eq:weighted-ip}, which, for a closed convex nonempty set ${\mathcal C}\subset\mathcal{H}_{{\bf C}^{-1}}$ and a point $\bm{\tau}\in \mathcal{C}$, is
\begin{equation}
N^{{\bf C}^{-1}}_{\mathcal C}(\bm{\tau})=\left\{\bm{y}\in \mathcal{H}_{{\bf C}^{-1}}: \text{ for any }\bm{c}\in \mathcal{C} \text{ we have }  \la\bm{y}, \bm{c}-\bm{\tau}\ra_{{\bf C}^{-1}}\leqslant 0 \right\}.
\label{eq:nc-C-inv-def}
\end{equation}
 \noindent On top of Proposition \ref{prop:nc-abstract-properties}, such normal cone has the following properties.

\begin{Proposition}
\label{prop:c-inv-normal-cone-properties}
Formula \eqref{eq:nc-C-inv-def} yields the following properties:
\begin{enumerate}[{\it i)}] 
\item For any $\bm{\tau}\in \mathcal{C}$
\begin{equation}
N^{{\bf C}^{-1}}_{\mathcal C}(\bm{\tau}) = {\bf C} N_{\mathcal C}(\bm{\tau}),
\label{eq:regular-and-weighted-nc}
\end{equation}
\item For any $\bm{v}\in \mathcal{V}$ we have
\begin{equation}
\mathcal{U}= N^{{\bf C}^{-1}}_{\mathcal{V}}(\bm{v})
\label{eq:orthogonal-space-nc}
\end{equation}
\end{enumerate}
\end{Proposition}

Now we are ready to derive the differential inclusions to solve the problem of Definition \ref{def:aepp}.

\begin{Theorem}
\label{th:aepp-to-sweeping-process}
Functions \eqref{eq:aepp-unknowns} solve the abstract problem of quasi-static evolution in elasticity-perfect plasticity if and only if the unknown
\begin{equation}
y = \sigma-\widetilde{\sigma}\in W^{1,\infty}(I, \mathcal{H}_{{\bf C}^{-1}})
\label{eq:yielding-to-sweeping}
\end{equation}
and the unknown $\varepsilon \in W^{1,\infty}(I, \mathcal{H})$ solve the differential inclusions
\begin{numcases}{}
-\frac{d}{dt}\bm{y}\in N^{{\bf C}^{-1}}_{\mathcal{C}(t)}(\bm{y}), \label{eq:sp-elastoplastic}\\
\frac{d}{dt} \bm{\varepsilon} \in {\bf C}^{-1}\left(\left(N^{{\bf C}^{-1}}_{\Sigma-\bm{\widetilde{\sigma}}(t)}(\bm{y})+ \frac{d}{dt} \bm{y}\right)\cap \mathcal{U}+ \frac{d}{dt} \bm{\widetilde{\sigma}}(t)\right) \label{eq:inclusion2-elastoplastic}
\end{numcases}
with the initial conditions
\begin{align}
\bm{y}(0) &= \bm{\sigma}_{\bf 0}- \bm{\widetilde{\sigma}}(0)\in \mathcal{C}(0), \label{eq:diff-incs-ic-y}\\
\bm{\varepsilon}(0)& = \bm{\varepsilon}_{\bf 0}\in {\bf C}^{-1}(\mathcal{U}+\bm{\widetilde{\sigma}}(0)),
\label{eq:diff-incs-ic-eps}
\end{align}
where the moving set is
\begin{equation}
\mathcal{C}(t)=\left(\Sigma-\bm{\widetilde{\sigma}}(t)\right)\cap \mathcal{V}.
\label{eq:aepp-sp-moving-set}
\end{equation}
\end{Theorem}
\begin{Remark} Notice, that differential inclusion \eqref{eq:sp-elastoplastic} is a sweeping process, which does not depend on $\varepsilon$, hence a unique $y$ can be found whenever the sweeping process is solvable.  Moreover, the right-hand side of \eqref{eq:inclusion2-elastoplastic} is independent of $\varepsilon$ as well, and, if we already know $y$ from  \eqref{eq:sp-elastoplastic}, \eqref{eq:diff-incs-ic-y}, we could find $\varepsilon$ by integrating a {\it measurable selection} from the right-hand side of \eqref{eq:inclusion2-elastoplastic}. We will discuss the existence and integrability of such selection for the solution of the sweeping process  \eqref{eq:sp-elastoplastic} after the proof of Theorem \ref{th:aepp-to-sweeping-process}.
\end{Remark}
\noindent{\bf Proof of Theorem \ref{th:aepp-to-sweeping-process}.} Since \eqref{eq:aepp-5} and \eqref{eq:ae-3} are the same equation on the unknown, by Proposition \ref{prop:elasticity-affine-planes} {\it ii)} we have equation \eqref{eq:aepp-5} equivalent to
\begin{equation*}
\bm{\sigma}\in \mathcal{V}+ Q\,\bm{f}.
\end{equation*}
Notice, that $G$ takes values in $\mathcal{V}$, therefore we can write
\begin{equation*}
\bm{\sigma}\in \mathcal{V}+ Q\,\bm{f} + G\, \bm{g},
\end{equation*}
i.e. (recall \eqref{eq:abstract_linear_solution})
\begin{equation}
\bm{\sigma}- \bm{\widetilde{\sigma}}\in \mathcal{V}.
\label{eq:stress-constraint}
\end{equation}
In particular, together with \eqref{eq:aepp-sigma-ic-compatible} this implies the compatibility of the initial value \eqref{eq:diff-incs-ic-y}.

In turn, since \eqref{eq:aepp-1} is also the same as \eqref{eq:ae-1}, by Proposition \ref{prop:elasticity-affine-planes} {\it i)} we have equation \eqref{eq:aepp-5} equivalent to
\begin{equation*}
{\bf C}\, \bm{\varepsilon}\in \mathcal{U}+ G \,\bm{g},
\end{equation*}
and, since $Q$ takes values in $\mathcal{U}$,
\[
{\bf C}\, \bm{\varepsilon}\in \mathcal{U}+ G \,\bm{g} + Q\, \bm{f},
\]
i.e. (recall \eqref{eq:abstract_linear_solution} again)
\begin{equation}
{\bf C}\, \bm{\varepsilon}\in \mathcal{U}+ \bm{\widetilde{\sigma}}.
\label{eq:compatibility-affine-constraint-proof}
\end{equation}
This yields the compatibility of the initial value \eqref{eq:diff-incs-ic-eps} and the equation on the strain rate 
\begin{equation}
\frac{d}{dt} {\bf C}\, \bm{\varepsilon}\in \mathcal{U}+ \frac{d}{dt} \bm{\widetilde{\sigma}}.
\label{eq:rate-strain-constraint}
\end{equation}
On the other hand, from \eqref{eq:aepp-2}--\eqref{eq:aepp-3} we get that
\[
{\bf C}\, \bm{\varepsilon} = \bm{\sigma}+ {\bf C} \, \bm{\varepsilon}_{\bf p},
\]
which also yields the following rate equation for a.a. $t\in I$ 
\begin{equation}
-\frac{d}{dt} \bm{\sigma}  = {\bf C} \,\frac{d}{dt} \bm{\varepsilon}_{\bf p} -\frac{d}{dt} {\bf C}\, \bm{\varepsilon},
\label{eq:rate-additive-decomposition}
\end{equation}
where ${\bf C}$ commutes with $\frac{d}{dt}$ due to Proposition \ref{prop:prelim-classical-derivatives-ae} and the continuity of ${\bf C}$. 

Observe that by applying ${\bf C}$ to \eqref{eq:aepp-4} and using the expression \eqref{eq:regular-and-weighted-nc} at the right-hand side, we get
\begin{equation}
{\bf C} \frac{d}{dt} \bm{\varepsilon}_{\bf p} \in N^{{\bf C}^{-1}}_{\Sigma}(\bm{\sigma}).
\label{eq:plastic-flow-rule-in-proof}
\end{equation}
Substitute \eqref{eq:rate-strain-constraint} and \eqref{eq:plastic-flow-rule-in-proof} into \eqref{eq:rate-additive-decomposition}:
\[
-\frac{d}{dt} \bm{\sigma}\in N^{{\bf C}^{-1}}_{\Sigma}(\bm{\sigma}) - \mathcal{U}-  \frac{d}{dt} \bm{\widetilde{\sigma}}.
\]
Recall that $\mathcal{U}=-\mathcal{U}$. Therefore,
\[
-\left(\frac{d}{dt} \bm{\sigma}-  \frac{d}{dt} \bm{\widetilde{\sigma}} \right)\in N^{{\bf C}^{-1}}_{\Sigma}(\bm{\sigma}) + \mathcal{U}.
\]
We use \eqref{eq:argument-sum-nc} in the normal cone term and substitute $\mathcal{U}$ as in \eqref{eq:orthogonal-space-nc} using \eqref{eq:stress-constraint}:
\[
-\left(\frac{d}{dt} \bm{\sigma}-  \frac{d}{dt} \bm{\widetilde{\sigma}} \right)\in N^{{\bf C}^{-1}}_{\Sigma-\bm{\widetilde{\sigma}}}(\bm{\sigma}-\bm{\widetilde{\sigma}}) + N^{{\bf C}^{-1}}_{\mathcal{V}}(\bm{\sigma}-\bm{\widetilde{\sigma}}).
\]
Next, we use \eqref{eq:nc-subadditivity} to get
\[
-\frac{d}{dt} \left( \bm{\sigma}- \bm{\widetilde{\sigma}} \right)\in N^{{\bf C}^{-1}}_{\left(\Sigma-\bm{\widetilde{\sigma}}\right)\cap \mathcal{V}}(\bm{\sigma}-\bm{\widetilde{\sigma}})
\]
and substitute \eqref{eq:yielding-to-sweeping} to obtain the sweeping process \eqref{eq:sp-elastoplastic}.

On the other hand, we can express the term with $\bm{\varepsilon}$ from \eqref{eq:rate-additive-decomposition}:
\[
\frac{d}{dt} {\bf C}\, \bm{\varepsilon} = \frac{d}{dt} \bm{\sigma}+ {\bf C} \,\frac{d}{dt} \bm{\varepsilon}_{\bf p},
\]
substitute there \eqref{eq:plastic-flow-rule-in-proof}:
\begin{equation}
\frac{d}{dt} {\bf C}\, \bm{\varepsilon} \in \frac{d}{dt} \bm{\sigma}+ N^{{\bf C}^{-1}}_{\Sigma}(\bm{\sigma}),
\label{eq:additive-dec-rate-nc-proof}
\end{equation}
and insert there $\bm{\widetilde{\sigma}}$:
\[
\frac{d}{dt} {\bf C}\, \bm{\varepsilon} \in N^{{\bf C}^{-1}}_{\Sigma-\bm{\widetilde{\sigma}}}(\bm{\sigma}-\bm{\widetilde{\sigma}})+ \frac{d}{dt} (\bm{\sigma}-\bm{\widetilde{\sigma}})+ \frac{d}{dt} \bm{\widetilde{\sigma}}.
\]
Combine this with \eqref{eq:rate-strain-constraint} to obtain
\begin{equation}
\frac{d}{dt} {\bf C}\, \bm{\varepsilon} \in \left(N^{{\bf C}^{-1}}_{\Sigma-\bm{\widetilde{\sigma}}}(\bm{\sigma}-\bm{\widetilde{\sigma}})+ \frac{d}{dt} (\bm{\sigma}-\bm{\widetilde{\sigma}})\right)\cap \mathcal{U}+ \frac{d}{dt} \bm{\widetilde{\sigma}},
\label{eq:near-ready-sp-proof}
\end{equation}
and then substitute \eqref{eq:yielding-to-sweeping} and apply ${\bf C}^{-1}$ to obtain the inclusion \eqref{eq:inclusion2-elastoplastic}.

Assume now that we have a solution 
\[
(y, \varepsilon)\in W^{1,\infty}(I, \mathcal{H}_{{\bf C}^{-1}})\times W^{1,\infty}(I, \mathcal{H})
\]
of \eqref{eq:sp-elastoplastic}--\eqref{eq:diff-incs-ic-eps}. Define $\sigma$ via $y$ according to \eqref{eq:yielding-to-sweeping}. Apply ${\bf C}$ to \eqref{eq:inclusion2-elastoplastic} and substitute \eqref{eq:yielding-to-sweeping} to get \eqref{eq:near-ready-sp-proof}, which is equivalent to \eqref{eq:additive-dec-rate-nc-proof} and \eqref{eq:rate-strain-constraint} together. Define $\varepsilon_{\rm el}$ and $\varepsilon_{\rm p}$ via both $\varepsilon$ and $\sigma$ according to \eqref{eq:aepp-2}--\eqref{eq:aepp-3}. We have that \eqref{eq:additive-dec-rate-nc-proof}, \eqref{eq:aepp-2} and \eqref{eq:regular-and-weighted-nc} imply \eqref{eq:aepp-4}. In turn, \eqref{eq:rate-strain-constraint} and the initial condition \eqref{eq:diff-incs-ic-eps} imply \eqref{eq:compatibility-affine-constraint-proof}, which is equivalent to \eqref{eq:aepp-1} by Proposition \ref{prop:elasticity-affine-planes} {\it i)}. Finally, from \eqref{eq:sp-elastoplastic} we have \eqref{eq:stress-constraint}, which means \eqref{eq:aepp-5} by Proposition \ref{prop:elasticity-affine-planes} {\it ii)}. $\blacksquare$

Now let us discuss the existence of solutions to \eqref{eq:sp-elastoplastic} and \eqref{eq:inclusion2-elastoplastic}.
\begin{Theorem}
\label{th:aepp-well-posedness-safe-load-strict} 
In the setting of Theorem \ref{th:aepp-to-sweeping-process} assume that Slater's constraint qualification \eqref{eq:nontrivial-intersection-condition} holds, i.e.
\begin{equation}
\left({\rm int}\, \Sigma-\bm{\widetilde{\sigma}}(t)\right)\cap \mathcal{V} \neq \varnothing \qquad \text{for all }t\in I.
\label{eq:safe-load-strict}
\end{equation}
Then 
\begin{enumerate}[{\it i)}]
\item \label{enum:th:aepp-well-posedness-sp} the sweeping process \eqref{eq:sp-elastoplastic}, \eqref{eq:diff-incs-ic-y}, \eqref{eq:aepp-sp-moving-set} satisfies the conditions of Theorem \ref{th:abstract-sp-existence} and it has a unique solution $y\in  W^{1,\infty}(I, \mathcal{H}_{{\bf C}^{-1}})$,
\item \label{enum:th:aepp-well-posedness-nonempty-rhs} the right-hand side of differential inclusion \eqref{eq:inclusion2-elastoplastic} is a nonempty set for a.a. $t\in I$.
\item  \label{enum:th:aepp-well-posedness-estimate-rhs} there exists at least one $\varepsilon \in W^{1,\infty}(I, \mathcal{H})$ which solves \eqref{eq:inclusion2-elastoplastic}, \eqref{eq:diff-incs-ic-eps}.
\end{enumerate}
\end{Theorem}
\noindent {\bf Proof.} 
\begin{enumerate}[{\it i)}]
\item For each $t\in I$ we have the moving set \eqref{eq:aepp-sp-moving-set} being closed, convex and nonempty. Moreover, the set $\Sigma-\bm{\widetilde{\sigma}}(t)$ is a time-dependent translation, therefore it is Lipschitz-continuous with the constant $\left\|\widetilde{\sigma}\right\|_{W^{1,\infty}(I, \mathcal{H}_{{\bf C}^{-1}})}$ (see \eqref{eq:ae-time-regularity}) with respect to Hausdorff distance. In turn, condition \eqref{eq:safe-load-strict} and the boundness of $\Sigma$ are sufficient for the general argument of \cite[Section 5.c, pp. 274--278]{Moreau1973}, \cite{Moreau1973-2}, \cite{Moreau1975}. From there we deduce that $\mathcal{C}(t)$, being the intersection \eqref{eq:aepp-sp-moving-set}, is also Lipschitz-continuous with respect to Hausdorff distance. 

Note: a simpler and more explicit way to show the Lipschitz continuity is possible under additional assumption that $\mathcal{C}(t)$ is a {\it polyhedron with fixed normal directions to its faces}. In such a case we can estimate the Lipschitz constant of $\mathcal{C}(t)$ by estimating the velocities of its vertices, of which there should be a finite number, see \cite[Lemma A.3, pp. 36--37]{Gudoshnikov2023ESAIM}. This will be the case for some of our further examples, see Remarks \ref{rem:ex11-plasticity-vertices-velocity}, \ref{rem:ex12-plasticity-vertices-velocity} below.

We have shown that all of the conditions of Theorem \ref{th:abstract-sp-existence} hold for the sweeping process \eqref{eq:sp-elastoplastic}, \eqref{eq:diff-incs-ic-y}, \eqref{eq:aepp-sp-moving-set}, which means that the statement {\it i)} of the current proposition also holds.
\item 
 Condition \eqref{eq:safe-load-strict} now plays the role of \eqref{eq:nontrivial-intersection-condition} and it guarantees that for all $t\in I$
\begin{equation}
N^{{\bf C}^{-1}}_{\left(\Sigma-\bm{\widetilde{\sigma}}(t)\right)\cap \mathcal{V}}(\bm{y}(t)) = N^{{\bf C}^{-1}}_{\Sigma-\bm{\widetilde{\sigma}}(t)}(\bm{y}(t))+ N^{{\bf C}^{-1}}_{\mathcal{V}}(\bm{y}(t)).
\label{eq:nc-additivity-in-aepp-proof}
\end{equation}
Hence for a.a. $t\in I$ we have 
\[
-\frac{d}{dt}\bm{y}(t)\in N^{{\bf C}^{-1}}_{\Sigma-\bm{\widetilde{\sigma}}(t)}(\bm{y}(t)) + \mathcal{U}.
\]
Thus for a.a. $t\in I$ there exists an element $\bm{\omega}(t)\in -\mathcal{U}=\mathcal{U}$ such that
\[
\bm{\omega}(t) \in  N^{{\bf C}^{-1}}_{\Sigma-\bm{\widetilde{\sigma}}(t)}(\bm{y}(t)) + \frac{d}{dt}\bm{y}(t),
\]
i.e. 
\begin{equation}
\bm{\omega}(t) \in  \left(N^{{\bf C}^{-1}}_{\Sigma-\bm{\widetilde{\sigma}}(t)}(\bm{y}(t)) + \frac{d}{dt}\bm{y}(t)\right)\cap \mathcal{U},
\label{eq:aepp-omega-inclusion}
\end{equation}
from where the statement about \eqref{eq:inclusion2-elastoplastic} follows.
\item
It is enough to show that there exists a function $\omega\in L^{\infty}(I, \mathcal{H}_{{\bf C}^{-1}})$ such that \eqref{eq:aepp-omega-inclusion} holds. Indeed, then we can set
\begin{equation}
\bm{\varepsilon} (t) = \bm{\varepsilon_0}+{\bf C}^{-1}\left(\bm{\widetilde{\sigma}}(t)-\bm{\widetilde{\sigma}}(0)+\int\limits_0^t \bm{\omega}(s)\, ds\right),
\label{eq:aepp-eps-from-omega}
\end{equation}
which then is from $W^{1, \infty}(I, \mathcal{H})$, since so is $\widetilde{\sigma}$ by Remark \ref{rem:ae-time-regularity}.

To find such $\omega\in L^{\infty}(I, \mathcal{H}_{{\bf C}^{-1}})$ we follow the ideas  originally published in \cite[Section 6.d, pp. 315--318]{Moreau1973}, \cite[Section 2.e, pp. 12.14--12.18]{Moreau1973-3}, \cite[Chapter III]{Castaing1975}, and which also can be found in \cite[pp. 20--23]{Kunze2000}. 

\noindent\underline{Claim: The right-hand side of \eqref{eq:aepp-omega-inclusion} admits a measurable selection.} Nowadays, the topic of selections from measurable multifunctions is present in many monographs, of which we will use \cite[Chapter 8]{Aubin2009}. Denote by $J\subset I$ the set of full measure where the right-hand side of \eqref{eq:aepp-omega-inclusion} is well-defined and nonempty. By Proposition \ref{prop:nc-equivalent} \ref{enum:abstract-nc-in-proj-notation}, the normal cone there can be written as
\begin{equation}
N^{{\bf C}^{-1}}_{\Sigma-\bm{\widetilde{\sigma}}(t)}(\bm{y}(t)) = \left\{\bm{\tau}\in \mathcal{H}_{{\bf C}^{-1}}: {\rm proj}^{{\bf C}^{-1}}\left(\bm{\tau}+\bm{y}(t), \Sigma-\bm{\widetilde{\sigma}}(t)\right)=\bm{y}(t) \right\},
\label{eq:aepp-proof-nc-measurable}
\end{equation}
where by ${\rm proj}^{{\bf C}^{-1}}(\bm{\tau}, \mathcal{C})$ we denote the metric projection in the sense of the inner product \eqref{eq:weighted-ip} of a point $\bm{\tau}$ on a closed convex nonempty set $\mathcal{C}$. Such projection is nonexpansive w.r. to $\bm{\tau}$ in the corresponding norm, hence it is continuous. Overall, the mapping
\begin{align*}
\pi: J\times \mathcal{H}_{{\bf C}^{-1}}& \to \mathcal{H}_{{\bf C}^{-1}},\\ 
\pi: (t, \bm{\tau})& \mapsto {\rm proj}^{{\bf C}^{-1}}(\bm{\tau}+\bm{y}(t), \Sigma-\bm{\widetilde{\sigma}}(t))
\end{align*}
is Caratheodory \cite[p. 311]{Aubin2009}, and, by the Inverse Image Theorem  \cite[Th. 8.2.9, p. 315]{Aubin2009} used with $g=\pi$, $F:t\mapsto \mathcal{H}_{{\bf C}^{-1}}$ and $G:t \mapsto \{\bm{y}(t)\}$, the normal cone \eqref{eq:aepp-proof-nc-measurable} is a measurable multimapping \cite[Def. 8.1.1, p. 307]{Aubin2009} from the domain $J$ endowed with Lebesgue measure.

Moreover, the mapping $t\mapsto \frac{d}{dt}\bm{y}(t)$ is measurable by the choice of $y$. Therefore, by \cite[Th. 8.2.8, p. 314]{Aubin2009} the sum $N^{{\bf C}^{-1}}_{\Sigma-\bm{\widetilde{\sigma}}(t)}(\bm{y}(t)) + \frac{d}{dt}\bm{y}(t)$ is also a measurable multimapping from $J$. By \cite[Th. 8.2.4, p. 312]{Aubin2009} the intersection, which makes the right-hand side of \eqref{eq:aepp-omega-inclusion}, is also a measurable multimapping on $J$, thus, by \cite[Th. 8.1.3, p. 308]{Aubin2009}, there exists a measurable selection $\omega$. 

Finally, we must estimate $\|\omega\|_{L^{\infty}(I,\mathcal{H}_{{\bf C}^{-1}})}<\infty$, for which we need to prove the following.

\noindent \underline{Claim: under condition \eqref{eq:safe-load-strict} the following statement is true:}
\begin{multline}
\text{there exists}\quad \rho>0 \quad \text{such that for any}\quad t\in I \quad\text{there exists} \quad \bm{h}_t\in\mathcal{C}(t) \\ \text{such that} \quad B_\rho (\bm{h}_t)\subset \Sigma -\bm{\widetilde{\sigma}}(t),
\label{eq:aepp-safe-load-bound-claim}
\end{multline}
where $B_\rho (\bm{h}_t)$ is an open ball in $\mathcal{H}_{{\bf C}^{-1}}$ with a center $\bm{h}_t$ and a radius $\rho$, and $\mathcal{C}(t)$ is as in \eqref{eq:aepp-sp-moving-set}.
Assume the contrary, namely that
\begin{multline*}
\text{for any}\quad \rho>0 \quad \text{there exists}\quad t\in I \quad \text{such that for all} \quad \bm{h}\in\mathcal{C}(t) \\ \text{ there exists }\quad \bm{\tau}_t\in B_\rho (\bm{h}) \quad \text{such that}\quad \bm{\tau}_t \notin \Sigma -\bm{\widetilde{\sigma}}(t).
\end{multline*}
Take a sequence $\rho_i\to 0, i\in \mathbb{N}$ and redenote it with a subsequence such that $t_i\to t^*$ for some $t^*\in I$. We have that
\begin{equation}
\text{for all} \quad i\in \mathbb{N},\quad \bm{h}\in\mathcal{C}(t_i)\quad \text{ there exists }\quad \bm{\tau}_{i}\in B_{\rho_i} (\bm{h}) \quad \text{such that}\quad \bm{\tau}_{i} \notin \Sigma -\bm{\widetilde{\sigma}}(t_i).
\label{eq:aepp-safe-load-bound-claim-contrary}
\end{equation}

On the other hand, condition \eqref{eq:safe-load-strict} means that
\begin{equation}
\text{there exist}\qquad \rho^*>0\quad \text{and}\quad \bm{h}^*\in\mathcal{C}(t^*) \quad\text{such that} \quad B_{\rho^*}(\bm{h}^*) \subset \Sigma-\bm{\widetilde{\sigma}}(t^*).
\label{eq:eq:safe-load-strict-aepp-proof-implication}
\end{equation}
However, as we discussed in the part {\it i)} of the proof, condition \eqref{eq:safe-load-strict} implies Lipschitz-continuity of $\mathcal{C}(t)$ with respect to the Hausdorff distance. Thus, 
\begin{multline}
\text{for any small number}\quad \delta>0 \quad \text{we can find} \quad N\in\mathbb{N}\\ \text{ such that for any}\quad i\geqslant N\quad \text{we will have}\quad \bm{h}^*\in B_\delta(0)+\mathcal{C}(t_i).
\label{eq:aepp-proof-haus-implication}
\end{multline}
We proceed further in three steps. 
\begin{enumerate}[1.]
\item Take $\delta<\frac{\rho^*}{3}$ and use \eqref{eq:aepp-proof-haus-implication} to find the corresponding $N\in \mathbb{N}$. Denote 
\[
\bm{h}_{i}={\rm proj}^{{\bf C}^{-1}}(\bm{h}^*, \mathcal{C}(t_i)) \in B_\delta(\bm{h}^*)
\]
 for $i\geqslant N$. 
\item Use the Lipschitz-continuity of $\widetilde{\sigma}$ (Remark \ref{rem:ae-time-regularity}) to make sure that $N$ is large enough for 
\[
\|\bm{\widetilde{\sigma}}(t_i)-\bm{\widetilde{\sigma}}(t^*)\|_{{\bf C}^{-1}}<\frac{\rho^*}{3}.
\]
Redefine $N$, if necessary.
\item 
 Now choose $i\geqslant N$ such that $\rho_i<\frac{\rho^*}{3}$, and use \eqref{eq:aepp-safe-load-bound-claim-contrary} to find the corresponding $\bm{\tau}_{i}\in B_{\rho_i}(\bm{h}_i)$.
\end{enumerate}
 Observe that
\begin{multline*}
\left\|\bm{\tau}_i+\bm{\widetilde{\sigma}}(t_i)-\bm{\widetilde{\sigma}}(t^*)-\bm{h}^*\right\|_{{\bf C}^{-1}}\leqslant \|\bm{\tau}_i-\bm{h}_i\|_{{\bf C}^{-1}}+\|\bm{\widetilde{\sigma}}(t_i)-\bm{\widetilde{\sigma}}(t^*)\|_{{\bf C}^{-1}}+\|\bm{h}_i-\bm{h}^*\|_{{\bf C}^{-1}}<\\<\rho_i +\frac{\rho^*}{3}+\delta<\rho^*,
\end{multline*}
i.e.
\begin{equation}
\bm{\tau}_i+\bm{\widetilde{\sigma}}(t_i)-\bm{\widetilde{\sigma}}(t^*) \in B_{\rho^*}(\bm{h}^*).
\label{eq:aepp-proof-in-ball-to-contadict}
\end{equation}
But, by the construction in \eqref{eq:aepp-safe-load-bound-claim-contrary},
\[
\bm{\tau}_i\notin \Sigma -\bm{\widetilde{\sigma}}(t_i),
\]
i.e.
\[
\bm{\tau}_i+\bm{\widetilde{\sigma}}(t_i)-\bm{\widetilde{\sigma}}(t^*) \notin \Sigma -\bm{\widetilde{\sigma}}(t^*),
\]
which, together with \eqref{eq:aepp-proof-in-ball-to-contadict}, contradicts \eqref{eq:eq:safe-load-strict-aepp-proof-implication}. Therefore, condition \eqref{eq:safe-load-strict} yields our claim \eqref{eq:aepp-safe-load-bound-claim}.

\noindent \underline{Claim: $\|\omega\|_{L^{\infty}(I,\mathcal{H}_{{\bf C}^{-1}})}<\infty$.} From \eqref{eq:aepp-safe-load-bound-claim} we get that for every $t\in I$ and every $\bm{\tau}\in \mathcal{H}_{{\bf C}^{-1}}$
\begin{equation}
\delta^*_{B_\rho(\bm{h}_t)}(\bm{\tau})\leqslant \delta^*_{\Sigma -\bm{\widetilde{\sigma}}(t)}(\bm{\tau}),
\label{eq:support-functions-embedded-proof}
\end{equation}
where $\delta^*$ denotes the support function \eqref{eq:def-support-function} and we use the property \eqref{eq:sp-subset-property}.
Denote 
\[
\bm{\zeta}(t)=\bm{\omega}(t)- \frac{d}{dt}\bm{y}(t)
\] for $t\in J$ and plug it into \eqref{eq:support-functions-embedded-proof} as $\bm{\tau}$ to  obtain
\[
\rho \left\|\bm{\zeta}(t) \right\|_{{\bf C}^{-1}} + \left\la\bm{\zeta}(t) , \bm{h}_t\right\ra_{{\bf C}^{-1}} \leqslant \delta^*_{\Sigma -\bm{\widetilde{\sigma}}(t)}\left(\bm{\zeta}(t) \right).
\]
On the other hand, from \eqref{eq:aepp-omega-inclusion} we have that
\[
\bm{\zeta}(t) \in N^{{\bf C}^{-1}}_{\Sigma-\bm{\widetilde{\sigma}}(t)}(\bm{y}(t)),
\]
i.e. by Proposition \ref{prop:nc-equivalent} \ref{enum:abstract-nc-in-support-func-notation}
\[
\left\la\bm{\zeta}(t), \bm{y}(t)\right\ra_{{\bf C}^{-1}} = \delta^*_{\Sigma -\bm{\widetilde{\sigma}}(t)}\left(\bm{\zeta}(t), \bm{y}(t)\right).
\]
Therefore,
\[
\rho \left\|\bm{\zeta}(t) \right\|_{{\bf C}^{-1}} + \left\la\bm{\zeta}(t) , \bm{h}_t\right\ra_{{\bf C}^{-1}} \leqslant \left\la\bm{\zeta}(t), \bm{y}(t)\right\ra_{{\bf C}^{-1}},
\]
i.e.
\[
\rho \left\|\bm{\omega}(t)- \frac{d}{dt}\bm{y}(t) \right\|_{{\bf C}^{-1}} \leqslant \left\la\bm{\omega}(t)- \frac{d}{dt}\bm{y}(t), \bm{y}(t)-\bm{h}_t\right\ra_{{\bf C}^{-1}}.
\]
Recall, however, that $\bm{\omega}(t)\in \mathcal{U}$ by \eqref{eq:aepp-omega-inclusion}, and both $\bm{y}(t)$ and $\bm{h}_t$ belong to $\mathcal{C}(t)\subset \mathcal{V}$ by \eqref{eq:sp-elastoplastic} and \eqref{eq:aepp-safe-load-bound-claim}, respectively. Since $\mathcal{U}$ and $\mathcal{V}$ are linear ${\bf C}^{-1}$-orthogonal subspaces, the corresponding term vanishes and we get that
\[
\rho \left\|\bm{\omega}(t)- \frac{d}{dt}\bm{y}(t) \right\|_{{\bf C}^{-1}} \leqslant \left\la- \frac{d}{dt}\bm{y}(t), \bm{y}(t)-\bm{h}_t\right\ra_{{\bf C}^{-1}}\leqslant \left\|\frac{d}{dt}\bm{y}(t)\right\|_{{\bf C}^{-1}}\left\|\bm{y}(t)-\bm{h}_t\right\|_{{\bf C}^{-1}}.
\]
Furthermore, the fact that $\bm{\omega}(t)\in \mathcal{U}$ and $\frac{d}{dt}\bm{y}(t)\in \mathcal{V}$ yields the second equality in the following estimate:
\begin{multline*}
\rho \left\|\bm{\omega}(t)\right\|_{{\bf C}^{-1}} = \rho \left(\left\la \bm{\omega}(t),  \bm{\omega}(t) \right\ra_{{\bf C}^{-1}}\right)^{\frac{1}{2}} \leqslant \rho \left(\left\la \bm{\omega}(t),  \bm{\omega}(t) \right\ra_{{\bf C}^{-1}} + \left\la \frac{d}{dt}\bm{y}(t),  \frac{d}{dt}\bm{y}(t) \right\ra_{{\bf C}^{-1}} \right)^{\frac{1}{2}}=
\\[2mm]
= \rho \left(\left\la \bm{\omega}(t),  \bm{\omega}(t) \right\ra_{{\bf C}^{-1}} -2\left\la \bm{\omega}(t),  \frac{d}{dt}\bm{y}(t) \right\ra_{{\bf C}^{-1}}+ \left\la \frac{d}{dt}\bm{y}(t),  \frac{d}{dt}\bm{y}(t) \right\ra_{{\bf C}^{-1}} \right)^{\frac{1}{2}}= \\[2mm]
=\rho \left\|\bm{\omega}(t)- \frac{d}{dt}\bm{y}(t) \right\|_{{\bf C}^{-1}} \leqslant \left\|\frac{d}{dt}\bm{y}(t)\right\|_{{\bf C}^{-1}}\left\|\bm{y}(t)-\bm{h}_t\right\|_{{\bf C}^{-1}}.
\end{multline*}
Since set $\Sigma$ is bounded (say, by a constant $L_\Sigma>0$) and $\bm{h}_t\in \Sigma - \bm{\widetilde{\sigma}}(t)$ by \eqref{eq:aepp-safe-load-bound-claim}, we have the estimate
\[
 \left\|\bm{\omega}(t)\right\|_{{\bf C}^{-1}} \leqslant \frac{1}{\rho} \|y\|_{W^{1,\infty}(I, \mathcal{H}_{\bf{C}^{-1}})} \left(\|y\|_{W^{1,\infty}(I, \mathcal{H}_{\bf{C}^{-1}})}+ L_\Sigma +  \|\widetilde{\sigma}\|_{W^{1,\infty}(I, \mathcal{H}_{\bf{C}^{-1}})}\right)
\]
for a.a. $t\in I$. Therefore, we indeed have $\omega\in L^{\infty}(I,\mathcal{H}_{{\bf C}^{-1}})$ and we can use \eqref{eq:aepp-eps-from-omega} to construct $\varepsilon \in W^{1,\infty}(I, \mathcal{H})$.
\end{enumerate}
$\blacksquare$
\begin{Corollary}
\label{cor:aepp-full-solution}
 Given the data as in Definition \ref{def:aepp} such that the condition \eqref{eq:safe-load-strict} holds (with $\widetilde{\sigma}$ found via \eqref{eq:abstract_linear_solution}), the problem \eqref{eq:sp-elastoplastic}--\eqref{eq:aepp-sp-moving-set} has at least one solution $(y, \varepsilon)\in W^{1,\infty}(I, \mathcal{H}_{{\bf C}^{-1}})\times W^{1,\infty}(I, \mathcal{H})$, and the abstract problem of quasi-static evolution in elasticity-perfect plasticity (Definition \ref{def:aepp}) also has at least one solution. The unknowns $y, \sigma$ and $\varepsilon_{\rm el}$ can be determined uniquely, while the unknowns $\varepsilon$ and $\varepsilon_{\rm p}$ may be solvable non-uniquely.
\end{Corollary}
\begin{Remark}
\label{rem:aepp-sp-in-V}
Notice, that the moving set \eqref{eq:aepp-sp-moving-set} always lays in the subspace $\mathcal{V}$, i.e. one can think of the sweeping process \eqref{eq:sp-elastoplastic} as being defined in Hilbert space $\mathcal{V}$ with the inner product induced from $\mathcal{H}_{{\bf C}^{-1}}$. Moreover, recall that in \eqref{eq:abstract_linear_solution} the first and the second terms take values in orthogonal subspaces $\mathcal{V}$ and $\mathcal{U}$, respectively. Therefore, the change of the load $\bm{f}(t)$ causes the change of the {\it shape} of the moving set \eqref{eq:aepp-sp-moving-set}, while the change of the load $\bm{g}(t)$ only causes a parallel translation of $\mathcal{C}(t)$ within $\mathcal{V}$. In particular, this leads to the following corollaries.
\end{Remark}
\begin{Corollary}
\label{cor:constant-force}
Assume that during a time-interval $I$ the load $\bm{f}(t)$ remains constant and such that the moving set \eqref{eq:aepp-sp-moving-set} is nonempty (for at least one $t$, which implies nonemptiness for all $t\in I$). If the load $\bm{g}(t)$ has regularity as in \eqref{eq:loads-abstract-def}, then the time-dependent set $\mathcal{C}(t)$ only moves by a parallel translation within $\mathcal{V}$. Hence $\mathcal{C}(t)$ is Lipschitz-continuous with respect to Hausdorff distance with Lipschitz constant $ \|G\|_{\rm op}\, \|g\|_{W^{1, \infty}(I, \mathcal{H})}$ and the conditions of Theorem \ref{th:abstract-sp-existence} are satisfied for the sweeping process \eqref{eq:sp-elastoplastic}, \eqref{eq:diff-incs-ic-y}, \eqref{eq:aepp-sp-moving-set} regardless of condition \eqref{eq:safe-load-strict}.
\end{Corollary}
\begin{Corollary}
Since in \eqref{eq:abstract_linear_solution} the operator $G$ takes values in $\mathcal{V}$, condition \eqref{eq:safe-load-strict} is equivalent to
\begin{equation}
\left({\rm int}\, \Sigma- Q\,{\bm f}(t)\right)\cap \mathcal{V} \neq \varnothing \qquad \text{for all }t\in I,
\label{eq:safe-load-forces-strict}
\end{equation}
while the condition in Theorem \ref{th:abstract-sp-existence} to have nonempty $\mathcal{C}(t)$ of the form \eqref{eq:aepp-sp-moving-set} is equivalent to 
\begin{equation}
\left(\Sigma- Q\,{\bm f}(t)\right)\cap \mathcal{V} \neq \varnothing \qquad \text{for all }t\in I. \label{eq:safe-load-forces-nonstrict}
\end{equation}
\end{Corollary}
\noindent Note that both \eqref{eq:safe-load-forces-strict} and \eqref{eq:safe-load-forces-nonstrict} are conditions on given external force $\bm{f}$, and they appear due to perfect plasticity. These conditions shall be considered along with the resolvability requirement \eqref{eq:resolvable-load-def}, which also applies to $\bm{f}$, but has a different nature and comes from quasi-static elasticity. In the context of perfect plasticity condition \eqref{eq:safe-load-forces-strict} and those similar to it are sometimes called {\it safe load conditions}, see e.g. Assumption~3 in \cite[p.~311]{Moreau1973}, \cite[p.~76]{Moreau1976}.

\subsection{Example 1 --- a discrete elastic-perfectly plastic model with one-dimensional sweeping process}
Now we modify the example of Section \ref{ssect:ex11-elasticity} to fit in the sweeping process framework. Recall the model of Fig. \ref{fig:discrete-models-elastic} {\it a)} and assume now that each of the two springs is elastic-perfectly plastic, i.e. that Hooke's law is valid for each spring only as long as the stress of the spring remains inside a given interval $(\sigma^-_i, \sigma^+_i), i\in \overline{1,2}$, which is called the {\it elastic range} of the spring. The threshold values $\sigma^-_i, \sigma^+_i$ are called {\it yield limits}, and we assume that for  $i\in \overline{1,2}$ 
\begin{equation}
\sigma_i^-<0<\sigma_i^+.
\label{eq:ex11-epp-yield-limits}
\end{equation}
Our assumption on the mechanical behavior of the springs is that 
\begin{equation*} 
\sigma_i(t) \in [\sigma^-_i, \sigma^+_i]
\end{equation*}
always holds. If a deformation is imposed upon a spring in such a manner, that Hooke's law would require the stress to go beyond the yield limits, then a plastic deformation occurs. Accumulated plastic deformation does not influence the evolution of stress, which is characteristic of {\it perfect plasticity}.  

The equations in Definition \ref{def:aepp} model such behavior. For the discrete model the unknowns $\varepsilon, \varepsilon_{\rm el}$ and $\varepsilon_{\rm p}$ are called, respectively, {\it total elongations}, {\it elastic elongations} (or the elastic component of the elongations) and {\it plastic elongations} (or the plastic component of the elongations). On top of the relations \eqref{eq:aepp-1}, \eqref{eq:aepp-3}, \eqref{eq:aepp-5}, which are similar to the elasticity problem and bear the same physical meaning, we have the {\it additive decomposition} \eqref{eq:aepp-2} of elongations into elastic and plastic components and the {\it normal form of the plastic flow rule} \eqref{eq:aepp-4}. The latter is the constitutive law for the plastic component, just like Hooke's law \eqref{eq:aepp-3} is the constitutive law for the elastic component, see Fig. \ref{fig:elasticity-pp-scheme}.

\begin{figure}[H]\center
\includegraphics{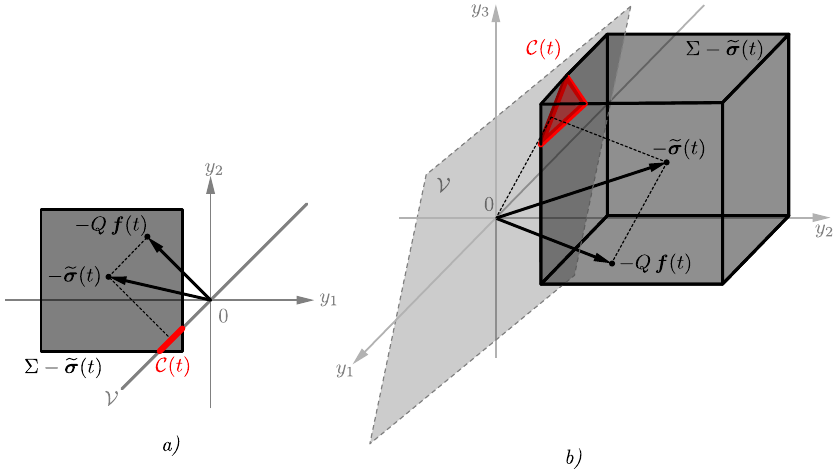}
\caption{\footnotesize The construction of the moving set in the discrete models of Examples 1.1 ({\it a}) and 1.2 ({\it b}). The figure shows the situation with the stiffness parameters $k_i$ equal to $1$.
} 
\label{fig:discrete-models-moving-sets}
\end{figure}

We take spaces $\mathcal{H}, \mathcal{X}, \mathcal{W}_0$, operators ${\bf E}, {\bf D}, {\bf C}$ and functions $\bm{g}, \bm{f}$ as in Section \ref{ssect:ex11-elasticity}, and the set
\begin{equation}
\Sigma = [\sigma^-_1, \sigma^+_1]\times[\sigma^-_2, \sigma^+_2]\subset \mathbb{R}^2\cong \mathcal{H}_{{\bf C}^{-1}}.
\label{eq:ex11-epp-sigma}
\end{equation}
Stress solution for elasticity is explicitly written as \eqref{eq:ex11-linear-solution}, and the moving set \eqref{eq:aepp-sp-moving-set} is 
\begin{multline}
\mathcal{C}(t) = \left([\sigma^-_1, \sigma^+_1]\times[\sigma^-_2, \sigma^+_2] -  \frac{l(t)}{k_1^{-1}+ k_2^{-1}} \begin{pmatrix}1 \\ 1\end{pmatrix} - \frac{F_2(t)}{k_1+k_2}\begin{pmatrix} k_1\\-k_2\end{pmatrix}\right)\cap {\rm Im}\, \begin{pmatrix}1\\1\end{pmatrix} = \\[2mm] \begin{aligned}
=\left\{\begin{pmatrix}y\\y\end{pmatrix}: y\in \mathbb{R} \text{ such that }\right. &\sigma^-_1\leqslant y+\frac{l(t)}{k_1^{-1}+ k_2^{-1}}+\frac{F_2(t)\, k_1}{k_1+k_2}\leqslant \sigma^+_1 \\  \text{ and } &\sigma^-_2\leqslant \left. y+\frac{l(t)}{k_1^{-1}+ k_2^{-1}}-\frac{F_2(t)\, k_2}{k_1+k_2}\leqslant \sigma^+_2 \right\},\end{aligned}
\label{eq:ex11-epp-moving-set}
\end{multline}
see Fig. \ref{fig:discrete-models-moving-sets} {\it a)}. In turn, safe load condition \eqref{eq:safe-load-forces-strict} (equivalent to Slater's condition \eqref{eq:safe-load-strict}) is written as
\begin{equation}
\begin{aligned}
\text{for each}\quad t\in I \quad\text{there exists}\quad y\in \mathbb{R} \quad \text{such that}\quad &\sigma^-_1<y+\frac{F_2(t)\, k_1}{k_1+k_2}<\sigma^+_1 \\ \text{and}\quad  &\sigma^-_2< y-\frac{F_2(t)\, k_2}{k_1+k_2}< \sigma^+_2.
\end{aligned}
\label{eq:ex11-epp-safe-load-strict}
\end{equation}
Notice from Fig.~\ref{fig:discrete-models-moving-sets} {\it a)} that the subspace $\mathcal{V}$ (written as \eqref{eq:ex11-space-V-exp}) will never intersect the axis-aligned rectangle $\Sigma$ along its edge, regardless of its parallel translation. Therefore, condition \eqref{eq:safe-load-strict} exactly means that the magnitude of the load $\bm{f}(t)$ is not too large, so that the segment $\mathcal{C}(t)$ does not degenerate to a point or an empty set. 

We also suggest an interested reader to illustrate the proof of Theorem \ref{th:aepp-well-posedness-safe-load-strict} \ref{enum:th:aepp-well-posedness-nonempty-rhs} and \ref{enum:th:aepp-well-posedness-estimate-rhs} by sketching the right-hand sides of \eqref{eq:inclusion2-elastoplastic} and \eqref{eq:aepp-omega-inclusion}, since in the context of the current example those are subsets of $\mathbb{R}^2$. Moreover, from such an illustration it is possible to observe non-uniqueness of the evolution of $\varepsilon$ when stresses in multiple springs reach the yield limits. 

\begin{Remark}
\label{rem:ex11-plasticity-vertices-velocity}
Since $\mathcal{C}(t)$ in the current example is a segment, i. e. it is a bounded {\it polyhedron with fixed normal directions to its faces} in a one-dimensional space  $\mathcal{V}$, we can estimate the Lipschitz constant of $\mathcal{C}(t)$ by estimating the velocities of its endpoints, see the note in the proof of Theorem \ref{th:aepp-well-posedness-safe-load-strict} \ref{enum:th:aepp-well-posedness-sp}.
\end{Remark}
We can now interpret Corollary \ref{cor:aepp-full-solution} as the following.
\begin{Corollary} Let us be given the data of Section \ref{ssect:ex11-elasticity} and set $\Sigma$ as in \eqref{eq:ex11-epp-sigma}, such that the condition \eqref{eq:ex11-epp-safe-load-strict} holds. Then the problem of quasi-static evolution in the model of Fig. \ref{fig:discrete-models-elastic} {\it a)} with elastic-perfectly plastic springs has at least one solution. In particular, the evolution of stresses $\bm{\sigma}(t)\in \mathbb{R}^2$ can be found as
\begin{equation}
\bm{\sigma}(t) = \bm{y}(t) + \bm{\widetilde{\sigma}}(t), 
\label{eq:finding-stresses-epp}
\end{equation}
where $\bm{\widetilde{\sigma}}$ is the stress solution for elasticity \eqref{eq:ex11-linear-solution} and ${\bm{y}}$ is the solution to the sweeping process \eqref{eq:sp-elastoplastic} in $\mathbb{R}^2$ with the inner product \eqref{eq:ex11-inner-product}. The moving set \eqref{eq:ex11-epp-moving-set} of the sweeping process is a time-dependent segment within the one-dimensional subspace $\mathcal{V}$, see  Fig. \ref{fig:discrete-models-moving-sets} {\it a)}.
\end{Corollary} 

\subsection{Example 2 --- a discrete elastic-perfectly plastic model with two-dimensional sweeping process}
Consider now the model of Fig. \ref{fig:discrete-models-elastic} {\it b)} and assume now that all three springs are elastic-perfectly plastic, i.e. that they are characterized by stiffness parameters $k_i$ and elastic ranges $(\sigma^-_i, \sigma^+_i), \, i\in \overline{1,3}$, for which \eqref{eq:ex11-epp-yield-limits} holds.

We take spaces $\mathcal{H}, \mathcal{X}, \mathcal{W}_0$, operators ${\bf E}, {\bf D}, {\bf C}$ and functions $\bm{g}, \bm{f}$ as in Section \ref{ssect:ex12-elasticity}, and the set
\begin{equation}
\Sigma = [\sigma^-_1, \sigma^+_1]\times[\sigma^-_2, \sigma^+_2]\times[\sigma^-_3, \sigma^+_3]\subset \mathbb{R}^3\cong \mathcal{H}_{{\bf C}^{-1}}.
\label{eq:ex12-epp-sigma}
\end{equation}
Stress solution for elasticity is explicitly written as \eqref{eq:ex12-linear-solution}, and the moving set \eqref{eq:aepp-sp-moving-set} is
\begin{multline}
\mathcal{C}(t) = \left([\sigma^-_1, \sigma^+_1]\times[\sigma^-_2, \sigma^+_2]\times[\sigma^-_3, \sigma^+_3] - \phantom{\begin{pmatrix} k_2^{-1}+k_3^{-1} & -k_2^{-1} \\[1mm]k_3^{-1} & k_1^{-1} \\[1mm]-k_2^{-1}& k_1^{-1}+ k_2^{-1} \end{pmatrix}}\right. \\ \left.  -\frac{1}{k_1^{-1}k_2^{-1}+k_2^{-1}k_3^{-1}+k_1^{-1}k_3^{-1}}\begin{pmatrix} k_2^{-1}+k_3^{-1} & -k_2^{-1} \\[1mm]k_3^{-1} & k_1^{-1} \\[1mm]-k_2^{-1}& k_1^{-1}+ k_2^{-1} \end{pmatrix}\begin{pmatrix}  l_1(t) \\  l_2(t)\end{pmatrix} - \frac{F_1(t)+F_3(t)}{k_1+k_2+k_3}\begin{pmatrix}-k_1 \\ k_2\\ -k_3\end{pmatrix}\right)\cap\\[2mm]
\cap {\rm Ker}\, \begin{pmatrix} -1 & 1 & -1\end{pmatrix} = \\[2mm]
\begin{aligned}=\left\{\begin{pmatrix}y_1\\y_2\\y_3\end{pmatrix}\in \mathbb{R}^3\right. &\text{ such that } -y_1+y_2-y_3=0\\
\text{ and } &\sigma^-_1\leqslant y_1+\frac{\left(k_2^{-1}+k_3^{-1}\right)l_1(t) -k_2^{-1} l_2(t)}{k_1^{-1}k_2^{-1}+k_2^{-1}k_3^{-1}+k_1^{-1}k_3^{-1}} -\frac{k_1\left(F_1(t)+F_3(t)\right)}{k_1+k_2+k_3}\leqslant \sigma^+_1\\[2mm]
 \text{ and } & \sigma^-_2\leqslant y_2+\frac{k_3^{-1}\, l_1(t)+ k_1^{-1}\, l_2(t)}{k_1^{-1}k_2^{-1}+k_2^{-1}k_3^{-1}+k_1^{-1}k_3^{-1}}+\frac{ k_2 \left(F_1(t)+F_3(t)\right)\,}{k_1+k_2+k_3}\leqslant \sigma^+_2\\[2mm]
\text{ and } &\sigma^-_3\leqslant y_3+\frac{-k_2^{-1} l_1(t)+\left(k_1^{-1}+k_2^{-1}\right) l_2(t)}{k_1^{-1}k_2^{-1}+k_2^{-1}k_3^{-1}+k_1^{-1}k_3^{-1}} -\frac{k_3\left(F_1(t)+F_3(t)\right)}{k_1+k_2+k_3}\leqslant \sigma^+_3
 \left.\vphantom{\begin{pmatrix}y_1\\y_2\\y_3\end{pmatrix}}\right\},\end{aligned}
\label{eq:ex12-epp-moving-set}
\end{multline}
see Fig. \ref{fig:discrete-models-moving-sets} {\it b)}. In turn, safe load condition \eqref{eq:safe-load-forces-strict} (equivalent to Slater's condition \eqref{eq:safe-load-strict}) is written as
\begin{equation}
\begin{aligned} \text{for each}\quad t\in I\quad \text{there exists}\quad \begin{pmatrix}y_1\\y_2\\y_3\end{pmatrix}\in \mathbb{R}^3\quad \text{such that}\quad & -y_1+y_2-y_3=0\\
\text{and}\quad &\sigma^-_1< y_1 -\frac{k_1\left(F_1(t)+F_3(t)\right)}{k_1+k_2+k_3}< \sigma^+_1\\[2mm]
 \text{and}\quad & \sigma^-_2< y_2+\frac{ k_2 \left(F_1(t)+F_3(t)\right)\,}{k_1+k_2+k_3}< \sigma^+_2\\[2mm]
\text{and}\quad &\sigma^-_3< y_3 -\frac{k_3\left(F_1(t)+F_3(t)\right)}{k_1+k_2+k_3}< \sigma^+_3.
\end{aligned}
\label{eq:ex12-epp-safe-load-strict}
\end{equation}
Again, notice from Fig.~\ref{fig:discrete-models-moving-sets} {\it b)} that the subspace $\mathcal{V}$ (written as \eqref{eq:ex12-space-V-exp}) will never intersect the axis-aligned rectangular cuboid $\Sigma$ along its edge or its facet, regardless of its parallel translation. 
Therefore, condition \eqref{eq:safe-load-strict} exactly means that the magnitude of the load $\bm{f}(t)$ is not too large, so that the polygon $\mathcal{C}(t)$ does not degenerate to a point or an empty set. 

\begin{Remark} 
\label{rem:ex12-plasticity-vertices-velocity}
Since $\mathcal{C}(t)$ in the current example is a polygon, i. e. it is again a bounded {\it polyhedron with fixed normal directions to its faces} in a two-dimensional space  $\mathcal{V}$, we again can estimate the Lipschitz constant of $\mathcal{C}(t)$ by estimating the velocities of its vertices, see \cite[Lemma A.3, pp. 36--37]{Gudoshnikov2023ESAIM}, cf. Remark \ref{rem:ex11-plasticity-vertices-velocity} and the note in the proof of Theorem \ref{th:aepp-well-posedness-safe-load-strict} \ref{enum:th:aepp-well-posedness-sp}.
\end{Remark}
We can now interpret Corollary \ref{cor:aepp-full-solution} as the following.
\begin{Corollary} Let us be given the data of Section \ref{ssect:ex12-elasticity} and set $\Sigma$ as in \eqref{eq:ex12-epp-sigma}, such that the condition \eqref{eq:ex12-epp-safe-load-strict} holds. Then the problem of quasi-static evolution in the model of Fig. \ref{fig:discrete-models-elastic} {\it b)} with elastic-perfectly plastic springs has at least one solution. In particular, the evolution of stresses $\bm{\sigma}(t)\in \mathbb{R}^3$ can be found as \eqref{eq:finding-stresses-epp}, where $\bm{\widetilde{\sigma}}$ is the stress solution for elasticity \eqref{eq:ex12-linear-solution} and ${\bm{y}}$ is the solution to the sweeping process \eqref{eq:sp-elastoplastic} in $\mathbb{R}^3$ with the inner product \eqref{eq:ex12-inner-product}. The moving set \eqref{eq:ex12-epp-moving-set} of the sweeping process is a time-dependent polygon (with fixed normal directions to its faces) within the two-dimensional subspace $\mathcal{V}$, see Fig. \ref{fig:discrete-models-moving-sets} {\it b)}.
\end{Corollary} 

The sweeping process framework of Definition \ref{def:aepp}, Theorem \ref{th:aepp-to-sweeping-process} and Slater's condition \eqref{eq:safe-load-strict} is well-applicable to the networks of elastic-perfectly plastic springs in general, of which the models in Fig. \ref{fig:discrete-models-elastic} are particular cases. For the general construction with one spatial dimension we refer to \cite{Gudoshnikov2021ESAIM}, \cite{Gudoshnikov2020}, and for many spatial dimensions we refer to \cite{Gudoshnikov2023preprint}. The former references focus on the models with the boundary conditions in terms of elongations, such as \eqref{eq:ex12-constraint-matrix-form} in Example 2, while in the latter reference we considered the boundary conditions in terms of displacements, i.e  of the type \eqref{eq:ex11-bc}. In \cite[Fig. 8, p. 3378]{Rachinskii2018}, \cite[Fig. 2, p. 16]{Gudoshnikov2023ESAIM} and \cite[Fig. 1, p.2]{Oleg2024preprint} the reader can find a very interesting network model, known in engineering as the {\it bridge circuit} \cite[Fig. 2.31, p. 58]{Sundararajan2020}, \cite[Fig. 11.15, p. 216]{Bracke2024reliability}. This model results in a sweeping process with three-dimensional $\mathcal{V}$ and $\mathcal{C}(t)$ (see \cite[Figs. 1, 3, pp. 16--17]{Gudoshnikov2023ESAIM}) and five-dimensional $\mathcal{H}$. Its quantities can also be written explicitly, but the formulas are even more cumbersome then those of the discrete examples in the present paper.

\section{Regularity lost (Example 3 with elasticity-perfect plasticity)}
\label{sect:regularity-lost}
\subsection{The sweeping process in the continuum model}
\label{ssect:continuum_rod_epp_general}
Let us now consider the continuous rod of Section \ref{ssect:ex21-elasticity}, endowed with the elastic-perfectly plastic behavior at each point $x\in \Omega$. Such rod is characterized by its stiffness function $\bm{C}\in L^\infty(\Omega)$ (see Section \ref{sssect:ex21-Hooke-law}) and its {\it local yield limits} $\sigma^-, \sigma^+\in L^{\infty}(\Omega)$, so that, analogously to \eqref{eq:ex11-epp-yield-limits} we require
\begin{equation}
\sigma^-(x)<0<\sigma^+(x) \qquad \text{for a.a. }x\in\Omega.
\label{eq:local_yield_limits}
\end{equation}
Functions $\sigma^-, \sigma^+$ are the boundaries of the {\it elastic range} at each point $x\in \Omega$, so that the one-sided constraint of the stress now has the form
\[
\sigma(t,x)\in \left[\sigma^-(x), \sigma^+(x)\right]\qquad \text{for all }t\in I \text{ and a.a. }x\in \Omega.
\] 

In the Definition \ref{def:aepp} we now use the data of Section \ref{ssect:ex21-elasticity} with 
\begin{equation}
\Sigma = \left\{\bm{\sigma}\in L^2(\Omega): \text{ for a.a. }x\in \Omega \text{ we have }\sigma^-(x)\leqslant\sigma(x)\leqslant\sigma^+(x) \right\}.
\label{eq:ex21-Sigma-def}
\end{equation}
The sweeping process \eqref{eq:sp-elastoplastic} happens within the subspace $\mathcal{V}$ of constant functions (see \eqref{eq:ex21-space-V}) in the space  $L^2_{{\bf C}^{-1}}(\Omega)$, which is $ L^2(\Omega)$ with the inner product \eqref{eq:ex21-ip}. We use the corresponding  stress solution for elasticity \eqref{eq:w0-const-def},\eqref{eq:ex21-linear-solution} to write the moving set \eqref{eq:aepp-sp-moving-set} as
\begin{multline}
\mathcal{C}(t) =\left\{\vphantom{ L^2_{{\bf C}^{-1}}}\right.\bm{y} \in L^2_{{\bf C}^{-1}}(\Omega): \text{ there exists }c\in \mathbb{R} \text{ such that for a.a. }x\in \Omega \text{ we have}\\ c\equiv y(x)\quad\text{and}\quad \sigma^-(x)-\widetilde{\sigma}(t,x) \leqslant c\leqslant \sigma^+(x)-\widetilde{\sigma}(t,x) \left.\vphantom{ L^2_{{\bf C}^{-1}}} \right\},
\label{eq:ex21-moving-set}
\end{multline}
which is closed and convex. For a reader's convenience we collect the new quantities, related to plasticity in all three examples, in Table \ref{tab:elasticity_perfect_plasticity_examples}.
\begin{table}[H]
\begin{center}
\renewcommand{\arraystretch}{1.2} 
\begin{tabular}{|c|c|c|c|}
\hline
\textbf{Quantity} & \textbf{Example 1 (Fig. \ref{fig:discrete-models-elastic} {\it a})} & \textbf{Example 2 (Fig. \ref{fig:discrete-models-elastic} {\it b})} & \textbf{Example 3 (Fig. \ref{fig:continuous-rod-model-elastic})}\\
\hline
$\Sigma\subset \mathcal{H}$ & $ [\sigma^-_1, \sigma^+_1]\times[\sigma^-_2, \sigma^+_2]$ & $ [\sigma^-_1, \sigma^+_1]\times[\sigma^-_2, \sigma^+_2]\times[\sigma^-_3, \sigma^+_3]$ &  $ \begin{array}{c}\left\{\bm{\sigma}\in L^2(\Omega): \text{ for a.a. }x\in \Omega \right.\\
\left.\sigma^-(x)\leqslant\sigma(x)\leqslant\sigma^+(x) \right\}\end{array}$
\\
\hline
$\begin{array}{c}
\text{Moving set}\\
\mathcal{C}(t)
\end{array}$
 & formula \eqref{eq:ex11-epp-moving-set} & formula \eqref{eq:ex12-epp-moving-set} & formula \eqref{eq:ex21-moving-set}
\\
\hline
\end{tabular}
\end{center}
\caption{\footnotesize The quantities in the examples of elasticity-perfect of the current paper, in addition to Table \ref{tab:linear_elasticity_examples}.}  
\label{tab:elasticity_perfect_plasticity_examples}
\end{table}

\subsection{The continuum model with particular data and its stress solution}
\label{ssect:ex21-concrete-data}
Let us examine the implications of the problem \eqref{eq:sp-elastoplastic}--\eqref{eq:aepp-sp-moving-set} with particular data describing the continuous rod. Consider the space- and time-domains
\[
\Omega=(-1,1), \qquad I=[0,T],\, T=3,
\]
the parameters of the rod
\[
 \bm{C}(x)\equiv 1,\qquad \sigma^-(x)\equiv -1, \qquad \sigma^+(x)\equiv 1,
\]
the initial stress state
\[
\sigma_0(x) \equiv 0 \qquad \text{for all }x\in \Omega
\]
and the following loading programme:
\[
F(t,x) = \begin{cases}
2xt &\text{when } t\in[0,1),\\
2x &\text{when } t\in[1, T],
\end{cases} \qquad u_a(t)\equiv 0 \text{ for all }t\in I, \qquad u_b=\begin{cases}0,& \text{when }t\in[0, 1),\\
t-1 &\text{when } t\in [1,T],
\end{cases}  
\]
i.e. we gradually increase the body force load from zero to $F(1,x) = 2x$ during time-interval $[0,1)$, and then keep the body force load constant, while monotonically increasing the required length in the displacement boundary condition.
We compute the stress solution for elasticity \eqref{eq:ex21-linear-solution} as
\begin{equation}
\widetilde{\sigma}(t,x) = \begin{cases} t\left(\frac{1}{3}-x^2\right) & \text{when } t\in [0,1),\\ \frac{1}{2}(t-1)+ \frac{1}{3}-x^2& \text{when }t\in[1,T], \end{cases}
\label{eq:ex21-linear-solution-concrete-loads}
\end{equation}
see \cite[\nolinkurl{Example_3_sigma_tilde.mp4}]{SupplMat} for the animation.

\begin{figure}[H]\center
\includegraphics{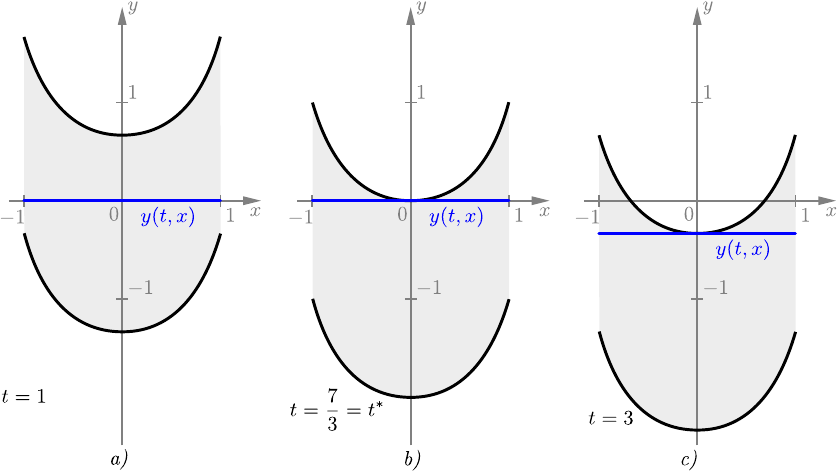}
\caption{\footnotesize With the particular data of Section \ref{ssect:ex21-concrete-data} the moving set $\mathcal{C}(t) \subset L^2_{{\bf C}^{-1}}(\Omega)$ consists of all functions constant a.e. on $\Omega =(-1,1)$  which fit in the gray area. {\it a)} illustrates the instant $t=1$, after which the body force load and the shape of $\mathcal{C}(t)$ will remain constant, and the evolution will be driven by the change of the displacement boundary condition. {\it b)} illustrates the instant $t=t^*$ of the beginning of the plastic deformation. {\it c)} illustrates the instant $t=3$ during the plastic deformation.
} 
\label{fig:ex21-moving-set}
\end{figure}

The sweeping process \eqref{eq:sp-elastoplastic} in this case is defined in terms of the regular inner product of $L^2(\Omega)$, its initial condition is
\begin{equation}
y(0,x) = \sigma_0(x) - \widetilde{\sigma}(0,x) \equiv 0 \qquad \text{for all }x\in \Omega
\label{eq:ex21-ic-data}
\end{equation}
and its moving set is \eqref{eq:ex21-moving-set}.  In particular for $t\in[1,T]$, the inequalities in \eqref{eq:ex21-moving-set} take the form
\[
-\frac{4}{3} - \frac{1}{2}(t-1)+x^2\leqslant c\leqslant \frac{2}{3}-\frac{1}{2}(t-1)+x^2,
\]
see Fig. \ref{fig:ex21-moving-set}.
Observe, that all the way through the interval $t\in\left[0, t^*\right],\, t^*=\frac{7}{3},$ the stationary constant-zero function is the solution of the sweeping process \eqref{eq:sp-elastoplastic} with initial condition \eqref{eq:ex21-ic-data}. Indeed, by plugging $c=0$ and \eqref{eq:ex21-linear-solution-concrete-loads} into \eqref{eq:ex21-moving-set} we can check the condition \ref{enum:abstract-sp-solution-def-constraint} of Definition \ref{def:abstract-sp}, and zero velocity clearly belongs to the normal cone \eqref{eq:nc-C-inv-def}, therefore the condition \ref{enum:abstract-sp-solution-def-nc} of Definition \ref{def:abstract-sp} also holds.

The time-instant $t^*=\frac{7}{3}$ is special, as that is when the {\it purely elastic deformation} ends and the {\it plastic deformation} of the rod starts.  Figs. \ref{fig:ex21-moving-set} {\it b)} and {\it c)} suggest, that for the interval $(t^*, T)$ the solution of the sweeping process is
\[
y(t,x) = -\frac{1}{2}(t-t^*) \qquad \text{with}\qquad \frac{d}{dt} y(t,x) =- \frac{1}{2}.
\]
It can be verified analytically that such function $\bm{y}(t)$ belongs to $\mathcal{C}(t)$ for $t\in[t^*, T]$. Moreover, for any $t\in[t^*, T]$ one can also see from Figs. \ref{fig:ex21-moving-set}  {\it b)} and {\it c)} that any constant function $c(x)$ that fits into $\mathcal{C}(t)$ has to be such that
\[
c(x)-y(t,x)\leqslant 0,
\] 
i.e. for the velocity we have
\[
-\left(\frac{d}{dt} y(t,x)\right)\,\left(c(x)-y(t,x)\right)\leqslant 0 \qquad \text{for any }x\in \Omega, \,t\in[t^*, T].
\]
This means that  condition \ref{enum:abstract-sp-solution-def-nc} of Definition \ref{def:abstract-sp} holds with the normal cone \eqref{eq:nc-C-inv-def} defined in terms of the regular integral inner product of $L^2(\Omega)$. To summarize, the solution to the sweeping process \eqref{eq:sp-elastoplastic} with the data of the current section is
\[
y(t,x) = \begin{cases}0&\text{when }t\in[0, t^*),\\
 -\frac{1}{2}(t-t^*)& \text{when } t\in [t^*, T],\end{cases} \qquad \text{with}\qquad \frac{d}{dt} y(t,x) =\begin{cases}0&\text{when }t\in[0, t^*),\\ - \frac{1}{2}& \text{when } t\in (t^*, T],\end{cases}
\]
see \cite[\nolinkurl{Example_3_y.mp4}]{SupplMat} for the animation.

By formula \eqref{eq:finding-stresses-epp} this gives the evolution of stress
\[
\sigma(t,x) = \begin{cases}
t\left(\frac{1}{3}-x^2\right) & \text{when } t\in [0,1)\\
\frac{1}{2}(t-1)+ \frac{1}{3}-x^2& \text{when }t\in[1,t^*)\\
1-x^2& \text{when } t\in [t^*, T]\end{cases}
\]
which, indeed, belongs to $\Sigma$, as required by perfect plasticity; see also \cite[\nolinkurl{Example_3_sigma.mp4}]{SupplMat} for the animation.

\subsection{Non-existence of a strain solution as a \texorpdfstring{$L^2$}{L2}-function}
\label{ssect:ex21-no-strain}
Now let us consider differential inclusion \eqref{eq:inclusion2-elastoplastic} and its component \eqref{eq:aepp-omega-inclusion} in particular. During time-interval $[0, t^*]$, when $\frac{d}{dt} \bm{y}=0$, we clearly can solve  \eqref{eq:aepp-omega-inclusion} with $\omega\equiv 0$ and \eqref{eq:inclusion2-elastoplastic} with $\frac{d}{dt} \bm{\varepsilon}=\frac{d}{dt} \bm{\widetilde{\sigma}}$, which agrees with the mechanical expectation of purely elastic evolution.

\begin{figure}[H]\center
\includegraphics{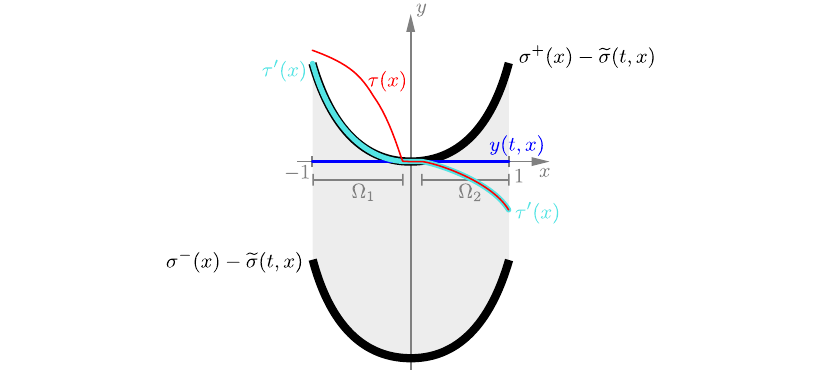}
\caption{\footnotesize An illustration of the argument for Observation \ref{obs:nonsupport}. 
} 
\label{fig:ex21-moving-set-projection}
\end{figure}

However, during time-interval $(t^*, T)$ we encounter the following phenomenon.
\begin{Observation}
\label{obs:nonsupport}
For all $t\in (t^*, T)$ in  \eqref{eq:aepp-omega-inclusion} we have
\begin{equation}
N^{{\bf C}^{-1}}_{\Sigma-\bm{\widetilde{\sigma}}(t)}(\bm{y}(t)) = \{0\}.
\label{eq:ex21-singleton-zero-normal-cone}
\end{equation}
\end{Observation}
\noindent This is the same as to say that the only function $\bm{\tau}\in L_{{\bf C}^{-1}}^2(\Omega)$ for which 
\[
{\rm proj}^{{\bf C}^{-1}}(\bm{\tau}, \Sigma-\bm{\widetilde{\sigma}}(t))  = {\bm y}(t) \qquad \text{i.e.} \qquad \bm{\tau} - {\bm y}(t) \in N^{{\bf C}^{-1}}_{\Sigma-\bm{\widetilde{\sigma}}(t)}(\bm{y}(t))
\]
is $\bm{y}(t)$ itself. We can show this with the following argument. Given $\bm{\tau}$ let $\bm{\tau'}$ be the projection
\begin{equation}
\bm{\tau'}= {\rm proj}^{{\bf C}^{-1}}(\bm{\tau}, \Sigma-\bm{\widetilde{\sigma}}(t)) = \underset{\bm{\zeta}\in \Sigma-\bm{\widetilde{\sigma}}(t)}{\rm arg\, min} \left(\intOm \frac{\left|\tau(x)-\zeta(x)\right|^2}{\bm{C}(x)}\, dx\right)^{\frac{1}{2}}
\label{eq:L2-projection}
\end{equation}
and assume that on a subset $\Omega_1\subset \Omega$ of positive measure $\tau(x)\neq y(t,x)$ and $\tau(x)$ violates the upper threshold of \eqref{eq:ex21-moving-set}, i.e. $\tau(x)>\sigma^+(x)-\widetilde{\sigma}(t,x)$, see Fig.~\ref{fig:ex21-moving-set-projection}. The projection $\tau'(x)$ has to coincide with the threshold $\sigma^+(x)-\widetilde{\sigma}(t,x)$ on $\Omega_1$, as any other function $\bm{\zeta}$ below the threshold, including $y(t,x)$, would give a strictly larger difference in the denominator of \eqref{eq:L2-projection} on $\Omega_1$, and a strictly larger value of the integral. We reason similarly for the lower threshold $\sigma^-(x)-\widetilde{\sigma}(t,x)$. 
\noindent Moreover, if $\tau(x)\neq y(t,x)$ and $\tau(x)$ stays between the thresholds on a subset $\Omega_2\subset \Omega$ of positive measure, then the projection $\tau'(x)$ has to coincide with $\tau(x)$ on $\Omega_2$, as it will give zero in the denominator of \eqref{eq:L2-projection}, which cannot be minimized further. Thus $\bm{\tau'}=\bm{y}(t)$ is only possible when $\tau(x) = y(t,x)$ a.e. on $\Omega$.

We have shown \eqref{eq:ex21-singleton-zero-normal-cone}, i.e. that the set $\Sigma-\bm{\widetilde{\sigma}}(t)$ in space $L^2(\Omega)$ {\it does not possess a nonzero supporting vector} at a point $\bm{y}(t)$, although it clearly lays on the boundary of the set.  This situation is in a striking contrast with a typical picture of convex analysis in finite dimensions (see Fig. \ref{fig:fig-normal-cone}), where the normal cone degenerates to the singleton-zero set if and only if $\bm{y}$ belongs to the interior of $\Sigma$, see Fig. \ref{fig:fig-normal-cone} {\it a)} specifically. In general, such points $\bm{y}$ with a degenerate normal cone fall into the category of {\it nonsupport points}, see \cite[p. 184]{Peressini1967}, \cite[p. 97]{Jeyakumar1992} , \cite[Def. 2.15, p. 21]{Borwein1992}, \cite[Def. 2.6 (b), p. 2544]{Borwein2003}, \cite[p. 155]{Mordukhovich2022}.

\begin{figure}[H]\center
\includegraphics{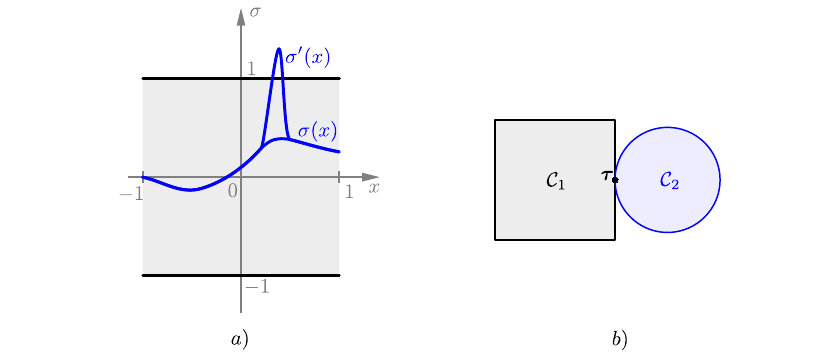}
\caption{\footnotesize {\it a)} An illustration of the argument for Observation \ref{obs:ex21-Slater-cq-fails}.  {\it b)} A typical finite-dimensional situation where $N_{\mathcal{C}_1}(\bm{\tau})+N_{\mathcal{C}_2}(\bm{\tau}) \subsetneqq N_{\mathcal{C}_1 \cap \mathcal{C}_2}(\bm{\tau})$.
}
\label{fig:misc-1}
\end{figure}

From Observation \ref{obs:nonsupport} it follows that for all $t\in (t^*, T)$ the first term in the intersection  \eqref{eq:aepp-omega-inclusion} is the singleton set
\begin{equation}
N^{{\bf C}^{-1}}_{\Sigma-\bm{\widetilde{\sigma}}(t)}(\bm{y}(t)) + \frac{d}{dt}\bm{y}(t) =\{0\}+ \frac{d}{dt}\bm{y}(t) =\left\{\frac{d}{dt}\bm{y}(t)\right\} =\left\{-\frac{1}{2}\right\},
\label{eq:ex21-with-data-nonzero-component}
\end{equation}
where $-\frac{1}{2}$ is meant as a constant function over $\Omega$. But, recall that  $\bm{y}(t)\in \mathcal{V}$ for all $t$ (see Remark \ref{rem:aepp-sp-in-V}), hence $\frac{d}{dt}\bm{y}(t)\in \mathcal{V}$ as well. As $\mathcal{U}$ and $\mathcal{V}$ are orthogonal subspaces, their only intersection is zero. But $\frac{d}{dt}\bm{y}(t)$ in \eqref{eq:ex21-with-data-nonzero-component} is not zero, from which we deduce the following.
\begin{Observation} For all $t\in(t^*, T]$ the right-hand sides of   \eqref{eq:aepp-omega-inclusion} and \eqref{eq:inclusion2-elastoplastic} are empty sets. Therefore, the problem \eqref{eq:sp-elastoplastic}--\eqref{eq:aepp-sp-moving-set} as a whole has no solution $(y, \varepsilon)\in W^{1, \infty}(I, \mathcal{H}_{{\bf C}^{-1}})\times W^{1, \infty}(I, \mathcal{H})$ for any $\bm{\varepsilon_0}$. By the statement of Theorem \ref{th:aepp-to-sweeping-process}, the problem of Definition \ref{def:aepp} also has no solution \eqref{eq:aepp-unknowns} with the data of Section \ref{ssect:ex21-concrete-data}.
\end{Observation}
This obviously contradicts the conclusion \ref{enum:th:aepp-well-posedness-nonempty-rhs} of Theorem \ref{th:aepp-well-posedness-safe-load-strict}. Therefore, Slater's condition \eqref{eq:safe-load-strict} must be failing. Indeed, any function $\bm{\sigma}\in \Sigma$ can be perturbed to another function $\bm{\sigma}'\notin \Sigma$ with 
\[
\|\bm{\sigma}-\bm{\sigma}'\|_{{\bf C}^{-1}} = \left(\intOm \frac{\left|\sigma(x)-\sigma'(x)\right|^2}{\bm{C}(x)}\, dx\right)^{\frac{1}{2}}
\]
arbitrarily small, see Fig. \ref{fig:misc-1} {\it a)}. I.e. the following fact takes place.
\begin{Observation} \label{obs:ex21-Slater-cq-fails} In the setting of the model of the continuous rod 
\[
{\rm int}\, \Sigma =\varnothing
\]
and Slater's condition \eqref{eq:safe-load-strict} is not satisfied.
\end{Observation}

It may appear somewhat surprising that the sweeping process \eqref{eq:sp-elastoplastic} with the set \eqref{eq:aepp-sp-moving-set} has a well-defined solution, and yet the normal cone \eqref{eq:ex21-singleton-zero-normal-cone} is degenerate and the differential inclusion \eqref{eq:inclusion2-elastoplastic} is ill-posed. To understand the situation, revisit the proof of Theorem \ref{th:aepp-well-posedness-safe-load-strict} \ref{enum:th:aepp-well-posedness-nonempty-rhs} and observe that the failure of the condition \eqref{eq:safe-load-strict} means that the additivity of normal cones \eqref{eq:nc-additivity-in-aepp-proof} is no longer guaranteed. And, indeed, we merely have one-sided relation
\begin{equation}
N^{{\bf C}^{-1}}_{\Sigma-\bm{\widetilde{\sigma}}(t)}(\bm{y}(t))+ N^{{\bf C}^{-1}}_{\mathcal{V}}(\bm{y}(t)) = \{0\} + \mathcal{U}\subsetneqq N^{{\bf C}^{-1}}_{\left(\Sigma-\bm{\widetilde{\sigma}}(t)\right)\cap \mathcal{V}}(\bm{y}(t)),
\label{eq:ex21-nc-additivity-failure}
\end{equation}
cf. \eqref{eq:nc-subadditivity}, because we already had found the vector $-\frac{d}{dt}\bm{y}(t)$, which belongs the right-hands side, but not to the left-hand side, see the argument we discussed after \eqref{eq:ex21-with-data-nonzero-component}.
We conclude with the following.
\begin{Observation}
\label{obs:ex21-cq-importance}
When the additivity of normal cones \eqref{eq:nc-additivity-in-aepp-proof} fails to hold with equality, in the  problem of Definition \ref{def:aepp} and Theorem \ref{th:aepp-to-sweeping-process}  it may be impossible to extract the evolution of $\varepsilon$ and $\varepsilon_{\rm p}$ from the evolution of $\sigma,\, \varepsilon_{\rm el}$ and $y$. Condition \eqref{eq:safe-load-strict} serves as a {\it constraint qualification} to ensure such additivity, see Proposition \ref{prop:nc-abstract-properties} \ref{enum:prop:nc-abstract-properties-nc-subadditivity}, but condition \eqref{eq:safe-load-strict} fails when $\Sigma$ is described via inequalities on the state variable in a space with an integral norm. If  for a similar problem we could find a less demanding constraint qualification which ensures the additivity of normal cones then we could solve the problem in full.   
\end{Observation}
Note, that a failure of constraint qualification \eqref{eq:nontrivial-intersection-condition} can lead to the loss of additivity of normal cones in a finite-dimensional setting as well (see Fig. \ref{fig:misc-1} {\it b)}), but not when the involved sets are polyhedral (see \cite[Corollary 23.8.1, p. 224]{Rockafellar1970}).

\subsection{Concluding remarks on the continuum elastic-perfectly plastic model}
As we have demonstrated, the lack of additivity of the normal cones explains why the problem of quasi-static evolution in elasticity-perfect plasticity, formulated in Definition \ref{def:aepp} cannot be solved in general through the means of Theorem \ref{th:aepp-to-sweeping-process} when we search for the unknown variable
\[
\bm{\varepsilon}\in \mathcal{H}= L^2(\Omega),
\]
although this approach works in a finite-dimensional setting. In fact, J.-J. Moreau encountered this issue and shown in \cite{Moreau1976} that for the problem of elastic-perfectly plastic rod one must search for a displacement field $\bm{u}$ and the corresponding strain field $\bm{\varepsilon}$  so that
\[
\bm{u}\text{ has bounded variation on }\overline{\Omega}, \qquad \bm{\varepsilon}=\frac{\partial}{\partial x} \bm{u} \text{ is a bounded Borel signed measure on } \overline{\Omega},
\]
where the derivative is meant in the appropriate sense, see e.g \cite[Def. 3.1, p. 117]{Ambrosio2000}. Unlike the approach of this text with a Hilbert space $\mathcal{H}$ for both strain and stress, Moreau interpreted a measure-valued strain as a functional applied to $\bm{\sigma}\in C(\overline{\Omega})$ by Riesz theorem \cite[Th. 5.15, pp. 347]{Pugachev1999}, \cite[Th. 1.196, p. 126]{Fonseca2007}. In such a case the interior of $\Sigma$ in Slater's condition \eqref{eq:safe-load-strict} is taken in the sense of the uniform norm on continuous functions, so that the condition becomes useful and allows to solve the problem. Specifically, for the continuous rod problem of Section \ref{ssect:continuum_rod_epp_general} the equivalent safe load condition \eqref{eq:safe-load-forces-strict} in $C(\overline{\Omega})$ is 
\begin{multline}
\text{for each}\quad t\in I \quad\text{there exists}\quad c\in \mathbb{R} \quad \text{such that for all}\quad x\in \overline{\Omega}\quad \text{we have}\\
\sigma^-(x)<c+ \omega_0 \intOm \frac{1}{\bm{C}(z)}\int\limits_a^z F(y)\, dy\, dz  -\int\limits_a^x F(y)\, dy < \sigma^+(x)
\label{eq:safe-load-continuous-topology}
\end{multline}
given upper semicontinuous on $\overline{\Omega}$ function $\sigma^-$ and lower semicountinuous $\sigma^+$. Notice, that this implies a uniform lower estimate $\delta$ on the size of the elastic range, i.e.
\[
\sigma^+(x)-\sigma^-(x)\geqslant \delta>0 \qquad \text{for all }x\in \overline{\Omega}
\] 
on top of \eqref{eq:local_yield_limits}. The safe load conditions of the type \eqref{eq:safe-load-continuous-topology} also lead to improved regularity in the related problems, such as Norton-Hoff matrtials (viscoelasticity, also called creep), see \cite[(4.2)]{Bensoussan1993}, \cite[(3.3)]{Bensoussan1996}.

When passing from one-dimensional domain $\Omega$ to two- and three-dimensional $\Omega\subset \mathbb{R}^n$ the displacement field $\bm{u}$ must be considered with the values in the so-called space of bounded deformation $BD(\Omega)$, i.e. 
\[
\bm{u}\in L^1(\Omega, \mathbb{R}^n), \qquad \bm{\varepsilon}=\frac{1}{2}\left(D\bm{u}+(D\bm{u})^T\right)\text{ is a bounded Borel vector measure on }\Omega.
\]
when $\Omega$ is bounded. We will leave this theory beyond the scope of the current text, but recommend the references from Section \ref{sect:intro}. 

\section{Constraint qualifications and elastoplasticity}
\label{sect:constraint-qualification-and-hardening}
\subsection{Constraint qualifications}
\label{ssect:constraint-qualifications}
Let us focus now on the following property (we repeat Proposition \ref{prop:nc-abstract-properties} \ref{enum:prop:nc-abstract-properties-nc-subadditivity} for presentation purposes):
\begin{Proposition}
Let $\mathcal{H}$ be a Hilbert space, $\mathcal{C}_1, \mathcal{C}_2\subset \mathcal{H}$ be closed, convex, nonempty sets and $\bm{x}\in \mathcal{C}_1\cap \mathcal{C}_2$. Then
\begin{equation}
N_{\mathcal{C}_1}(\bm{x}) + N_{\mathcal{C}_2}(\bm{x}) \subseteqq N_{\mathcal{C}_1 \cap \mathcal{C}_2}(\bm{x}).
\label{eq:nc-subadditivity-general}
\end{equation}
\end{Proposition}
In the proof of Theorem \ref{th:aepp-to-sweeping-process} we used \eqref{eq:nc-subadditivity-general} to justify that the evolution of stress in problems with elastoplasticity can be found from the corresponding sweeping process. We also observed, that to extract the evolution of the plastic component and of the total strain from the evolution of stresses, we need the  equality of sets in \eqref{eq:nc-subadditivity-general}, which is the property called {\it additivity of normal cones}:
\begin{equation}
N_{\mathcal{C}_1}(\bm{x}) + N_{\mathcal{C}_2}(\bm{x}) = N_{\mathcal{C}_1 \cap \mathcal{C}_2}(\bm{x}).
\label{eq:nc-additivity-general}
\end{equation}
To justify such equality we can use conditions on $\mathcal{C}_1$ and $\mathcal{C}_2$, which are called {\it constraint qualifications}. Constraint qualifications are usually discussed in the context of optimization problems, written in a specific form, and we reformulate our question of additivity of the normal cones in terms of {\it Fenchel-Rockafellar dual problems}.
\begin{Proposition} 
\label{prop:from-nc-to-dual-problems}
Let $\mathcal{H}$ be a Hilbert space, $\mathcal{C}_1, \mathcal{C}_2\subset \mathcal{H}$ be closed, convex, nonempty sets and $\bm{x}\in \mathcal{C}_1\cap \mathcal{C}_2$. Let 
\begin{equation}
\bm{v} \in N_{\mathcal{C}_1 \cap \mathcal{C}_2}(\bm{x}).
\label{eq:nc-additivity-to-primal-dual-v-taken-1}
\end{equation}
and consider two optimization problems with optimal solutions, respectively, $p^*$ and $d^*$:
\begin{align}
p^*&=- \delta^*_{\mathcal{C}_1\cap \mathcal{C}_2}(\bm{v})&&=\inf_{\bm{c}\in \mathcal{H}}\left(\delta_{\mathcal{C}_1}(\bm{c})+\delta_{\mathcal{C}_2}(\bm{c})-\la \bm{v},\bm{c}\ra\right),&
\label{eq:abstract-primal-problem}
\\[2mm]
d^*&= - \inf_{\bm{y}\in \mathcal{H}} \left(\delta^*_{\mathcal{C}_1}(\bm{y})+ \delta^*_{\mathcal{C}_2}(\bm{v}-\bm{y}) \right) &&= \sup_{\bm{y}\in \mathcal{H}}\left(-\delta^*_{\mathcal{C}_1}(\bm{y})- \delta^*_{\mathcal{C}_2}(\bm{v}-\bm{y}) \right),& 
\label{eq:abstract-dual-problem}
\end{align}
where $\delta$ and $\delta^*$ denote the indicator function \eqref{eq:def-indicator-function} and the support function \eqref{eq:def-support-function}, respectively.
Then
\begin{equation}
\bm{v} \in N_{\mathcal{C}_1}(\bm{x}) + N_{\mathcal{C}_2}(\bm{x})
\label{eq:nc-additivity-to-primal-dual-v-taken-2}
\end{equation}
if and only if the following equality, called strong duality, holds:
\begin{equation}
p^*=d^* \text{ and the optimal value } d^* \text{ is achieved in \eqref{eq:abstract-dual-problem} at some }\bm{y^*}\in \mathcal{H}. 
\label{eq:nc-additivity-to-primal-dual-strong-duality}
\end{equation}
\end{Proposition}
\begin{Remark} Optimization problems \eqref{eq:abstract-primal-problem} and \eqref{eq:abstract-dual-problem} are also sometimes written (see e.g. \cite{Attouch1986}) in the following forms, respectively,
\begin{align*}
-p^*&= (\delta_{\mathcal{C}_1}+\delta_{\mathcal{C}_2})^*(\bm{v}),\\
-d^*&= (\delta^*_{\mathcal{C}_1}\square \,\delta^*_{\mathcal{C}_2})(\bm{v}),
\end{align*}
where $f^*$ is the {\it Fenchel conjugate} (or {\it  Legendre transform}) of a function $f:\mathcal{H}\to \mathbb{R}\cup\{+\infty\}$, defined as 
\[
f^*: \mathcal{H}\to \mathbb{R}\cup\{+\infty\}, \qquad f^*(\bm{y}) = \sup_{\bm{x}\in\mathcal{H}}\left(\la\bm{y}, \bm{x}\ra - f(\bm{x})\right),
\]
and $f_1\square\, f_2$ is the {\it infimal convolution} of two functions $f_1, f_2:\mathcal{H}\to \mathbb{R}\cup\{+\infty\}$, defined as 
\[
(f_1\square \,f_2): \mathcal{H}\to \mathbb{R}\cup\{+\infty\}, \qquad (f_1\square \,f_2)(\bm{v}) = \inf_{\bm{y}\in\mathcal{H}}\left(f_1(\bm{y})+f_2(\bm{v}-\bm{y}))\right).
\]
\end{Remark}

\noindent Before we prove Proposition \ref{prop:from-nc-to-dual-problems}, let us verify the following statement.
\begin{Proposition} \label{prop:weak-duality} {\bf (weak duality)}
Let $\mathcal{H}$ be a Hilbert space and $\mathcal{C}_1, \mathcal{C}_2\subset \mathcal{H}$ be closed, convex, nonempty sets such that $\mathcal{C}_1\cap \mathcal{C}_2$ is also nonempty. For a given $\bm{v}\in \mathcal{H}$ let $p^*, d^*\in \mathbb{R}\cup \{-\infty\}$ be defined by \eqref{eq:abstract-primal-problem}--\eqref{eq:abstract-dual-problem}. Then
\[
p^*\geqslant d^*.
\]
\end{Proposition}
\noindent{\bf Proof of Proposition \ref{prop:weak-duality}.}
Clearly, for any $\bm{y} \in \mathcal{H}$ and any $\bm{c}\in \mathcal{C}_1\cap \mathcal{C}_2$
\[
\la \bm{y}, \bm{c}\ra + \la\bm{v}- \bm{y}, \bm{c}\ra  =\la \bm{v},\bm{c} \ra.
\]
Since $\mathcal{C}_1\cap\mathcal{C}_2$ is a smaller set than $\mathcal{C}_1$ and $\mathcal{C}_2$ each, 
\[
\sup_{\bm{c}\in \mathcal{C}_1} \la \bm{y}, \bm{c}\ra +\sup_{\bm{c}\in \mathcal{C}_2} \la \bm{v}- \bm{y}, \bm{c}\ra \geqslant \sup_{\bm{c}\in \mathcal{C}_1 \cap \mathcal{C}_2} \la \bm{v}, \bm{c}\ra,
\]
i.e. for any $\bm{y}\in \mathcal{H}$
\[
\delta^*_{\mathcal{C}_1}(\bm{y}) +\delta^*_{\mathcal{C}_2}(\bm{v}-\bm{y}) \geqslant \delta^*_{\mathcal{C}_1 \cap \mathcal{C}_2}(\bm{v}).
\]
Find the infimum w.r. to $\bm{y}$:
\[
\inf_{\bm{y}\in \mathcal{H}}  \left(\delta^*_{\mathcal{C}_1}(\bm{y}) +\delta^*_{\mathcal{C}_2}(\bm{v}-\bm{y})\right) \geqslant \delta^*_{\mathcal{C}_1 \cap \mathcal{C}_2}(\bm{v})
\] 
and multiply by $-1$ to obtain the claimed inequality. $\blacksquare$

\noindent {\bf Proof of Proposition \ref{prop:from-nc-to-dual-problems}.} \underline{Part ``if''.} From the statement of the Proposition we are given $\bm{x}\in \mathcal{C}_1\cap \mathcal{C}_2, \bm{v}\in \mathcal{H}$ such that \eqref{eq:nc-additivity-to-primal-dual-v-taken-1} holds, and $\bm{y^*}\in \mathcal{H}$, at which the optimal value $d^*$ is achieved in \eqref{eq:abstract-dual-problem}. Since $p^*=d^*$ we have 
\begin{equation*}
0=p^*-d^* = - \delta^*_{\mathcal{C}_1\cap \mathcal{C}_2}(\bm{v}) + \delta^*_{\mathcal{C}_1}(\bm{y^*}) + \delta^*_{\mathcal{C}_2}(\bm{v}-\bm{y^*}) =-\la\bm{v},\bm{x}\ra + \sup_{\bm{c}\in \mathcal{C}_1}\,\la \bm{y^*}, \bm{c} \ra+ \sup_{\bm{c}\in \mathcal{C}_2}\,\la \bm{v}-\bm{y^*}, \bm{c} \ra,
\end{equation*}
where the last equality is due to \eqref{eq:nc-additivity-to-primal-dual-v-taken-1} and the basic property of the normal cone, see Proposition \ref{prop:nc-equivalent} {\it i), iii)}. Therefore,
\begin{equation*}
\la \bm{y^*},\bm{x}\ra+ \la\bm{v} - \bm{y^*},\bm{x}\ra= \sup_{\bm{c}\in \mathcal{C}_1}\,\la \bm{y^*}, \bm{c} \ra+ \sup_{\bm{c}\in \mathcal{C}_2}\,\la \bm{v}-\bm{y^*}, \bm{c} \ra,
\end{equation*}
Observe that, since $\bm{x}\in \mathcal{C}_1$ and $\bm{x}\in \mathcal{C}_2$ simultaneously, this is only possible when
\[
\la \bm{y^*},\bm{x}\ra  = \sup_{\bm{c}\in \mathcal{C}_1}\,\la \bm{y^*}, \bm{c} \ra \qquad\text{and} \qquad  \la\bm{v}- \bm{y^*},\bm{x}\ra = \sup_{\bm{c}\in \mathcal{C}_2}\,\la \bm{v}-\bm{y^*}, \bm{c} \ra,
\]
i.e.
\begin{equation}
\la \bm{y^*},\bm{x}\ra  = \delta^*_{\mathcal{C}_1}(\bm{y^*}) \qquad\text{and} \qquad  \la\bm{v}- \bm{y^*},\bm{x}\ra =  \delta^*_{\mathcal{C}_2}(\bm{v}-\bm{y^*}).
\label{eq:nc-additivity-to-primal-dual-proof-separate-support-funcs}
\end{equation}
Again, by  Proposition \ref{prop:nc-equivalent} {\it i), iii)} this means that
\begin{equation}
\bm{y^*}\in N_{\mathcal{C}_1}(\bm{x})\qquad\text{and} \qquad \bm{v}-\bm{y^*}\in N_{\mathcal{C}_2}(\bm{x}),
\label{eq:nc-additivity-to-primal-dual-proof-separate-nc}
\end{equation}
from which \eqref{eq:nc-additivity-to-primal-dual-v-taken-2} follows.
\newline\noindent \underline{Part ``only if''.} Assume now that both \eqref{eq:nc-additivity-to-primal-dual-v-taken-1} and \eqref{eq:nc-additivity-to-primal-dual-v-taken-2} hold for some $\bm{v}, \bm{x}$. From \eqref{eq:nc-additivity-to-primal-dual-v-taken-2} we deduce the existence of $\bm{y^*}$, such that \eqref{eq:nc-additivity-to-primal-dual-proof-separate-nc} holds, which yields \eqref{eq:nc-additivity-to-primal-dual-proof-separate-support-funcs} and 
\[
\la \bm{y^*},\bm{x}\ra+ \la\bm{v} - \bm{y^*},\bm{x}\ra= \delta^*_{\mathcal{C}_1}(\bm{y^*}) +  \delta^*_{\mathcal{C}_2}(\bm{v}-\bm{y^*}),
\]
i.e.
\[
\la\bm{v}, \bm{x}\ra = \delta^*_{\mathcal{C}_1}(\bm{y^*}) +  \delta^*_{\mathcal{C}_2}(\bm{v}-\bm{y^*}).
\]
On the other hand, \eqref{eq:nc-additivity-to-primal-dual-v-taken-1} yields
\[
\la\bm{v}, \bm{x}\ra=\delta^*_{\mathcal{C}_1\cap \mathcal{C}_2}(\bm{v}),
\]
hence we have 
\[
\delta^*_{\mathcal{C}_1\cap \mathcal{C}_2}(\bm{v}) = \delta^*_{\mathcal{C}_1}(\bm{y^*}) +  \delta^*_{\mathcal{C}_2}(\bm{v}-\bm{y^*}),
\]
i. e.
\[
p^* = -\delta^*_{\mathcal{C}_1}(\bm{y^*}) -  \delta^*_{\mathcal{C}_2}(\bm{v}-\bm{y^*}).
\]
As we have shown in Proposition \ref{prop:weak-duality}, the supremum in \eqref{eq:abstract-dual-problem} cannot be larger than $p^*$, and we have presented the specific $\bm{y}=\bm{y}^*$ at which the function inside the supremum in \eqref{eq:abstract-dual-problem} attains the value $p^*$. This is exactly what we claimed in \eqref{eq:nc-additivity-to-primal-dual-strong-duality}. $\blacksquare$

Proposition \ref{prop:from-nc-to-dual-problems} gives the connection between the additivity of normal cones and the dual optimization problems. It allows to guarantee the additivity of normal cones when at least one of the constraint qualifications below holds.

\begin{Proposition}
\label{prop:constraint-qualifications}
Let $\mathcal{H}$ be a Hilbert space, $\mathcal{C}_1, \mathcal{C}_2\subset \mathcal{H}$ be closed, convex, nonempty sets. Assume that at least one of the conditions listed below holds. Then for any $\bm{x}, \bm{v}$ as in \eqref{eq:nc-additivity-to-primal-dual-v-taken-1} the strong duality \eqref{eq:nc-additivity-to-primal-dual-strong-duality} holds , and,  therefore, the additivity of normal cones \eqref{eq:nc-additivity-general} holds for any $\bm{x}\in \mathcal{C}_1\cap\mathcal{C}_2$.
\begin{enumerate}[{\it i)}]
\item \label{cq-Slater-I} (Slater's constraint qualification I. \cite[p. 239]{Aubin2000}, \cite[Section 2.f, p. 206]{Moreau1973}, \cite[Lemma 1 (b), pp.~4--5]{Kunze2000})
\[
{\mathcal{C}_1}\cap {{\rm int}\, \mathcal{C}_2}\neq \varnothing.
\]
\item \label{cq-Slater-II} (Slater's constraint qualification II.  \cite[Th. 10.5.3, pp. 239--240]{Aubin2000}, \cite[Th. 2.1, p. 926]{Gowda1990})
\[
0\in {\rm int} \left({\mathcal{C}_1} - {\, \mathcal{C}_2}\right).
\]
\item  \label{cq-Rockafellar} (Rockafellar's constraint qualification \cite[Th. 17, Th. 18, pp. 41--42]{Rockafellar1974},  \cite[Th. 2.2, p. 926]{Gowda1990})
\[
\left\{\lambda \left(\bm{c_1}-\bm{c_2}\right): \lambda \geqslant 0, \bm{c_1}\in \mathcal{C}_1, \bm{c_2}\in \mathcal{C}_2 \right\}=\mathcal{H}.
\]

\item \label{cq-GCQ} (Attouch--Brezis's constraint qualification \cite{Attouch1986}, also known as the ``general constraint qualification'' \cite[Th. 3.5, Def. 3.3, Prop. 3.4, pp. 930--932]{Gowda1990})
\[
\left\{\lambda \left(\bm{c_1}-\bm{c_2}\right): \lambda \geqslant 0, \bm{c_1}\in \mathcal{C}_1, \bm{c_2}\in \mathcal{C}_2 \right\}\quad \text{is a closed linear subspace of }\mathcal{H}.
\]
\end{enumerate}
\end{Proposition}
\noindent We should make several comments on Proposition \ref{prop:constraint-qualifications}. First, note, that \ref{cq-Slater-I}$\Rightarrow$\ref{cq-Slater-II}$\Rightarrow$\ref{cq-Rockafellar}$\Rightarrow$\ref{cq-GCQ}. Slater's constraint qualifications have their name from the work \cite{Slater1950}, where such type of conditions was introduced in 1950. 
Constraint qualifications \ref{cq-Rockafellar} and \ref{cq-GCQ} are related, respectively, to the concepts of the {\it core} and  the {\it strong quasi-relative interior} of a set. They are generalizations to Banach spaces of the relative interior of a set. A reader, interested in the topic, may also find useful the following works: \cite{Jeyakumar1992},\cite{Jeyakumar1992b},\cite{Borwein1992}, \cite{Zalinescu2002}, \cite{Borwein2003},\cite{Grad2010}, \cite{Daniele2014}.

\subsection{Attempting to apply constraint qualifications to the continuum elastic-perfectly plastic model}
\label{ssect:ex2-pp-cq}
In Section \ref{ssect:ex21-concrete-data} we have verified directly that in the continuum elastic-perfectly plastic model the additivity of normal cones can fail (see \eqref{eq:ex21-nc-additivity-failure}), which, in turn, leads to the impossibility of solving \eqref{eq:aepp-omega-inclusion} and \eqref{eq:inclusion2-elastoplastic}. We have also checked, that the constraint qualification with the strongest condition in our list, Proposition \ref{prop:constraint-qualifications} \ref{cq-Slater-I} fails, see Observation \ref{obs:ex21-Slater-cq-fails}.    For completeness, let us confirm the constraint qualification with the weakest condition (Proposition \ref{prop:constraint-qualifications} \ref{cq-GCQ}) also fails in the same situation.

Recall, that we apply the constraint qualifications to the pair of sets
\[
\mathcal{C}_1=\Sigma - \bm{\widetilde{\sigma}}(t), \qquad \mathcal{C}_2=\mathcal{V},
\]
in the Hilbert space 
\[
\mathcal{H} = L^2_{{\bf C}^{-1}}(\Omega)= L^2(\Omega), \qquad \Omega=(-1, 1),
\]
where $\Sigma$ is given by \eqref{eq:ex21-Sigma-def} with $\sigma^+\equiv-\sigma^-\equiv 1$, vector $\bm{\widetilde{\sigma}}(t)$ is given by \eqref{eq:ex21-linear-solution-concrete-loads} (and we are interested in $t\in\left[t^*, T\right] = \left[\frac{7}{3}, 3\right]$ in particular), and $\mathcal{V}$ is the subspace of constant functions \eqref{eq:ex21-space-V}, see Fig. \ref{fig:ex21-moving-set}.

\begin{figure}[H]\center
\includegraphics{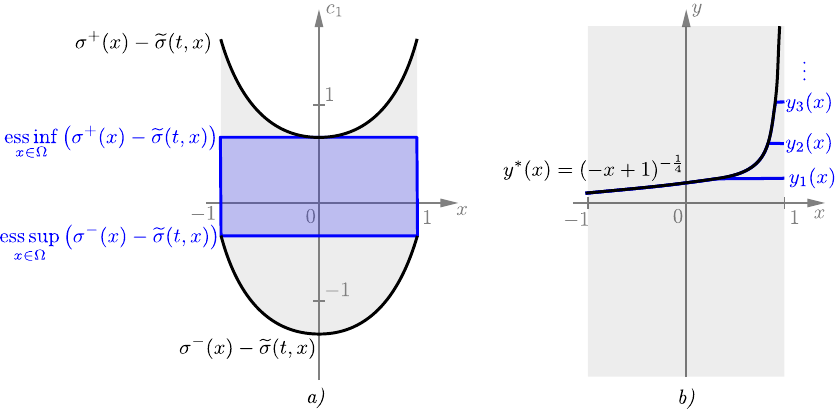}
\caption{\footnotesize {\it a)} As long as the blue rectangle is nondegenerate, we have \eqref{eq:ex21-pp-nondegeneracy} and \eqref{eq:ex21-qc-generate-L-infty} being true,  because any function $\bm{v}\in L^\infty(\Omega)$ can be rescaled and translated to fit in the blue rectangle, cf \eqref{eq:ex21-cq-L-infty-proof-affine-transformation}. {\it  b)} An explicit example showing that $L^\infty(\Omega)$ is not closed in the $L^2$-norm. Take $\bm{y^*}\in L^2(\Omega)$ as $y^*(x) = (-x+1)^{-\frac{1}{4}}$, and approximate it by $y_i(x) = \min(y^*(x), i),\,i \in \mathbb{N}$. We have $\|\bm{y^*}-\bm{y_}i\|_{L^2}(\Omega)\to 0, \bm{y_i}\in L^\infty(\Omega)$, yet $\bm{y^*}\notin L^\infty(\Omega)$.
} 
\label{fig:ex21-pp-gcq-fail}
\end{figure}

\begin{Lemma}
\label{lemma:ex21-qc-generate-L-infty}
With the data recalled above (in the current Section \ref{ssect:ex2-pp-cq}), we have
\begin{equation}
\left\{\lambda \left(\bm{c_1}-\bm{c_2}\right): \lambda \geqslant 0, \bm{c_1}\in \mathcal{C}_1, \bm{c_2}\in \mathcal{C}_2 \right\}  = L^\infty(\Omega).
\label{eq:ex21-qc-generate-L-infty}
\end{equation}
\end{Lemma}
\noindent{\bf Proof.} Indeed, in the particular case \eqref{eq:ex21-linear-solution-concrete-loads}, as well as in the general case \eqref{eq:ex21-linear-solution}, $\bm{\widetilde{\sigma}}(t)$ is an absolutely continuous function of $x$ on $\overline{\Omega}$ for any fixed $t$. Hence 
\[
\bm{\widetilde{\sigma}}(t)\in  L^\infty(\Omega),
\]
\[
\Sigma - \bm{\widetilde{\sigma}}(t)\subset  L^\infty(\Omega)\]
due to the construction of $\Sigma$ in \eqref{eq:ex21-Sigma-def}, and, of course,
\[
\left\{\lambda \left(\bm{c_1}-c_2\right): \lambda \geqslant 0, \bm{c_1}\in \Sigma - \bm{\widetilde{\sigma}}(t), c_2\in \mathbb{R} \right\} \subset L^\infty(\Omega),
\]
cf. \eqref{eq:ex21-qc-generate-L-infty}. On the other hand, as long as
\begin{equation}
\underset{x\in \Omega}{\rm ess\, inf} \left(\sigma^+(x)-\widetilde{\sigma}(t,x)\right)> \underset{x\in \Omega}{\rm ess\, sup} \left(\sigma^-(x)-\widetilde{\sigma}(t,x)\right),
\label{eq:ex21-pp-nondegeneracy}
\end{equation}
i.e. as long as the rectangle on Fig. \ref{fig:ex21-pp-gcq-fail} {\it a)} does not degenerate to a horizontal segment or an empty set, any $\bm{v}\in L^\infty(\Omega)$ can be represented as 
\begin{equation}
{\bm v}  =\lambda(\bm{c_1}-c_2)
\label{eq:ex21-cq-L-infty-proof-affine-transformation}
\end{equation}
for some $\lambda>0, c_2\in \mathbb{R}$ and 
\begin{multline*}
\bm{c_1}\in \left\{\bm{c}\in L^\infty(\Omega): \vphantom{\underset{x\in \Omega}{\rm ess\, inf}} \text{for a.a. }y\in \Omega \text{ we have }\right.\\ \left.\underset{x\in \Omega}{\rm ess\, sup} \left(\sigma^-(x)-\widetilde{\sigma}(t,x)\right) \leqslant c(y)\leqslant \underset{x\in \Omega}{\rm ess\, inf} \left(\sigma^+(x)-\widetilde{\sigma}(t,x)\right) \right\} \subset \Sigma - \bm{\widetilde{\sigma}}(t).
\end{multline*}
This finishes the proof of \eqref{eq:ex21-qc-generate-L-infty}. $\blacksquare$

Recall, however, that $L^\infty(\Omega)$ is not a closed subspace of $L^2(\Omega)$, see Fig.   \ref{fig:ex21-pp-gcq-fail} {\it b)}  for an illustration. Thus Lemma \ref{lemma:ex21-qc-generate-L-infty} means that even the constraint qualification with the weakest condition (Proposition \ref{prop:constraint-qualifications} \ref{cq-GCQ}) fails for the  the continuum elastic-perfectly plastic model. This was expected though, as we have already shown that \eqref{eq:ex21-nc-additivity-failure} does not hold as equality.

\subsection{The abstract framework, suitable for elasticity-hardening plasticity}
\begin{figure}[H]\center
\includegraphics{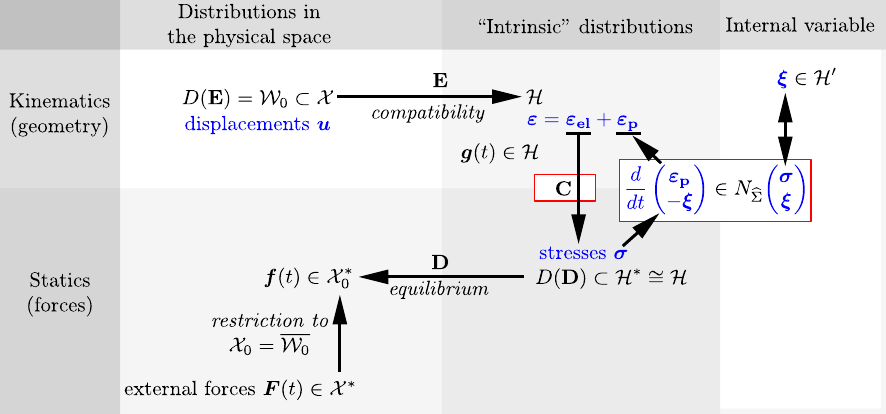}
\caption{\footnotesize Schematic representation of the problem of Definition \ref{def:aehp}. The unknown variables are indicated by blue color. In the problem of Definition \ref{def:aehp} we are only looking for the unknowns $\varepsilon, \varepsilon_{\rm el}, \varepsilon_{\rm p}, \sigma$ and $\xi$. Red rectangles indicate the constitutive relations.
} 
\label{fig:elasticity-hp-scheme}
\end{figure}

Let us modify the Definition \ref{def:aepp} of an abstract plasticity problem to allow for a different kind of constraint $\Sigma$.
\begin{Definition} 
\label{def:aehp}
Let spaces $\mathcal{H}, \mathcal{X}, \mathcal{W}_0$, operators ${\bf E}, {\bf D}, {\bf C}$ and functions $\bm{g}, \bm{f}$ be as in Section \ref{ssect:adjoint-operatorsED} and Definition \ref{def:ae}, but we require $\mathcal{H}$ to be separable. In addition, let us be given another separable Hilbert space $\mathcal{H}'$, which we call the {\it space of internal variables}. We denote the entire configuration space by
\[
\widehat{\mathcal{H}} =\mathcal{H}\times\mathcal{H}'.
\]

Let us be given a closed convex nonempty set $\widehat{\Sigma}\subset \widehat{\mathcal{H}}$, possibly unbounded. We say that the unknown variables
\begin{equation}
\varepsilon, \varepsilon_{\rm el}, \varepsilon_{\rm p}, \sigma \in W^{1, \infty}(I, \mathcal{H}), \qquad \xi\in W^{1, \infty}(I, \mathcal{H}')
\label{eq:aehp-unknowns}
\end{equation}
solve the {\it abstract problem of quasi-static evolution in elasticity-monotone plasticity} if they satisfy
\begin{align}
\bm{\varepsilon}& \in {\rm Im}\, {\bf E}+\bm{g} , \label{eq:aehp-1} \tag{EMP1}\\
\bm{\varepsilon}& = \bm{\varepsilon_{\bf el}}+\bm{\varepsilon_{\bf p}}, \label{eq:aehp-2} \tag{EMP2}\\
\bm{\sigma}& = {\bf C} \, \bm{\varepsilon_{\bf el}}, \label{eq:aehp-3} \tag{EMP3}\\
\frac{d}{dt} \begin{pmatrix}\bm{\varepsilon_{\bf p}}\\ -\bm{\xi} \end{pmatrix}& \in N_{\widehat{\Sigma}} \begin{pmatrix}\bm{\sigma}\\ \bm{\xi}\end{pmatrix}, \label{eq:aehp-4} \tag{EMP4}\\ 
\bm{\sigma}\in D(\bf D)\quad \text{ and }\quad {\bf D}\, \bm{\sigma}&=\bm{f}, \label{eq:aehp-5} \tag{EMP5}
\end{align}
and the initial condition
\[
(\bm{\varepsilon}(0), \bm{\varepsilon}_{\bf  el}(0), \bm{\varepsilon}_{\bf  p}(0), \bm{\sigma}(0), \bm{\xi}(0)) =(\bm{\varepsilon}_{\bf 0}, \bm{\varepsilon}_{\bf  el0}, \bm{\varepsilon}_{\bf  p0}, \bm{\sigma}_{\bf 0}, \bm{\xi}_{\bf 0})
\]
with some given right-hand side from $\mathcal{H}^4\times \mathcal{H}'$ satisfying \eqref{eq:aehp-1}--\eqref{eq:aehp-3}, \eqref{eq:aehp-5} at $t=0$ and 
\begin{equation*}
\begin{pmatrix}\bm{\sigma}_{\bf 0}\\ \bm{\xi}_{\bf 0}\end{pmatrix} \in \widehat{\Sigma}.
\label{eq:aehp-sigma-ic-compatible}
\end{equation*}
\end{Definition}
The problem of Definition \ref{def:aehp} can be represented as the diagram of Fig. \ref{fig:elasticity-hp-scheme}. Similarly to Theorem \ref{th:aepp-to-sweeping-process}, we will convert the problem to a sweeping process and a differential inclusion with a known right-hand side. Both elasticity-perfect plasticity and elasticity-hardening plasticity can be modeled via the framework of Definition \ref{def:aehp}, and we use the term elasticity-{\it monotone} plasticity due to the stress-strain relation being monotone (non-decreasing) in both cases throughout combined elastic and plastic regimes, see e.g. \cite[Fig.~3.2 b), p.~43 vs. Fig.~3.5, p.~46]{Han2012}. This is not true for the phenomenon of {\it softening} (also known as ``negative hardening''), which cannot be modeled in the framework, and which requires a further nontrivial modification, see \cite{Gudoshnikov2025softening}. 

Before we derive the sweeping process from Definition \ref{def:aehp}, we will need the counterpart of space $\widehat{H}$ with the weighted inner product.
\begin{Definition} Consider the operator
\[
\widehat{{\bf C}}:\widehat{\mathcal{H}} \to \widehat{\mathcal{H}}, 
\qquad 
\widehat{{\bf C}}:\begin{pmatrix}\bm{\sigma}\\ \bm{\xi}\end{pmatrix} \mapsto \begin{pmatrix}{\bf C}\, \bm{\sigma}\\ \bm{\xi}\end{pmatrix}
\]
and its inverse
\[
\widehat{{\bf C}}^{-1}:\widehat{\mathcal{H}} \to \widehat{\mathcal{H}}, 
\qquad 
\widehat{{\bf C}}^{-1}:\begin{pmatrix}\bm{\sigma}\\ \bm{\xi}\end{pmatrix} \mapsto \begin{pmatrix}{\bf C}^{-1}\, \bm{\sigma}\\ \bm{\xi}\end{pmatrix}.
\]
Define the space $\widehat{\mathcal{H}}_{{\bf C}^{-1}}$ as
\[
\widehat{\mathcal{H}}_{{\bf C}^{-1}}= \mathcal{H}_{{\bf C}^{-1}} \times \mathcal{H}',
\]
i.e. $\widehat{\mathcal{H}}_{{\bf C}^{-1}}$ is the space  $\widehat{H}$ with the following inner product
\begin{multline}
\left\la\begin{pmatrix}\bm{\sigma}\\ \bm{\xi}\end{pmatrix}, \begin{pmatrix}\bm{\tau}\\ \bm{\zeta}\end{pmatrix}\right\ra_{\widehat{{\bf C}}^{-1}}=
\left\la\begin{pmatrix}\bm{\sigma}\\ \bm{\xi}\end{pmatrix}, \widehat{{\bf C}}^{-1}\begin{pmatrix}\bm{\tau}\\ \bm{\zeta}\end{pmatrix}\right\ra_{\widehat{\mathcal{H}}}= \la \bm{\sigma}, {\bf C}^{-1} \bm{ \tau} \ra_\mathcal{H}+\la\bm{\xi}, \bm{\zeta}\ra_{\mathcal{H}'} \\
\text{for any } \bm{\sigma}, \bm{\tau}\in \mathcal{H},\, \bm{\xi}, \bm{\zeta}\in \mathcal{H}'
\label{eq:weighted-ip-hardening} 
\end{multline}
We denote by $N^{\widehat{{\bf C}}^{-1}}$ the normal cone in the Hilbert space $\widehat{\mathcal{H}}_{{\bf C}^{-1}}$, defined by formula \eqref{eq:nc-abstract-def} with the inner product \eqref{eq:weighted-ip-hardening}. Finally, we denote by $P_{1}$ the projection from the space $\widehat{\mathcal{H}}_{{\bf C}^{-1}}$ to its subspace $\mathcal{H}_{{\bf C}^{-1}}$, i.e. the extraction of the first component:
\[
P_1:\widehat{\mathcal{H}}_{{\bf C}^{-1}} \mapsto \mathcal{H}_{{\bf C}^{-1}}, \qquad P_1:\begin{pmatrix}\bm{\sigma}\\ \bm{\xi}\end{pmatrix} \mapsto \bm{\sigma}.
\]
\end{Definition}
\begin{Theorem}
\label{th:aehp-to-sweeping-process}
Functions \eqref{eq:aehp-unknowns} solve the abstract problem of quasi-static evolution in elasticity-monotone plasticity if and only if the unknown $y$ as in \eqref{eq:yielding-to-sweeping} and the unknowns $\varepsilon \in W^{1,\infty}(I, \mathcal{H}),\, \xi \in W^{1,\infty}(I, \mathcal{H}')$ solve the differential inclusions
\begin{numcases}{}
-\frac{d}{dt}\begin{pmatrix}\bm{y}\\ \bm{\xi} \end{pmatrix}\in N^{\widehat{{\bf C}}^{-1}}_{\widehat{\mathcal{C}}(t)}\begin{pmatrix}\bm{y}\\ \bm{\xi}\end{pmatrix}, \label{eq:sp-elastoplastic-hardening}\\
\frac{d}{dt}\bm{\varepsilon} \in {\bf C}^{-1}\left(P_1\left(\left(N^{\widehat{{\bf C}}^{-1}}_{\widehat{\Sigma}-\begin{pmatrix}\bm{\widetilde{\sigma}}(t)\\0\end{pmatrix}}\begin{pmatrix}\bm{y}\\\bm{\xi}\end{pmatrix}+ \frac{d}{dt} \begin{pmatrix}\bm{y}\\ \bm{\xi}\end{pmatrix}\right)\cap \left(\mathcal{U}\times\{0\}\right)\right)+ \frac{d}{dt} \bm{\widetilde{\sigma}}(t)\right) \label{eq:inclusion2-elastoplastic-hardening}
\end{numcases}
with the initial conditions
\begin{align*}
\bm{y}(0) &= \bm{\sigma}_{\bf 0}- \bm{\widetilde{\sigma}}(0)\\
\begin{pmatrix} \bm{y}(0) \\ \bm{\xi}(0)\end{pmatrix}& \in \mathcal{C}(0), \\
\bm{\varepsilon}(0)& = \bm{\varepsilon}_{\bf 0}\in {\bf C}^{-1}(\mathcal{U}+\bm{\widetilde{\sigma}}(0)),
\end{align*}
where the moving set is
\begin{equation}
\widehat{\mathcal{C}}(t)=\left(\widehat{\Sigma}-\begin{pmatrix}\bm{\widetilde{\sigma}}(t)\\ 0\end{pmatrix}\right)\cap \left(\mathcal{V}\times \mathcal{H}'\right).
\label{eq:aehp-sp-moving-set}
\end{equation}
\end{Theorem}
\noindent {\bf Proof.} We follow the lines of the proof of Theorem \ref{th:aepp-to-sweeping-process} and by the same reasoning we obtain \eqref{eq:stress-constraint}, \eqref{eq:rate-strain-constraint} and \eqref{eq:rate-additive-decomposition} from \eqref{eq:aehp-1}--\eqref{eq:aehp-3}, \eqref{eq:aehp-5}. Apply $\widehat{\bf C}$ to \eqref{eq:aehp-5} to obtain
\begin{equation}
\widehat{\bf C} \frac{d}{dt} \begin{pmatrix}\bm{\varepsilon_{\bf p}}\\ -\bm{\xi} \end{pmatrix} \in N^{\widehat{\bf C}^{-1}}_{\widehat{\Sigma}} \begin{pmatrix}\bm{\sigma}\\ \bm{\xi}\end{pmatrix},
\label{eq:plastic-flow-rule-in-proof-hardening}
\end{equation}
i.e.
\[
\frac{d}{dt} \begin{pmatrix}{\bf C}\, \bm{\varepsilon}_{\bf p}\\ -\bm{\xi} \end{pmatrix} \in N^{\widehat{\bf C}^{-1}}_{\widehat{\Sigma}} \begin{pmatrix}\bm{\sigma}\\ \bm{\xi}\end{pmatrix}.
\]
Substitute there \eqref{eq:rate-additive-decomposition} to get
\[
-\frac{d}{dt} \begin{pmatrix}\bm{\sigma} \\ \bm{\xi} \end{pmatrix} \in N^{\widehat{\bf C}^{-1}}_{\widehat{\Sigma}} \begin{pmatrix}\bm{\sigma}\\ \bm{\xi}\end{pmatrix} - \frac{d}{dt} \begin{pmatrix}{\bf C}\, \bm{\varepsilon}\\ 0\end{pmatrix}
\]
and use \eqref{eq:rate-strain-constraint} to obtain
\[
-\frac{d}{dt} \begin{pmatrix}\bm{\sigma} \\ \bm{\xi} \end{pmatrix} \in N^{\widehat{\bf C}^{-1}}_{\widehat{\Sigma}} \begin{pmatrix}\bm{\sigma}\\ \bm{\xi}\end{pmatrix} - \frac{d}{dt}\begin{pmatrix} \bm{\widetilde{\sigma}}\\ 0\end{pmatrix} - (\mathcal{U}\times \{0\}),
\]
i.e.
\[
-\frac{d}{dt} \begin{pmatrix}\bm{\sigma}- \bm{\widetilde{\sigma}} \\ \bm{\xi} \end{pmatrix} \in N^{\widehat{\bf C}^{-1}}_{\widehat{\Sigma}} \begin{pmatrix}\bm{\sigma}\\ \bm{\xi}\end{pmatrix} + (\mathcal{U}\times \{0\}).
\]
On the other hand, from  \eqref{eq:stress-constraint} we have 
\[
\begin{pmatrix}
\bm{\sigma}-\bm{\widetilde{\sigma}}\\ \bm{\xi}
\end{pmatrix}\in \mathcal{V}\times \mathcal{H}',
\]
where the right-hand side is the orthogonal complement of $\mathcal{U}\times \{0\}$ in the space $\widehat{\mathcal{H}}_{{\bf C}^{-1}}$. Therefore,
\[
-\frac{d}{dt} \begin{pmatrix}\bm{\sigma}- \bm{\widetilde{\sigma}} \\ \bm{\xi} \end{pmatrix} \in N^{\widehat{\bf C}^{-1}}_{\widehat{\Sigma} - \begin{pmatrix}\bm{\widetilde{\sigma}}\\0\end{pmatrix}} \begin{pmatrix}\bm{\sigma}-\bm{\widetilde{\sigma}}\\ \bm{\xi}\end{pmatrix}+ N^{\widehat{\bf C}^{-1}}_{\mathcal{V}\times \mathcal{H}'} \begin{pmatrix}\bm{\sigma}-\bm{\widetilde{\sigma}}\\ \bm{\xi}\end{pmatrix}.
\]
By the subadditivity of the normal cones \eqref{eq:nc-subadditivity-general} we have
\[
-\frac{d}{dt} \begin{pmatrix}\bm{\sigma}- \bm{\widetilde{\sigma}} \\ \bm{\xi} \end{pmatrix} \in N^{\widehat{\bf C}^{-1}}_{\left(\widehat{\Sigma} - \begin{pmatrix}\bm{\widetilde{\sigma}}\\0\end{pmatrix}\right)\cap \left(\mathcal{V}\times \mathcal{H}'\right)} \begin{pmatrix}\bm{\sigma}-\bm{\widetilde{\sigma}}\\ \bm{\xi}\end{pmatrix}.
\]
Substitute \eqref{eq:yielding-to-sweeping} to obtain the desired sweeping process \eqref{eq:sp-elastoplastic-hardening}.

To derive \eqref{eq:inclusion2-elastoplastic-hardening} observe, that \eqref{eq:rate-additive-decomposition} implies 
\[
\widehat{\bf C} \frac{d}{dt} \begin{pmatrix}\bm{\varepsilon}\\0\end{pmatrix} = \frac{d}{dt} \begin{pmatrix}\bm{\sigma}\\\bm{\xi}\end{pmatrix}+\widehat{\bf C}\frac{d}{dt} \begin{pmatrix}\bm{\varepsilon}_{\bf p}\\ -\bm{\xi}\end{pmatrix}.
\]
Substitute \eqref{eq:plastic-flow-rule-in-proof-hardening} to get 
\[
\widehat{\bf C} \frac{d}{dt} \begin{pmatrix}\bm{\varepsilon}\\0\end{pmatrix} \in \frac{d}{dt} \begin{pmatrix}\bm{\sigma}\\\bm{\xi}\end{pmatrix}+N^{\widehat{\bf C}^{-1}}_{\widehat{\Sigma}} \begin{pmatrix}\bm{\sigma}\\ \bm{\xi}\end{pmatrix}
\]
i.e. 
\[
\widehat{\bf C} \frac{d}{dt} \begin{pmatrix}\bm{\varepsilon}\\0\end{pmatrix} \in N^{\widehat{\bf C}^{-1}}_{\widehat{\Sigma} - \begin{pmatrix}\bm{\widetilde{\sigma}}\\ 0\end{pmatrix}} \begin{pmatrix}\bm{\sigma}- \bm{\widetilde{\sigma}}\\ \bm{\xi}\end{pmatrix}+ \frac{d}{dt} \begin{pmatrix}\bm{\sigma}-\bm{\widetilde{\sigma}}\\\bm{\xi}\end{pmatrix}+\frac{d}{dt}\begin{pmatrix}\bm{\widetilde{\sigma}}\\0 \end{pmatrix}.
\]
But from \eqref{eq:rate-strain-constraint} we have
\[
\widehat{\bf C} \frac{d}{dt} \begin{pmatrix}\bm{\varepsilon}\\0\end{pmatrix} \in \left(\mathcal{U}\times \{0\}\right)+ \frac{d}{dt}\begin{pmatrix}\bm{\widetilde{\sigma}}\\0 \end{pmatrix},
\]
hence
\[
\widehat{\bf C} \frac{d}{dt} \begin{pmatrix}\bm{\varepsilon}\\0\end{pmatrix} \in \left( N^{\widehat{\bf C}^{-1}}_{\widetilde{\Sigma} - \begin{pmatrix}\bm{\widetilde{\sigma}}\\ 0\end{pmatrix}} \begin{pmatrix}\bm{\sigma}- \bm{\widetilde{\sigma}}\\ \bm{\xi}\end{pmatrix}+ \frac{d}{dt} \begin{pmatrix}\bm{\sigma}-\bm{\widetilde{\sigma}}\\\bm{\xi}\end{pmatrix}\right)\cap \left(\mathcal{U}\times \{0\}\right)+\frac{d}{dt}\begin{pmatrix}\bm{\widetilde{\sigma}}\\0 \end{pmatrix}.
\]
Substitute \eqref{eq:yielding-to-sweeping} and apply $\widetilde{\bf C}^{-1}$ to obtain
\[
\frac{d}{dt}\begin{pmatrix} \bm{\varepsilon}\\0 \end{pmatrix} \in \widehat{{\bf C}}^{-1}\left(\left(N^{\widehat{{\bf C}}^{-1}}_{\widehat{\Sigma}-\begin{pmatrix}\bm{\widetilde{\sigma}}\\0\end{pmatrix}}\begin{pmatrix}\bm{y}\\\bm{\xi}\end{pmatrix}+ \frac{d}{dt} \begin{pmatrix}\bm{y}\\ \bm{\xi}\end{pmatrix}\right)\cap \left(\mathcal{U}\times\{0\}\right)+ \frac{d}{dt} \begin{pmatrix}\bm{\widetilde{\sigma}}\\0 \end{pmatrix}\right),
\]
which is equivalent to the desired inclusion \eqref{eq:inclusion2-elastoplastic-hardening}.

The derivation of \eqref{eq:aehp-1}--\eqref{eq:aehp-5} from \eqref{eq:yielding-to-sweeping}, \eqref{eq:sp-elastoplastic-hardening}--\eqref{eq:aehp-sp-moving-set} is similar to the corresponding part of the proof of Theorem \ref{th:aepp-to-sweeping-process}. $\blacksquare$

\subsection{The continuum model with elasticity-hardening plasticity}
\label{ssect:aehp-rod}
Consider again the continuous rod of Section \ref{ssect:ex21-elasticity}, endowed with the elastic-hardening plastic behavior at each point $x\in \Omega$. The rod is characterized by its stiffness function $\bm{C}\in L^\infty(\Omega)$ (see Section \ref{sssect:ex21-Hooke-law}) and its yield limits $\sigma^-_x(\xi), \sigma^+_x(\xi)$, which in this paper we consider as functions
\[
(x,\xi)\in \Omega \times\mathbb{R} \longmapsto \sigma^-_x(\xi)\in \mathbb{R},
\]
\[
(x,\xi)\in \Omega \times\mathbb{R} \longmapsto \sigma^+_x(\xi) \in \mathbb{R},
\]
depending on the point $x\in \Omega$ and the value of the internal variable $\xi$ at that point.
In the Definition \ref{def:aehp} we now use the data of Section \ref{ssect:ex21-elasticity} with 
\[
\mathcal{H}'= L^2(\Omega), \qquad \widehat{H} = L^2(\Omega)\times L^2(\Omega),
\]
\begin{equation}
\widehat\Sigma = \left\{\begin{pmatrix} \bm{\sigma}\\ \bm{\xi} \end{pmatrix}\in L^2(\Omega)\times L^2(\Omega): \text{ for a.a. }x\in \Omega \text{ we have }\sigma^-_x(\xi(x))\leqslant \sigma(x)\leqslant \sigma^+_x(\xi(x))\right\}.
\label{eq:aehp-general-set-sigma}
\end{equation}
The constraint on the stress to belong to the elastic range now takes the form 
\[
\sigma(t,x) \in [\sigma^-_x(\xi(t,x)), \sigma^+_x(\xi(t,x))] \qquad  \text{for all }t\in I \text{ and a.a. }x\in \Omega,
\]
where $\xi(t,x)$ is the new internal unknown variable we introduced in \eqref{eq:aehp-unknowns}, and it allows to model the change of the elastic range during the plastic evolution. We will only consider the cases for which we can ensure that 
\begin{equation}
\sigma^-_x(\xi(t,x))< \sigma^+_x(\xi(t,x)) \qquad \text{for a.a. }x\in \Omega
\label{eq:aehp-nondegenerate-elastic-range}
\end{equation}
during the entire evolution. In particular, we assume that \eqref{eq:aehp-nondegenerate-elastic-range} holds for the initial value $\bm{\xi}(t)=\bm{\xi_0}$.

We must make several further assumptions on the data. At first, for a. a. $x\in \Omega$
\begin{align}
\text{ the map}\quad& \xi\longmapsto \sigma^-_x(\xi) \quad\text{is convex, strictly monotone and surjective onto }\mathbb{R} \label{eq:aehp-top-yield-limit-map}
\\
\text{ the map}\quad& \xi\longmapsto \sigma^+_x(\xi) \quad \text{is concave, stricitly monotone and surjective onto }\mathbb{R}. \label{eq:aehp-bottom-yield-limit-map}
\end{align}

We consider the following two simple cases, in which we can guarantee the non-degeneracy of the elastic range \eqref{eq:aehp-nondegenerate-elastic-range} and a constraint qualification of the intersection in \eqref{eq:aehp-sp-moving-set}.

\subsubsection{Expansion of the elastic range, e.g. nonlinear isotropic hardening}
\label{sssect:aehp-rod-isotropic}
\begin{figure}[H]\center
\includegraphics{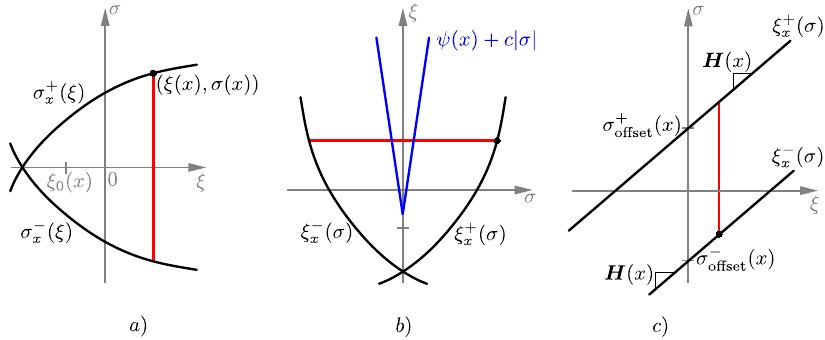}
\caption{\footnotesize Elastic range (red interval), depending on the internal variable $\xi$ for a particular $x\in \Omega$. {\it a)} Elastic range expands nonlinearly in the plastic regime (isotropic hardening in particular). {\it b)} Linear growth condition \eqref{eq:aehp-isotropic-hardening-linear-growth}. {\it c)} Elastic range translates linearly, which is called linear kinematic hardening.
} 
\label{fig:hardening-rules}
\end{figure}

As the first case, consider the situation when the maps \eqref{eq:aehp-top-yield-limit-map}--\eqref{eq:aehp-bottom-yield-limit-map} have the opposite type of monotonicity. For definiteness, we set that for a. a. $x\in \Omega$  the map \eqref{eq:aehp-top-yield-limit-map} is increasing, and \eqref{eq:aehp-bottom-yield-limit-map} is decreasing, see Fig. \ref{fig:hardening-rules} {\it a)}. Under our assumptions both maps are invertible, and we denote their respective inverses by $\xi^-_x(\sigma),\,  \xi^+_x(\sigma)$, i.e. we have the maps
\begin{align*}
\sigma\in \mathbb{R} \longmapsto \xi^-_x(\sigma) \in \mathbb{R},\\
\sigma \in \mathbb{R}  \longmapsto \xi^+_x(\sigma) \in \mathbb{R}.
\end{align*}

The maps \eqref{eq:aehp-top-yield-limit-map}--\eqref{eq:aehp-bottom-yield-limit-map} are continuous on the entire $\mathbb{R}$ as convex and concave functions with finite values \cite[Th. 10.1, p. 82]{Rockafellar1970}. Hence their inverses are also continuous. We assume that for each $\xi\in \mathbb{R}$ the maps
\[
x\longmapsto\sigma^-_x(\xi), \qquad x\longmapsto\sigma^+_x(\xi)
\]
are Lebesgue-measurable, so that 
\begin{align}
(x,\sigma)\in \Omega \times\mathbb{R} \longmapsto \xi^-_x(\sigma)\in \mathbb{R}, \label{eq:eph-bottom-yield-limit-map-inverse}
\\
(x,\sigma)\in \Omega \times\mathbb{R} \longmapsto \xi^+_x(\sigma) \in \mathbb{R}
\label{eq:eph-top-yield-limit-map-inverse}
\end{align}
 are Carath\'{e}odory maps. Finally, we must assume their uniform linear growth with respect to $\sigma$:
\begin{multline}
\text{there exist } \psi\in L^2(\Omega) \text{ and } c>0 \text{ such that for a.a. }x\in \Omega \text{ and all } \sigma\in \mathbb{R}\\
\max\left(|\xi^-_x(\sigma)|, |\xi^+_x(\sigma)|\right)\leqslant \psi(x)+c|\sigma|,
\label{eq:aehp-isotropic-hardening-linear-growth}
\end{multline}
see Fig. \ref{fig:hardening-rules} {\it b)}. Under such assumptions we can write the set $\widehat \Sigma$ in \eqref{eq:aehp-general-set-sigma} as
\begin{equation}
\widehat\Sigma = \left\{\begin{pmatrix} \bm{\sigma}\\ \bm{\xi} \end{pmatrix}\in L^2(\Omega)\times L^2(\Omega): \text{ for a.a. }x\in \Omega \text{ we have }\xi(x)\geqslant \max\left(\xi^-_x(\sigma(x)), \xi_x^+(\sigma(x))\right)\right\}
\label{eq:ehp-isortopic-set-sigma-with-inverses}
\end{equation}
and observe that $\widehat\Sigma$ is a closed, convex, nonempty set. 

We use the set $\widehat\Sigma$ to construct the sweeping process \eqref{eq:sp-elastoplastic-hardening}, which happens within the subspace $\mathcal{V}\times L^2(\Omega)$, cf.  \eqref{eq:ex21-space-V}, in the space  $L^2_{{\bf C}^{-1}}(\Omega)\times L^2(\Omega)$, which is $L^2(\Omega)\times L^2(\Omega)$ with the inner product \eqref{eq:weighted-ip-hardening}. In particular, the moving set \eqref{eq:aehp-sp-moving-set} of the sweeping process is 
\begin{multline*}
\mathcal{C}(t) =\left\{\vphantom{ L^2_{{\bf C}^{-1}}}\right.(\bm{y}, \bm{\xi}) \in L^2_{{\bf C}^{-1}}(\Omega)\times L^2(\Omega): \text{ there exists }c\in \mathbb{R} \text{ such that for a.a. }x\in \Omega \text{ we have}\\ c\equiv y(x)\quad\text{and}\quad \sigma_x^-(\xi(x)) \leqslant c+ \widetilde{\sigma}(t,x)\leqslant \sigma^+_x(\xi(x))\left.\vphantom{ L^2_{{\bf C}^{-1}}} \right\},
\end{multline*}
where $\widetilde{\sigma}$ is the corresponding stress solution for elasticity \eqref{eq:w0-const-def},\eqref{eq:ex21-linear-solution}.

The particular symmetric situation when 
\[
\sigma_x^-(\xi)=-\sigma_x^+(\xi) \qquad \text{for a.a. }x\in \Omega \text{ and all } \xi\in \mathbb{R}
\]
is known as {\it nonlinear isotropic hardening} (see e.g. \cite[Sect. 3.5, pp.~70--72]{Han2012}) for the one-dimensional rod.

\subsubsection{Linear translation of the elastic range, i.e. linear kinematic hardening}
\label{sssect:aehp-rod-kinematic}
We would like to also include the case of {\it kinematic hardening} (see e.g. \cite[Sect. 5.4.3, 5.4.4, pp.~205--240]{Lemaitre1994}), but in order to simultaneously maintain the well-posedness \eqref{eq:aehp-nondegenerate-elastic-range}, the convexity of $\widehat{\Sigma}$ and the simplicity of the model (i.e. without considering the multiyield Mroz-type models, see e.g. \cite[pp.~216--218]{Lemaitre1994}, \cite[pp.~15--18]{Krejci1996}) we will only consider the case of {\it linear kinematic hardening}, see Fig. \ref{fig:hardening-rules} {\it c)}.

In such a case we specifically consider the yield limits of the form
\[
\sigma^-_x(\xi)  = \sigma^-_{\rm offset}(x)+ \bm{H}(x)\,\xi,
\]
\[
\sigma^+_x(\xi)  = \sigma^+_{\rm offset}(x)+ \bm{H}(x)\, \xi,
\]
where $\sigma_{\rm offset}^-, \sigma_{\rm offset}^+\in L^\infty (\Omega)$ we call the {\it initial yield limits} and $\bm{H}\in L^\infty(\Omega)$ is the {\it hardening modulus} 
with a uniform lower bound, which prevents the degeneration to perfect plasticity: 
\begin{equation}
\text{there exist } \eta\in \mathbb{R} \text{ such that for a.a. } x\in \Omega \text{ we have}\qquad 0<\eta\leqslant\bm{H}(x).
\label{eq:hardening-modulus-lower-estimate}
\end{equation}

Observe from \eqref{eq:aehp-general-set-sigma} that the set $\widehat \Sigma$ is also closed, convex and nonempty in such a case. The sweeping process \eqref{eq:sp-elastoplastic-hardening} happens again within the subspace $\mathcal{V}\times L^2(\Omega)$, cf.  \eqref{eq:ex21-space-V}, in the space  $L^2_{{\bf C}^{-1}}(\Omega)\times L^2(\Omega)$, which is $L^2(\Omega)\times L^2(\Omega)$ with the inner product \eqref{eq:weighted-ip-hardening}. The moving set \eqref{eq:aehp-sp-moving-set} of the sweeping process is 
\begin{equation*}
\begin{aligned}
\mathcal{C}(t) =&\left\{\vphantom{ L^2_{{\bf C}^{-1}}}\right.(\bm{y}, \bm{\xi}) \in L^2_{{\bf C}^{-1}}(\Omega)\times L^2(\Omega): \text{ there exists }c\in \mathbb{R} \text{ such that for a.a. }x\in \Omega \text{ we have}\\
&c\equiv y(x)\quad\text{and}\quad \sigma_{\rm offset}^-(x) \leqslant c+\widetilde{\sigma}(t,x) - \bm{H}(x)\,\xi(x)\leqslant \sigma_{\rm offset}^+(x)\left.\vphantom{ L^2_{{\bf C}^{-1}}} \right\},
\end{aligned}
\end{equation*}
where $\widetilde{\sigma}$ is given by \eqref{eq:w0-const-def},\eqref{eq:ex21-linear-solution} again.

\subsection{Constraint qualifications for the continuum model with elasticity-hardening plasticity }
Let us now verify that a constraint qualification is satisfied for the continuum model with hardening. Specifically, we claim the following. 
\begin{Proposition}
For any $\bm{\widetilde{\sigma}} \in L^2_{{\bf C}^{-1}}(\Omega)$ denote
\begin{equation}
\mathcal{C}_1=\widehat{\Sigma} -\begin{pmatrix}\bm{\widetilde{\sigma}}\\ 0\end{pmatrix}, \qquad \mathcal{C}_2=\mathcal{V}\times L^2(\Omega),
\label{eq:ih-check-cq}
\end{equation}
where $\widehat{\Sigma}$ and $\mathcal{V}$ are given by \eqref{eq:aehp-general-set-sigma} and \eqref{eq:ex21-space-V} respectively. Under the assumptions of Section \ref{ssect:aehp-rod} (including either Section \ref{sssect:aehp-rod-isotropic} or \ref{sssect:aehp-rod-kinematic}) Proposition \ref{prop:constraint-qualifications} applies to $\mathcal{C}_1$ and $\mathcal{C}_2$, so that constraint qualification \ref{cq-Slater-I} fails, but constraint qualifications \ref{cq-Slater-II}--\ref{cq-GCQ} are satisfied.
\end{Proposition}
\noindent {\bf Proof.} Constraint qualification of Proposition \ref{prop:constraint-qualifications} \ref{cq-Slater-I} fails, since both $\mathcal{C}_1$ and $\mathcal{C}_2$ in \eqref{eq:ih-check-cq} have empty interiors due to the same argument as in the case of perfect plasticity, see Observation \ref{obs:ex21-Slater-cq-fails} and Fig. \ref{fig:misc-1} {\it a)}. To show that Proposition \ref{prop:constraint-qualifications} \ref{cq-Slater-II} is true we demonstrate that
\[
\mathcal{C}_1-\mathcal{C}_2=\widehat{\mathcal{H}}.
\]
Take arbitrary 
\[
\begin{pmatrix}\bm{\sigma'}\\ \bm{\xi'}\end{pmatrix}\in\widehat{\mathcal{H}} = L^2_{{\bf C}^{-1}}(\Omega)\times L^2(\Omega).
\]
We need to show that there are 
\begin{equation*}
\begin{pmatrix}\bm{\sigma}\\ \bm{\xi}\end{pmatrix}\in \widehat{\Sigma},\qquad  c\in \mathbb{R}, \qquad \bm{\zeta}\in L^{2}(\Omega)
\end{equation*}
such that 
\begin{equation*}
\bm{\sigma'} = \bm{\sigma}- \bm{\widetilde{\sigma}} - c \qquad\text{and}\qquad \bm{\xi'} = \bm{\xi}-\bm{\zeta}.
\end{equation*}
Due to the arbitrariness of $\bm{\sigma'}$ we can absorb $\bm{\widetilde{\sigma}} + c$ into it, so that it is sufficient to find
\begin{equation}
\bm{\xi},\, \bm{\zeta}\in L^{2}(\Omega)
\label{eq:in-proof-cq-hardening-1}
\end{equation}
such that 
\begin{equation}
\begin{pmatrix}\bm{\sigma}'\\ \bm{\xi}\end{pmatrix}\in \widehat{\Sigma} \qquad\text{and}\qquad \bm{\xi'} = \bm{\xi}-\bm{\zeta}.
\label{eq:in-proof-cq-hardening-2}
\end{equation}

In the case of nonlinearly expanding elastic range (Section \ref{sssect:aehp-rod-isotropic}) choose
\begin{equation}
\xi(x) =  \max\left(\xi^-_x(\sigma'(x)), \xi_x^+(\sigma'(x))\right), \qquad \zeta(x) =  \xi(x) - \xi'(x)\qquad \text{for a.a. }x\in \Omega.
\label{eq:in-proof-isotropic-hardening-choice}
\end{equation}
Notice, that since the functions \eqref{eq:eph-bottom-yield-limit-map-inverse}--\eqref{eq:eph-top-yield-limit-map-inverse} are Carath\'{e}odory and satisfy the uniform linear growth condition \eqref{eq:aehp-isotropic-hardening-linear-growth}, the corresponding Nemytskii operator maps $L^2(\Omega)\cong L^2_{{\bf C}^{-1}}(\Omega)$ into $L^2(\Omega)$ \cite[Def. 5.2, Th. 5.1, pp.~62--63]{Precup2002}, so that we indeed can fulfill \eqref{eq:in-proof-cq-hardening-1} with such choice. It is evident from \eqref{eq:ehp-isortopic-set-sigma-with-inverses} that \eqref{eq:in-proof-cq-hardening-2} is satisfied as well.

In the case of linear kinematic hardening (Section \ref{sssect:aehp-rod-kinematic}) choose
\begin{equation}
\xi(x) = \frac{\sigma'(x)-\tau(x)}{\bm{H}(x)}, \qquad \zeta(x) =  \xi(x) - \xi'(x)\qquad \text{for a.a. }x\in \Omega,
\label{eq:in-proof-kinematic-hardening-choice}
\end{equation}
with an arbitrary measurable function $\bm{\tau}$ such that
\[
\tau(x) \in \left[ \sigma^-_{\rm offset}(x), \sigma^+_{\rm offset}(x)\right].
\]
Due to the assumption \eqref{eq:hardening-modulus-lower-estimate} we have \eqref{eq:in-proof-cq-hardening-1}--\eqref{eq:in-proof-cq-hardening-2} satisfied in this case as well.

Therefore, in both cases the set $\mathcal{C}_1-\mathcal{C}_2$ is the entire space, and the interior of $\mathcal{C}_1-\mathcal{C}_2$ is also the entire space, which includes the zero element. This fulfills the condition of Proposition \ref{prop:constraint-qualifications} \ref{cq-Slater-II}, from which \ref{cq-Rockafellar} and \ref{cq-GCQ} follow as well.  $\blacksquare$

\begin{Remark}
The uniform linear growth condition \eqref{eq:aehp-isotropic-hardening-linear-growth} is sharp in the sense of the constraint qualifications of Proposition \ref{prop:constraint-qualifications} \ref{cq-Slater-II}--\ref{cq-GCQ}. Indeed, \eqref{eq:aehp-isotropic-hardening-linear-growth} is not only a sufficient, but also a necessary condition for a Nemytskii operator to map $L^2(\Omega)$ into $L^2(\Omega)$, see \cite[Th. 3.4.4, p.~407]{Gasinski2005}. If the uniform linear growth condition does not hold, function $\bm{\xi}$ in \eqref{eq:in-proof-isotropic-hardening-choice} will not be in $L^2(\Omega)$ for some $\bm{\sigma'}\in L^2(\Omega)$, and even the translation on a constant function or a multiplication by a constant (which appears in Proposition \ref{prop:constraint-qualifications} \ref{cq-GCQ}) will not change the class of integrability. By the similar argument, condition \eqref{eq:hardening-modulus-lower-estimate} is necessary to have $\bm{\xi}\in L^2(\Omega)$ in \eqref{eq:in-proof-kinematic-hardening-choice}.
\end{Remark}
\begin{Remark} The uniform linear growth condition \eqref{eq:aehp-isotropic-hardening-linear-growth} implies the estimate
\[
\max\left( \|\xi^-_{\cdot}(\sigma(\cdot))\|_{L^2},\|\xi^+_{\cdot}(\sigma(\cdot))\|_{L^2}  \right)\leqslant  \|\psi\|_{L^2}+ c\|\bm{\sigma}\|_{L^2}\qquad \text{for any } \bm{\sigma}\in L^2(\Omega),
\]
see again \cite[Th. 5.1, p.~63]{Precup2002}. In \cite[Assumption 8.4, p. 230]{Han2012} one can find another way to state the uniform linear growth condition, which is very similar to \eqref{eq:aehp-isotropic-hardening-linear-growth} due to convexity of $\widehat{\Sigma}$.
\end{Remark}

Let us now state the consequence of applicability of the constraint qualifications to the problem with hardening.

\begin{Corollary} Under the assumptions of Section \ref{ssect:aehp-rod} (including either Section \ref{sssect:aehp-rod-isotropic} or \ref{sssect:aehp-rod-kinematic}), and provided with a solution of the sweeping process \eqref{eq:sp-elastoplastic-hardening}, the right-hand side of the inclusion \eqref{eq:inclusion2-elastoplastic-hardening} is nonempty for a.a. $t\in I$.
\label{cor:hardening}
\end{Corollary}
\noindent Since we have the additivity of the normal cones in the right-hand side of \eqref{eq:sp-elastoplastic-hardening} due to Proposition \ref{prop:constraint-qualifications}, the proof of Corollary \ref {cor:hardening} is analogous to the proof of Theorem \ref{th:aepp-well-posedness-safe-load-strict} \ref{enum:th:aepp-well-posedness-nonempty-rhs}. 

We have shown that constraint qualifications of Proposition \ref{prop:constraint-qualifications} \ref{cq-Slater-II}--\ref{cq-GCQ} can, indeed, capture the difference between solvable (e.g. hardening) and unsolvable (e.g. perfect plasticity) problems in $L_2(\Omega)$, which was our goal. We leave for the future a complete well-posedness proof of the problem with hardening in the framework of Definition \ref{def:aehp} and Theorem \ref{th:aehp-to-sweeping-process}, which could be done in the spirit of Theorem \ref{th:aepp-well-posedness-safe-load-strict}. The ideas for the proof and other plans will be discussed in the next section.
\section{Conclusions}
\label{sect:conclusions}
We have formulated rigorous frameworks for linear elasticity, elasticity-perfect plasticity and elasticity-hardening plasticity (Definitions \ref{def:ae}, \ref{def:aepp}, \ref{def:aehp} respectively) in terms of abstract adjoint linear operators (Section \ref{ssect:adjoint-operatorsED}) and converted them to equivalent formulations in terms of differential inclusions (Theorems \ref{th:elasticity-solution}, \ref{th:aepp-to-sweeping-process}, \ref{th:aehp-to-sweeping-process} respectively). These differential inclusions have convenient forms of, respectively, a sweeping process and an ``open loop'' problem, i.e. a problem with the known right-hand side.   Curiously, notice that even the formulation \eqref{eq:abstract_intersection} of the elasticity problem can be seen as a sweeping process with a singleton moving set. These frameworks are ready to be used for discrete models, as well as for continuous models.

Recall, that the ``dual'' approach, which we use, gives an additional insight about the behavior of the constraint set $\mathcal{C}(t)$. We mean that, if the force load remains constant in time, then $\mathcal{C}(t)$ moves only by a translation, see Remark \ref{rem:aepp-sp-in-V} and Corollary \ref{cor:constant-force}. At least in the discrete case, one can think about the displacement boundary conditions as being ``more respectful'' to perfect plasticity, or more compatible with it. 

In this line of thinking it may be tempting to conjecture that in elasticity-perfect plasticity any ``pathological'' situations, such as a measure-valued strain, may appear only during yielding caused by a prescribed change in a force load. This is not the case, as we have demonstrated by Example 3 in Section \ref{sect:regularity-lost}, where force load was only used to form a stress profile well within the elastic range at a.a. $x\in \Omega$, and after that displacement load was applied to drive the system into the plastic regime with no strain solution in $L^2(\Omega)$. By the same mechanism, any extremum of the stress or of the yield limit distributions (e.g. caused by a defect in the material) can lead to a measure-valued strain, even if the plastic deformation is due to a change in the displacement load and the motion of $\mathcal{C}(t)$ is just a parallel translation.  

The comparison of discrete and continuous models within the same framework is our main achievement. On one hand, discrete examples are well-solvable for both strain (elongation) and stress in the case of perfect plasticity. On the other hand, a continuous example may not have a strain solution in $L^2(\Omega)$, because  the lack of additivity of normal cones prevents us from extracting the evolution of total strain $\bm{\varepsilon}$ from the evolution of stress $\bm{\sigma}$. As one would expect, the classical condition (Slater I constraint qualification) to ensure the additivity of normal cones fails hopelessly for plasticity-type unilateral constrains due to properties of $L^2(\Omega)$. But this failure also includes elastoplasticity with hardening, which is known to be a well-solvable problem.

Using later advances in the infinite-dimensional optimization, we were able to show that elastoplasticity with hardening possesses additivity of the normal cones, although it requires the use of more sophisticated constraint qualifications and the assumption of a uniform linear growth. We can classify elastoplastic models to distinguish those, for which Theorems \ref{th:aepp-to-sweeping-process} and \ref{th:aehp-to-sweeping-process} can be used to extract the evolution of strain in a certain space, see Table \ref{tab:constraint-qualifications-plasticity}.

This work is a foundation for future developments, such as:
\begin{itemize}
\item Limit analysis for cyclical loads, previously performed for discrete models  \cite{Gudoshnikov2021ESAIM}, could be done for the current abstract framework in general, although there are certain challenges to overcome. Such analysis can include the model of a three-dimensional continuous medium \cite{Gudoshnikov2025}, which fits in the current framework.
\item It would be interesting to find examples of the similar kind, which would satisfy just one or two constraint qualifications from Proposition \ref{prop:constraint-qualifications}, cf. Table \ref{tab:constraint-qualifications-plasticity}, or an even weaker constraint qualification compared to those of Proposition \ref{prop:constraint-qualifications}. Although we admit that condition \cite[Assumption 8.4, p. 230]{Han2012} is available for the models with hardening, from which Proposition \ref{prop:constraint-qualifications} cannot offer much improvement, it seems that the natural relation between constraint qualifications in infinite-dimensional optimization and elastoplastic models is mostly unexplored. We only can name \cite{Daniele2014}, where the authors apply a constraint qualification to the static elastoplastic torsion problem.
\item To keep the volume of the text at least somewhat reasonable, for the Theorem \ref{th:aehp-to-sweeping-process} we focused on the aspect of additivity of the normal cones in \eqref{eq:inclusion2-elastoplastic-hardening}. In other words, for the models with hardening we have only proven the fact, which corresponds to the part \ref{enum:th:aepp-well-posedness-nonempty-rhs} of Theorem \ref{th:aepp-well-posedness-safe-load-strict}. Strictly speaking, to claim the solvability in the models with hardening, we also have to prove the facts, similar to parts \ref{enum:th:aepp-well-posedness-sp} and \ref{enum:th:aepp-well-posedness-estimate-rhs} of Theorem \ref{th:aepp-well-posedness-safe-load-strict}, i.e. the Lipschitz-continuity of $\widehat{\mathcal{C}}(t)$ with respect to Hausdorff distance in \eqref{eq:sp-elastoplastic-hardening} and an estimate on a selection from the right-hand side of \eqref{eq:inclusion2-elastoplastic-hardening}. 
While some results for particular cases are available, e.g. in \cite{Duvaut1972} and \cite[Lemma~8.8, p.~236]{Han2012}, it would be interesting to provide an answer on the Lipschitz-continuity of the moving set with respect to the Hausdorff distance, so that such answer would be general in the same sense as the additivity of normal cones and constraint qualifications provide a general answer to the question of nonemptiness of the right-hand sides in \eqref{eq:inclusion2-elastoplastic} and \eqref{eq:inclusion2-elastoplastic-hardening}. One might think about some kind of ``differential constraint qualifications''. 

Lipschitz estimates on $\mathcal{C}(t)$ are important not only for the problems with hardening, but for elastoplasticity of $d$-dimensional continuous bodies in general when $d\geqslant 2$. Indeed, in all such problems the admissible set $\Sigma$ or $\widehat{\Sigma}$  is unbounded: recall the continuous rod with hardening \eqref{eq:aehp-general-set-sigma} and von Mises, Tresca, Mohr-Coulomb etc. yield criteria \cite[Sect.~3.3,~3.4]{Han2012} respectively. Therefore Moreau's general argument in the proof of Theorem \ref{th:aepp-well-posedness-safe-load-strict} \ref{enum:th:aepp-well-posedness-sp} and the proof of the estimate in \ref{enum:th:aepp-well-posedness-estimate-rhs} cannot be used there in the current form. Instead, we must rely on the understanding of how the intersection of the subspace $\mathcal{V}$ with a time-dependent non-polyhedral infinite-dimensional unbounded set in \eqref{eq:aepp-sp-moving-set}, \eqref{eq:aehp-sp-moving-set} affects the velocity of the solution (i.e. the stress rate).
\end{itemize}

\begin{table}[H]
\begin{center}
\renewcommand{\arraystretch}{1.2} 
\begin{tabular}{|c|c|c|c|c|}
\hline
\multirow{3}{*}{\textbf{Constraint qualification}}
&\multirow{3}{*}{$\begin{array}{c} \textbf{Discrete}\\\textbf{models}\end{array}$}&
\multicolumn{3}{c|}{\textbf{Continuous models in $L^2(\Omega)$}}
\\\cline{3-5}
 && \multirow{2}{*}{ $\begin{array}{c} \textbf{perfect}\\\textbf{plasticity}\end{array}$} & \multicolumn{2}{c|}{\textbf{hardening}}\\\cline{4-5}
&&&
$\begin{array}{c} \textbf{nonuniform}\\ \textbf{or sublinear}\\\textbf{growth of the}\\ \textbf{elastic range,}\\\textbf{i.e.} \\\textbf{\eqref{eq:aehp-isotropic-hardening-linear-growth} / \eqref{eq:hardening-modulus-lower-estimate}}\\\textbf{fails}\end{array}$
&
$\begin{array}{c}\textbf{uniform}\\ \textbf{linear}\\ \textbf{growth of the} \\ \textbf{elastic range}\\ \textbf{i.e.} \\ \textbf{\eqref{eq:aehp-isotropic-hardening-linear-growth} / \eqref{eq:hardening-modulus-lower-estimate}}\\\textbf{holds}\end{array}$
\\
\hline
Slater I, Prop. \ref{prop:constraint-qualifications} \ref{cq-Slater-I} & {\color{ForestGreen}YES} &{\color{red}NO} &{\color{red}NO} &  {\color{red}NO}
\\
\hline
Slater II, Prop. \ref{prop:constraint-qualifications} \ref{cq-Slater-II} & {\color{ForestGreen}YES} &{\color{red}NO} & {\color{red}NO} &  {\color{ForestGreen}YES}
\\
\hline
Rockafellar, Prop. \ref{prop:constraint-qualifications} \ref{cq-Rockafellar} & {\color{ForestGreen}YES} &{\color{red}NO}&{\color{red}NO} &  {\color{ForestGreen}YES} \\
\hline
Attouch--Brezis, Prop. \ref{prop:constraint-qualifications} \ref{cq-GCQ}& {\color{ForestGreen}YES} &{\color{red}NO}&{\color{red}NO} &  {\color{ForestGreen}YES} \\
\hline
\end{tabular}
\end{center}
\caption{\footnotesize Constraint qualifications for different types of elastoplastic models.}  
\label{tab:constraint-qualifications-plasticity}
\end{table}

\appendix
\section{Some specific preliminaries}
\label{sect:appendix}
\subsection{Sobolev functions with values in Banach spaces}
\label{ssect:prelim-bochner-sobolev}
Let $I=[0,T], T>0$ be the time-domain of an evolution problem. For Bochner-Sobolev spaces $W^{1,p}(I,X)$ of functions with values in a Banach space $X$ we refer to \cite[Section 8.5, pp.~229--236]{Leoni2017}, \cite[Section II.5, pp.~92--110]{Boyer2012} and \cite[Section 2.4, pp.~37--49]{Migorsky2013}. We would like to remind the reader that the characterization of a function as  Bochner-Sobolev is equivalent to classical differentiability a.e. on $I$. Specifically, we mean the following facts.
\begin{Proposition}{\rm (\cite[Prop. 2.50, p.~47]{Migorsky2013}, cf. \cite[Th. 8.55, p.~230]{Leoni2017})} Let $X$ be a reflexive Banach space and $u:I\to X$ be an absolutely continuous function. Then $u$ is differentiable a.e. on $I$ in the classical sense with $u'\in L^{1}(I, X)$ and 
\[
u(t)=u(0) +\int\limits_{0}^{t} u'(s) ds \qquad \text{for all }t\in I.
\]
\end{Proposition}
\begin{Proposition}
\label{prop:prelim-classical-derivatives-ae}
{\rm (\cite[Prop. 2.51, p.~47]{Migorsky2013}, cf. \cite[Th. 8.57, p.~232]{Leoni2017})}
Let $1\leqslant p\leqslant \infty$, $k\in \mathbb{N}$ and $u\in L^{p}(I, X)$. Then the following properties are equivalent:
\begin{enumerate}[{\it i)}]
\item $u\in W^{k,p}(I, X)$,
\item there exists an absolutely continuous function $\widetilde{u}$, such that $u=\widetilde{u}$ a.e. on $I$ and for all $m\in \overline{0, k}$ the classical derivatives $\widetilde{u}^{(m)}$ of $\widetilde{u}$ exist a.e. on $I$ with $\widetilde{u}^{(m)} \in L^{p}(I, X)$.
\end{enumerate}
\end{Proposition}

\subsection{Moore-Penrose pseudoinverse matrix}
In this text we use (real) Moore-Penrose pseudoinverse matrix as it is an important practical tool to solve linear algebraic equations, and it is also available in many numerical packages. 
Here we remind the reader the definition of the Moore-Penrose pseudoinverse and some of its basic properties.
\begin{Proposition}{\bf \cite[p. 9]{Bapat2010}}
Let $A$ be an $m \times n$-matrix. Then there exists a unique $n \times m$-matrix $A^+$, called {\it Moore-Penrose pseudoinverse of $A$}, such that all of the following hold:
\begin{align}
AA^+A&=A, \label{eq:MP1}\\
A^+AA^+&=A^+, \label{eq:MP2}\\
(AA^+)^\top&=AA^+,\label{eq:MP3}\\
(A^+A)^\top&=A^+A.\label{eq:MP4}
\end{align}
\label{prop:MP1}
\end{Proposition}
\begin{Proposition}{\bf \cite[Def. 1.1.2]{Campbell2008}}
A matrix $A^+$ is a Moore-Penrose pseudoinverse of $A$ if and only if $AA^+$ and $A^+A$ are orthogonal projection matrices onto, respectively, ${\rm Im}\,A$ and ${\rm Im}\,A^\top$.
\label{prop:MP2}
\end{Proposition}
\noindent The following proposition can be verified by direct substitution into \eqref{eq:MP1}-\eqref{eq:MP4}:
\begin{Proposition}
If the columns of $A$ are linearly independent, then $A^+= (A^\top A)^{-1}A^\top$ and $A^+$ is a {\it left inverse} of $A$, i.e.
\begin{equation*}
A^+A=I_{n\times n}.
\end{equation*}
Similarly, if the rows of $A$ are linearly independent, then $A^+= A^\top(AA^\top)^{-1}$ and $A^+$ is a {\it right inverse} of $A$, i.e.
\begin{equation}
AA^+=I_{m\times m}.
\label{eq:mp-right-inverse}
\end{equation}
\label{prop:MP3}
\end{Proposition}

\begin{Corollary}
Let $A$ be an $m\times n$ matrix with linearly independent rows and $b\in \mathbb{R}^m$. Then for $x\in \mathbb{R}^n$
\begin{equation}
Ax = b \qquad \Longleftrightarrow \qquad x\in {\rm Ker}\, A + A^+ b,
\label{eq:mp-solving-linear-systems}
\end{equation}
where 
\[
A^+=A^\top(AA^\top)^{-1}.
\]
\label{cor:mp-solving-linear-systems}
\end{Corollary}
\noindent {\bf Proof.} Assume that the left part of \eqref{eq:mp-solving-linear-systems} holds. Apply $A^+$ to it to get
\begin{equation}
A^+ A x=A^+b.
\label{eq:mp-in-proof1}
\end{equation}
From Proposition \ref{prop:MP2} we know that $A^+A$ is the matrix of orthogonal projection onto ${\rm Im}\, A^\top$, therefore there is $y\in\left( {\rm Im}\, A^\top\right)^\perp = {\rm Ker}\, A$, such that
\[
x=A^+Ax + y.
\]
Add this to \eqref{eq:mp-in-proof1} to obtain the right part of \eqref{eq:mp-solving-linear-systems}, as claimed.

Vice versa, assume that the right part of \eqref{eq:mp-solving-linear-systems} holds. Apply $A$ to it to get
\[
Ax=AA^+b.
\]
Use \eqref{eq:mp-right-inverse} to derive the right part of \eqref{eq:mp-solving-linear-systems}. $\blacksquare$

\section*{Acknowledgements} The author thanks Pavel Krej{\v{c}}{\'{i}}, Giselle Antunes Monteiro and Michal K{\v{r}}{\'{i}}{\v{z}}ek from the Institute of Mathematics
CAS for helpful scientific discussions. The author thanks A. Phung for the copy of the book \cite{Aubin2009}, which came in useful for the current work. The author honors the memory of Prof. Zalman Balanov, from whom the author first learned about Nemytskii operator. The author also thanks the anonymous referees for their effort of carefully reading the manuscript and for their feedback, which helped to improve it.
\subsection*{Funding}
\noindent This research is supported by the GA{\v{C}}R project GA24-10586S and the Czech Academy of Sciences (RVO: 67985840).

\printbibliography
\end{document}